\documentclass[12pt]{amsart}
\usepackage{amsthm}
\usepackage{amstext}
\usepackage{amssymb}
\usepackage{mathrsfs}
\usepackage{leftidx}
\usepackage[english]{babel}
\usepackage{enumerate}
\usepackage{xfrac}
\usepackage{mathtools}
\usepackage{xcolor}
\usepackage[all]{xy}
\usepackage{tikz}
\usepackage{subcaption}

\usepackage[T1]{fontenc}

\usepackage{mathpazo}

\usepackage{hyperref}

\allowdisplaybreaks

\theoremstyle{plain}
\newtheorem{theorem}{Theorem}[section]
\newtheorem*{theorem*}{Theorem}
\newtheorem{proposition}[theorem]{Proposition}
\newtheorem{corollary}[theorem]{Corollary}
\newtheorem{lemma}[theorem]{Lemma}

\newtheorem{thmALPH}{Theorem}[section]
\newtheorem{corALPH}[thmALPH]{Corollary}

\theoremstyle{definition}
\newtheorem{definition}[theorem]{Definition}

\newtheorem{example}[theorem]{Example}
\newtheorem{remark}[theorem]{Remark}

\theoremstyle{remark}

\numberwithin{equation}{section}

\newcommand{\N}{\mathbb N}

\newcommand{\R}{\mathbb R}

\newcommand{\IR}{\mathrm{R}}

\newcommand{\dom}{\mathsf{Dom}}

\newcommand{\dd}{\mathrm{d}}
\newcommand{\dist}{\mathbf{d}}
\newcommand{\signdist}{\mathbf{sd}}
\newcommand{\ssubset}{\subset\joinrel\subset}

\newcommand{\ILL}{\mathscr{L}}
\newcommand{\A}{\mathsf{A}}

\DeclareMathOperator{\seccurv}{sec}
\DeclareMathOperator{\secfund}{I\mkern-1mu I}
\DeclareMathOperator{\scalsecfund}{h}

\DeclareMathOperator{\grad}{grad}
\DeclareMathOperator{\Hess}{Hess}

\newcommand{\e}{\mathrm{e}}

\newcommand{\sn}{\mathrm{sn}}
\newcommand{\ct}{\mathrm{ct}}
\DeclareMathOperator{\inj}{inj}
\DeclareMathOperator{\ninj}{ninj}
\DeclareMathOperator{\conv}{conv}

\newcommand{\Bsf}{\mathsf{B}}
\newcommand{\Ssf}{\mathsf{S}}
\newcommand{\Csf}{\mathsf{C}}

\newcommand{\Scal}{\mathcal S}
\newcommand{\Fcal}{\mathcal F}

\newcommand{\Kcal}{\mathcal K}

\newcommand{\Fct}{\mathsf{Fct}}

\DeclareMathOperator{\arccot}{arccot}

\DeclareMathOperator{\compR}{\mathsf{R}}

\title[the normal injectivity radius of a hypersurface]{Lower bounds on the normal injectivity radius of hypersurfaces and bounded geometries on manifolds with boundary}
\author{Sebastian Boldt}
\address{Fakultät f\"ur Mathematik\\Technische Universit\"at Chemnitz\\Germany}
\email{sebastian.boldt@mathematik.tu-chemnitz.de}
\author{Batu G\"uneysu}
\address{Fakultät f\"ur Mathematik\\Technische Universit\"at Chemnitz\\Germany}
\email{batu.gueneysu@mathematik.tu-chemnitz.de }
\author{Stefano Pigola }
\address{Dipartimento di Matematica\\Universita degli studi di Milano Bicocca\\Milano\\Italy}
\email{stefano.pigola@unimib.it}

\begin{document}

\begin{abstract} 
    We prove for the first time a pointwise lower estimate of the normal injectivity radius of an embedded hypersurface in an arbitrary Riemannian manifold. Main applications include:
    \begin{itemize}
        \item a pointwise lower estimate of the graphing radius of a properly embedded hypersurface,
        \item the construction of metrics of bounded geometry on arbitrary manifolds with boundary,
        \item the equivalence of the classical (topological) notion of orientation with that of the geometric notion (in the sense of metric measure spaces) on arbitrary Riemannian manifolds with boundary.
    \end{itemize}
In addition, we prove that every manifold with boundary admits a metric with bounded geometry such that boundary becomes convex. This result strengthens the justification of a recent notion of \emph{orientation} on finite dimensional RCD spaces. 
\end{abstract}


\maketitle
\tableofcontents

\section{Introduction}

The normal injectivity radius is an ubiquitous object in the theory of submanifolds. It is used to show e.g.\ the existence of normal neighbourhoods of submanifolds \cite[Theorem~5.25]{lee}, which leads to all kinds of fundamental results, such as the existence of Fermi-coordinates \cite[Proposition~5.26]{lee}, or the explicit description of distance functions to submanifolds \cite[Proposition~6.37]{lee}. In typical applications in geometric analysis, such as e.g. the regularity of geometric flows and of functional minimizers, or spectral theory  \cite{bousquet, white, mramor, post,Amann,dewitt}, the normal injectivity radius is usually assumed to be uniformly bounded from below, noting that this assumption is trivially implied by the compactness of the submanifold. Despite this prominent role, no pointwise lower estimate of the normal injectivity radius in terms of geometric data seems to exist. It is the aim of this work to fill in this gap in the case of hypersurfaces.\vspace{2mm}

In order to formulate our main results, we have to introduce some notation first: we let $\Sigma$ be an embedded hypersurface in an arbitrary Riemannian manifold $(M^n,g)$. With $\exp:TM\supseteq \mathcal D\to M$ the Riemannian exponential map associated with $g$, where the domain $\mathcal D$ is an open neighbourhood of the zero section of $TM$ which in general does not coincide with $TM$ since we do not assume our manifold $(M,g)$ to be complete, for $x\in \Sigma$ and $\delta>0$ we denote
\[
    U_\delta(x) = \left\{ (y,v)\in N\Sigma\,|\,\dist(x,y)< \delta, |v|<\delta  \right\}\,,
\]
the (two-sided) cylinder of length $2\delta$ over $\Bsf_\delta(x)\cap\Sigma$ in the normal bundle $N\Sigma$ of $\Sigma$. Our chosen object of study, the \emph{normal (injectivity) radius of $\Sigma$}, is the quantitiy $\ninj^\Sigma:\Sigma\to(0,\infty]$ given by
\begin{equation*}
        \ninj^\Sigma(x) = \sup\{\delta > 0\,|\, U_\delta(x)\subseteq\mathcal D\text{ and } \exp|U_\delta(x) \text{ is a diffeomorphism onto} \}\,.
\end{equation*}
Note that $\ninj$ is indeed everywhere positive since the differential of $\exp$, viewed as a map $\mathcal D\cap N\Sigma\to M$, at every $0_x\in N_x\Sigma$ is invertible and $\Sigma$ is embedded. By considering geodesic spheres in constant curvature spaces it is immediately clear that any pointwise lower bound on $\ninj^\Sigma$ must include the curvature of $M$ and the one of $\Sigma$. These two pieces of data are, however, not enough. For one, the topology of $M$ plays a role as one can see by considering e.g.\ geodesic spheres in a flat cylinder. When one studies the injectivity radius of a point, the classic Klingenberg lemma covers this phenomenon by assuming that the length of any self-intersecting geodesic in a given compact ball has a positive lower bound \cite[Lemma~3.1]{nardmann}.

\begin{figure}
    \begin{tikzpicture}
		\node at (3,2.5) {$M$};
		\draw[smooth,tension=0.7] plot coordinates {(3,1) (5,1.3) (7,1.6) (9,2) (11,2.5) (13,2) (13,0) (11,-.5) (9,0) (7,.5) (5,.7) (3,.8)} node at (7, 2) {$\Sigma$};
		\begin{scope}
			\clip (5,.7) circle (.7);
			\draw[thick, red, smooth, tension=1] plot coordinates {(3,1) (5,1.3) (7,1.6) (9,2) (11,2.5) (13,2) (13,0) (11,-.5) (9,0) (7,.5) (5,.7) (3,.8)} node at (7, 2.5) {$\Sigma$};
		\end{scope}
                \node at (6.1,0) {$\Bsf_\delta(x)$};
			\draw[thin, dashed] (5,.7) circle (.7);
			\filldraw[black] (5,.7) circle (1pt) node[below] {$x$};
	\end{tikzpicture}
 \caption{A connected hypersurface with arbitrarily small normal radius around flat points}
 \label{fig:global-geometry-Sigma}
\end{figure}
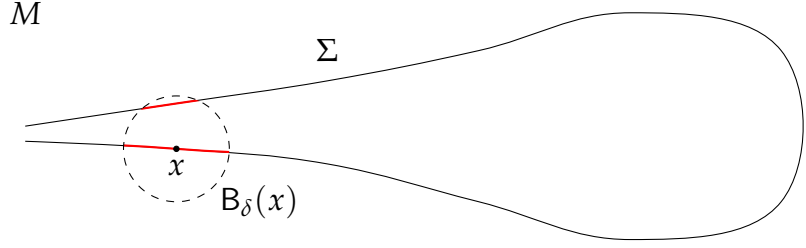

When dealing with a hypersurface $\Sigma$, its global geometry comes into play as well. This can be seen in Figure~\ref{fig:global-geometry-Sigma}. Here, $M$ is flat, has no topology and $\Sigma$ is connected. Two flat pieces of $\Sigma$ can be placed as close to each other as one likes, all the while the curvature of $\Sigma$ can be made arbitrarily small. This leads us to the notion of a $\Sigma$-\emph{slice ball}, a ball $\Bsf_r(x)$ about $x\in\Sigma$ such that $\Sigma_r(x)=\Sigma\cap\Bsf_r(x)$ is connected and separates $\Bsf_r(x)$ into exactly two components. Moreover, we call $\Bsf_r(x)$ a \emph{regular $\Sigma$-slice ball}, if $\Bsf_s(x)$ is a $\Sigma$-slice ball for all $s\leq r$ (see Definition~\ref{def:local-0}). To state our main result, we denote by $\conv(x)$ the convexity radius of $(M,g)$ at $x$ and by $\compR^c_\lambda$ the radius of a geodesic sphere with curvature $\lambda$ in a model space of constant curvature $c\in\R$. Note that this latter quantity can be explicitely calculated. Our main result reads as follows.\\

\begin{thmALPH}[Klingenberg's lemma for embedded hypersurfaces, Theorem~\ref{thm:Klingenberg-hypersurface}]\label{thm:intro-Klingenberg}
	Let $(M,g)$ be a Riemannian manifold, $\Sigma\subseteq M$ an embedded hypersurface and $x\in\Sigma$. Let $c,\lambda \geq 0$ and $s > 0$ such that
	\begin{itemize}
		\item $\Bsf_s(x)$ is a precompact, regular $\Sigma$-slice ball,
		\item $| \seccurv(\sigma) | \leq c $ for all $\sigma\in \mathrm{Gr}_2(T_yM), y\in \Bsf_s(x)$,
		\item $|\secfund^\Sigma| \leq \lambda$ on $\Sigma_s(x)$.
	\end{itemize}
	Then
	\begin{equation*}
		\ninj^{\Sigma}(x)\geq \frac{1}{6}\min\left\{s, \conv(x), \tfrac 65\compR_\lambda^c \right\}\,.
	\end{equation*}
\end{thmALPH}

The main part of the proof of Theorem~\ref{thm:intro-Klingenberg} deals with the injectivity of $\exp$ on $U_r(x)$, where $r = (1/6)\min\left\{s, \conv(x), (6/5)\compR_\lambda^c \right\}$, as the regularity of the differential of $\exp$ on $U_r(x)$ can be easily deduced from standard comparison theory estimates. We state the injectivity as a \emph{local interior Blaschke ball rolling result}.\\

\begin{thmALPH}[Local two-sided inner Blaschke ball rolling theorem, Theorem~\ref{thm:local-Blaschke}]\label{thm:intro-local-Blaschke}
    Let $(M,g)$ be a Riemannian manifold and $\Sigma\subseteq M$ an embedded hypersurface. Let $\Bsf_r(x)$ be a precompact, strongly convex regular $\Sigma$-slice ball. Let $c,\lambda \geq 0$ be real numbers such that
    \begin{itemize}
        \item $|\seccurv(\sigma)| \leq c$ for all $\sigma\in \mathrm{Gr}_2(T_yM), y\in \Bsf_r(x)$,
        \item $|\secfund^\Sigma|\leq \lambda$ on $\Sigma_r(x)$,
        \item $\frac 56 \cdot r \leq \compR_\lambda^c$.
    \end{itemize}
    Then any ball of radius $<r/3$ can roll freely on either side of $\Sigma_{r/3}(x)$ inside $\Bsf_{r/3}(x)$, i.e.\ if $\rho<r/3$, $o = \exp_q(v)\in \Bsf_{r/3}(x)$ with $q\in\Sigma_{r/3}$ and $v\in N_q\Sigma$ with $|v| = \rho$, then $\overline\Bsf_\rho(o)\cap \Sigma_{r/3} = \{q\}$.
\end{thmALPH}


Blascke's original result concerned strictly convex domains in $\R^2$:

\begin{theorem*}[Blaschke Inclusion Theorem \cite{Blaschke}]
    Let $D_1,D_2\subseteq \R^2$ be two strictly convex compact domains with smooth boundaries. Assume that $D_1$ and $D_2$ touch at $p\in\partial D_1\cap\partial D_2$ and that the curvatures of their boundaries, as functions of the inward-pointing normal $\nu\in \Ssf^1$, satisfy
    \[
            \scalsecfund^{D_1}(\nu)\leq \scalsecfund^{D_2}(\nu)\,.
    \]
    Then
    \[
        D_2\subseteq D_1\,.
    \]
\end{theorem*}

One then obtains outer respectively inner Blaschke ball rolling theorems by letting either $D_1$ or $D_2$ be a ball and says that the domain can roll freely inside the ball respectively the ball inside the domain. In this context, the corresponding rolling radius is the smallest respectively largest radius of a ball for which the above holds.

A long list of authors have worked on generalising the outer ball rolling theorem and the inclusion theorem. Instead of reproducing the account of that development here, we instead refer to the introduction of \cite{Dra24}, which represents the state of the art regarding the outer Blaschke ball rolling theorem in Riemannian manifolds.

Regarding the inner ball rolling theorem, the only work we are aware of is \cite{How99}. Here, complete Riemannian manifolds with a lower curvature bound, a compact boundary and various curvature assumptions on the boundary are considered. The rolling theorem with respect to the boundary is then stated as an equality of the focal radius of the boundary and the distance of boundary to its cut locus.

Apart from the local nature of Theorem~\ref{thm:intro-local-Blaschke}, it differs from all previous rolling-type results in that the hypersurface is not assumed to satisfy any type of convexity or positive curvature assumption.

The proof of Theorem~\ref{thm:intro-local-Blaschke} is based on a new \emph{local two-sided radial angle comparison theorem}, Theorem~\ref{thm:rad-angle-comp}. The radial angle comparison technique originated in the work of Borisenko with the aim to study strictly convex hypersurfaces with bounds on their principal curvatures \cite{Bor02, BD13, BD15, BG01, BM02}. Very recently, Drach has proved \cite{Dra24} a sharp version of the upper radial angle comparison theorem for $\lambda$-convex hypersurfaces in complete Riemannian manifolds to obtain an \emph{outer Blaschke ball rolling theorem} for such hypersurfaces. We partly review his work in Section~\ref{sec:intersection-hypersurfaces-convex-bodies}.

Our contribution to this technique is that we obtain two-sided bounds on the radial angle function and that our hypersurface $\Sigma$ need not be convex. The latter fact forces us to localise (to $\Sigma$-slice balls), which we need to anyways since our ambient manifold $M$ is neither assumed complete nor to possess any bounds on its curvature or its injectivity radius.\\

Given Theorem~\ref{thm:intro-Klingenberg}, we are interested in the construction of (regular) slice balls. This is not difficult: assuming for simplicity that $\Sigma$ is properly embedded, any ball $\Bsf_r(x)$ with $r\leq \inj(x)$ which is centred at a point $x\in\Sigma$ and which is contained in a tubular neighbourhood $U$ of $\Sigma$ is a regular slice ball. Of course, if one wishes to apply Theorem~\ref{thm:intro-Klingenberg} in this way, this construction is circular as one would need an estimate of the size of $U$ around $x$, at which point one can just use $U$ itself.

Hence, our point of view is as follows. We regard $g$ as a background metric and wish to construct slice balls for a different Riemannian metric $h$. In that case, a ball $\Bsf_r^h(x)$ with respect to $h$ simply becomes a domain $D$ of $M$ which intersects $\Sigma$. To start, we assume that $D$ is simply connected. One now needs to show that $\Sigma\cap D$ is connected and we try to enforce this situation by putting $D$ into a tubular neighbourhood $U$ of $\Sigma$. Unfortunately, without any further assumptions on $D$, this is far from the case even when $\Sigma$ is a hyperplane in $\R^n$, since a sufficiently bumpy domain $D$ will intersect $\Sigma$ in as many components as one wishes. A first assumption might be to require $D$ to be convex. However, even this assumption turns out to be too weak. Consider a sausage body $D$ in $\R^n$, i.e.\ the convex hull of two disjoint balls of radius, say, one. A hypersurface $\Sigma$ which is a plane except for an arbitrarily small dent can be placed such that it intersects $D$ in two components while $D$ is contained in a tubular neighbourhood of $\Sigma$. We thus turn to the class of \emph{$\lambda$-convex domains}, i.e. domains $D$ for which the scalar second fundamental form $\scalsecfund^D$ of their smooth boundary satisfies $\scalsecfund^D\geq \lambda$, where $\lambda$ is a positive real number, and prove the following:

\begin{thmALPH}[Theorem~\ref{thm:normal-radius-est-sigma-intersection-convex-domain}]\label{thm:intro-normal-radius-est-sigma-intersection-convex-domain}
    Let $\Sigma$ be a properly embedded hypersurface in $(M,g)$, $x\in M$, $r \leq \conv(x)$ and $\overline D\subseteq \overline\Bsf_{r/2}(x)$ be a $\lambda$-convex domain. Assume that there exists $c\in\R$ such that
   \begin{equation}\label{eqn:intro-lower-curvature}
        \phantom{\text{for all} \sigma\in \mathrm{Gr}_2(T_yM), y\in \Bsf_r(x)} c\leq\sec(\sigma)\quad\text{ for all }\quad \sigma\in \mathrm{Gr}_2(T_yM), y\in \Bsf_r(x)
    \end{equation}
    and
    \begin{equation}\label{eqn:intro-lower-radius}
        \compR^c_\lambda \leq r/4 
    \end{equation}    
    hold. Then, if $\Sigma\cap D\neq \emptyset$ and  $\Sigma\cap D$ has two components, there exists a geodesic sphere $\Ssf_\rho(p)\subseteq\Bsf_r(x)$ with $\rho\leq \compR^c_\lambda$ touching $\Sigma$ at two distinct points in $\Sigma\cap\Bsf_r(x)$.
\end{thmALPH}

Noting that a geodesic sphere which touches $\Sigma$ in two distinct points implies in particular the existence of two distinct normal geodesics to $\Sigma$ which meet at a point, we obtain the following corollary by putting the situation of Theorem~\ref{thm:intro-normal-radius-est-sigma-intersection-convex-domain} into a tubular neighbourhood of $\Sigma$. 

\begin{corALPH}[Corollary~\ref{cor:lambda-convex-body-in-Sigma-tubular-neighbourhood}]\label{cor:intro-lambda-convex-body-in-Sigma-tubular-neighbourhood}
    Let $x\in\Sigma$ and $\delta>0$ such that $\exp|U_\delta(x)$ is a diffeomorphism onto. Let $r\leq \min\{\conv(x),\delta/2\}$ and $\overline D\subseteq \overline\Bsf_{r/2}(x)$ be a $\lambda$-convex domain. Assume there exists $c\in\R$ such that \eqref{eqn:intro-lower-curvature} and \eqref{eqn:intro-lower-radius} hold. Then, $D\cap\Sigma$ is connected. In particular, $D$ is a $\Sigma$-slice domain.
\end{corALPH}

In case $D=\Bsf^h_r(x)$, the $g$-$\lambda$-convexity of $\partial D$ can be checked by choosing $r$ such that $\partial D$ is $h$-$(\lambda+\varepsilon)$-convex and a suitable local bound on the $\Csf^1$-distance of $g$ to $h$. In case $h=\e^{2u}g$ is conformal to $g$, we follow this strategy in Theorem~\ref{thm:tubular-neighbourhood-qfz}.\\

Our first application of Theorem~\ref{thm:intro-Klingenberg} is to the graphing radius of a properly embedded hypersurface $\Sigma$. As is well-known, a hypersurface in $\R^n$ is locally the graph of a function. Here, given a point $x\in\Sigma$ one has $T_x\Sigma = \R^{n-1}\subseteq \R^n$ in a canonical way and can therefore easily show the existence of an $r>0$ and a $\Csf^\infty$-function $f:T_x\Sigma\supseteq \Bsf_r(0)\to \R$ such that near $x$ one has $\Sigma = \mathrm{Graph}(f)$. The largest such $r$ is called the \emph{graphing radius} of $\Sigma$ at $x$. A lower bound on the graphing radius can be given solely in terms of bounds on the curvature of $\Sigma$. 

When we move from $\R^n$ to an arbitrary Riemannian manifold $(M,g)$, there is no canonical choice of ``base space'' on which a potential graphing function might live. Given $x\in\Sigma$, we choose as a base the embedded disc $\Gamma^{x,r} = \exp_x(T_x\Sigma\cap \Bsf_r(0))$ with $r \leq \inj(x)$. Note that $\Bsf_r(x)$ is by definition a regular $\Gamma^{x,r}$-slice ball, so that Theorem~\ref{thm:intro-Klingenberg} establishes a coordinate system over $\Gamma^{x,s}$ through parallel hypersurfaces for a suitable $s\leq r$. With this choice we define the graphing radius of $\Sigma$ at $x$ as the quantity
    \begin{align*}
        \mathrm{graph}^\Sigma(x)\coloneqq \sup \>&\big\{r\in (0,\inj^M_g(x)]\,\big|\,r\leq \ninj^{\Gamma^{x,\inj(x)}}(x) \text{ and there exists } f\in \Csf^\infty(\Gamma^{x,r})\text{ s.t. }\\\nonumber
        & -r<f<r \text{ and }  (\Sigma\cap\Bsf_r(x))^0 = \Gamma(f) \big\}\,,
    \end{align*}
    where $(\Sigma\cap\Bsf_r(x))^0$ denotes the connected component of $\Sigma\cap\Bsf_r(x)$ containing $x$ and $\Gamma(f) \coloneqq \left\{\exp_y(f(y)\cdot\nu_y)\,|\,y\in \Gamma^{x,r,T_xM}\right\}$ is the graph of $f$ with $\nu$ a unit normal to $\Gamma^{x,r}$.

\begin{thmALPH}[Theorem~\ref{thm:graphing-radius-estimate}]
    Let $\Sigma$ be a properly embedded hypersurface in a Riemannian manifold $(M,g)$ and $x\in\Sigma$. Assume there exist $r>0$ and $c,\lambda,\Lambda\geq 0$ such that
    \begin{itemize}
        \item The ball $\Bsf_r(x)$ is precompact,
        \item $|\seccurv(\sigma)|\leq c$ for all $\sigma\in\mathrm{Gr}_2(T_yM)$, $y\in\Bsf_r(x)$,
        \item $|\nabla \mathrm{\mathrm{Riem}}|\leq \Lambda$ on $\Bsf_r(x)$,
        \item $|\secfund^\Sigma|\leq\lambda$ on $\Bsf_r(x)$.
    \end{itemize}
    Then, there exists an explicitely computable constant $C = C(r,c,\lambda,\Lambda,\inj(x),\conv(x))>0$ such that
    \[
        \mathrm{graph}^\Sigma(x)\geq C\,.
    \]
\end{thmALPH}

We have recorded the proof of Theorem~E in such a way that anyone who wishes to compute the constant $C$ given concrete bounds as above can follow a concise algorithm to do so.\\

The second application of our main result concerns the construction of Riemannian metrics of bounded geometry (cf. Definition \ref{boundedgeom}) on manifolds with boundary. Bounded geometries on noncompact manifolds with boundary play a crucial role in geometric analysis, for example in the context of vanishing results for the $L^2$-cohomology, and the Hodge-de Rham-theorem for $L^2$-cohomology \cite{Schi98}, self-adjointness problems for the Laplace-Beltrami operator, and global elliptic estimates \cite{grosse}. A view towards the ultimate justification of the above cited analytic results raises the following question: \emph{does every noncompact manifold with boundary admit a Riemannian metric of bounded geometry?} We give an affirmative answer to this question in this paper. In fact, we prove the following much stronger result:

\begin{thmALPH}[Corollary \ref{cor:existence-metric-bounded-geometry-boundary}]\label{exist} Let $M$ be a manifold with boundary and $g_0$ a Riemannian metric on $M$. Then there exists a conformal factor $u\in\Csf^\infty(M)$ such that $(M,e^{2u}g_0)$ is of bounded geometry.
\end{thmALPH}

Theorem \ref{exist} connects to our above results as follows: we realize $(M,g_0)$ isometrically as a domain in a Riemannian manifold $(N,h)$ without boundary, using \cite{stefano}. The boundary $\Sigma=\partial M$ is now a properly embedded hypersurface in $N$, and the problem consists in finding a conformal factor $u\in\Csf^\infty(N)$ such that (i) $(N,\e^{2u}h)$ is of bounded geometry, (ii) $(\Sigma,h_{|\Sigma})$ is of bounded geometry, and such that (iii) $\Sigma$ possesses a uniform $h$-tubular neighbourhood in $N$. 

While (i) was first accomplished by \cite{greene}, property (iii) is the main challenge in proving Theorem~\ref{exist} above, as it corresponds to the uniform collar of the boundary in the definition of bounded geometry on manifolds with boundary. Here, we establish a theory of flatzoomers \cite{nardmann} for hypersurfaces, in which Theorem~\ref{thm:intro-Klingenberg} and Corollary~\ref{cor:intro-lambda-convex-body-in-Sigma-tubular-neighbourhood} yield the essential tools for the construction of the conformal factor $u$ such that (iii) holds. We note in passing that the constructions of \cite{greene, nardmann} do not apply directly to manifolds with boundary. There, a manifold is exhausted by smooth compact sets and the conformal factor is essentially constructed over uniform tubular neighbourhoods of the boundaries of the exhausting sets. In contrast, any sequence of compact sets exhausting a manifold with boundary will eventually hit the boundary, at which point there simply does not exist a tubular neighbourhood of the boundary of the corresponding compact set.

\vspace{2mm}

From another perspective, the convexity of the boundary of a Riemannian manifold is of fundamental importance in geometric analysis \cite{li,anton}, and also in building a bridge to the world of metric measure spaces: namely \cite{bxh}, given $K\in\R$, a Riemannian $n$-manifold $(M,g)$ with boundary is an $\mathsf{RCD}(K,n)$ space, if and only if $(M,g)$ is complete with a convex boundary and $\mathrm{Ric}_g\geq K$ (in short: $g$ is an $\mathsf{RCD}$ metric on $M$). In view of these observations, we improve Theorem \ref{exist} as follows:

\begin{thmALPH}[Theorem \ref{thm:bounded-geometry-convex-boundary}]\label{exist2} Let $M$ be a manifold with boundary and $g_0$ a Riemannian metric on $M$. Then there exists a conformal factor $u\in\Csf^\infty(M)$ such that $(M,e^{2u}g_0)$ is of bounded geometry and has a convex boundary.
\end{thmALPH}

While, given a uniform collar of the boundary of $M$ and a two-sided bound on its second fundamental form, it is not difficult to explicitly construct a $\Csf^\infty$-bounded conformal factor which convexifies the boundary, the difficulty in establishing Theorem~\ref{exist2} lies in assuring that the boundary still possesses a uniform collar in the new conformal metric. To this end, we have to establish fine properties of flatzoomers.\vspace{2mm}

Theorem \ref{exist2} plays a fundamental role in connecting $\mathsf{RCD}$ spaces with topology: in \cite{honda}, Honda establishes a theory of orientability on Ricci limit spaces (which in view of \cite{deng} extends to arbitrary finite dimensional $\mathsf{RCD}$ spaces). Now, Theorem \ref{exist2} directly implies that every noncompact manifold with boundary admits an $\mathsf{RCD}$-Riemannian metric. Ultimately, in Theorem \ref{main2}, we prove that Honda's geometric notion of orientability is compatible with the usual topological notion of orientability on mannifolds with boundary.\vspace{4mm}

\noindent
\emph{Notation and conventions}: We understand a \emph{manifold} to be smooth and without boundary. Given a Riemannian manifold $(M,g)$, we will often simply write $(\cdot,\cdot)$ for $g$ and $|\cdot|$ for the corresponding norm. Furthermore, we denote by $\dist  = \dist_g$ the Riemannian distance function induced by $g$ and by $\dist_p\coloneqq \dist(p,\cdot)$ the distance to a point $p\in M$. The symbols $\Bsf_r(x)$ resp. $\overline{\Bsf}_r(x)$ denote an open resp. closed ball of radius $r$ centred at $x$ with respect to $\dist$, and $\Ssf_r(x)$ denotes the corresponding metric sphere. For points $p,q\in M$ we denote by $\gamma_{pq}$ the geodesic (segment) with $\gamma_{pq}(0)=p$ and $\gamma_{pq}(1)=q$, and for $v\in T_pM$ by $\gamma_{pv}$ the geodesic with $\gamma_{pv}(0)=p$ and $\dot\gamma_{pv}(0)=v$.\vspace{1mm}

For an embedded submanifold $\Sigma$ we denote the second fundamental form by $\secfund=\secfund^\Sigma$, and the scalar second fundamental form with respect to a (local) unit normal vector field $\nu$ by $\scalsecfund = \scalsecfund^{\Sigma,\nu} = (\secfund,\nu)$. By the same symbol $\scalsecfund$ we denote also the induced quadratic form which allows us to write inequalities like
\[
    \scalsecfund \geq \lambda\,,
\]
which are, as usual, to be understood in the sense that $\scalsecfund(X)\geq\lambda\cdot (X,X)$ for all $X$.
A tubular neighbourhood of $\Sigma$ is the diffeomorphic image $U$ of a set $\{(x,v)\in N\Sigma\,|\,|v|<\delta(x)\}$ under $\exp$, where $\delta:\Sigma\to (0,\infty)$ is a continuous function which we call the \emph{size of $U$}. A uniform tubular neighbourhood of $\Sigma$ is a tubular neighbourhood of constant size.

\section*{Acknowledgments}
The authors are greatly indebted to Kostiantyn Drach for helpful discussions about \\Blaschke's rolling ball, and Shouhei Honda for bringing the reference \cite{deng} to our attention.

\section{Local radial angle comparison}\label{sec:radial-angle-comp}

Let $(M^n,g)$ be a not necessarily complete Riemannian manifold and $\Sigma\subseteq M$ an embedded hypersurface.

\begin{definition}\label{def:local-0}
	\begin{enumerate}[(i)]
        \item Let $U\subseteq M$ be a connected open set such that
        \[
            \Sigma_U\coloneqq U\cap\Sigma
        \]
        is connected, two-sided and separates $U$ into exactly two components. Then we call $U$ a \emph{$\Sigma$-slice set}.
        \item If $U=\Bsf_r(x)$ with $x\in \Sigma$ is a $\Sigma$-slice set, we call $\Bsf_r(x)$ a \emph{$\Sigma$-slice ball about $x$} and denote $\Sigma_r\coloneqq \Sigma_r(x)\coloneqq \Sigma_U$. If $\Bsf_s(x)\subseteq \Bsf_r(x)$ is a $\Sigma$-slice ball for all $s\leq r$, we call $\Bsf_r(x)$ a \emph{regular $\Sigma$-slice ball}.
        \item Given a $\Sigma$-slice set $U$ and a point $p\in U\setminus \Sigma$ such that $U$ is a normal neighbourhood of $p$ and for each $q\in \Sigma_U$ there exists a unique length minimising geodesic from $p$ to $q$ in $M$, we endow $\Sigma_U$ with the unit normal vector field $\nu$ pointing into the component of $U\setminus\Sigma$ in which $p$ lies and define the \emph{radial angle function with respect to the origin $p$} by
        \begin{gather}\label{def:angle-function}
					\begin{aligned}
						\phi \coloneqq \phi_p : \Sigma_U&\to [0,\pi]\\
						\phi(q) &\mapsto \sphericalangle(\dot{\gamma}_{pq}(1),-\nu_q)\,,
					\end{aligned}
				\end{gather}
				where $\gamma_{pq}$ is the unique geodesic connecting $p$ to $q$, and call $U$ the corresponding \emph{radial angle function domain with respect to $p$}.
	\end{enumerate}
\end{definition}

\begin{remark}\label{rem:continuity-seperation-property-and-slice-balls}
    \begin{enumerate}[(i)]
        \item Note that $\phi_p(q) = \sphericalangle(\grad(\dist_p)_q,-\nu_q)$, so that $\phi_p^{-1}(\{0,\pi\})$ is precisely the critical set of $\dist_p|\Sigma_r$.
        \item Note that $\phi_p$ is continuous, and smooth on $\phi_p^{-1}((0,\pi))$. Observe also that if $U$ is additionally strongly convex so that it is a radial angle function domain with respect to any origin $p \in U\setminus \Sigma$, then \[
            (U\setminus\Sigma_U)\times\Sigma_U\ni (p,q)\mapsto \phi_p(q)\in [0,\pi]
            \]
            is continuous too.
        \item If $\Sigma$ is properly embedded and $U$ is a simply connected open set such that $\Sigma_U$ is connected, then it follows from \cite[Theorem~4.6]{Hir76} that $U$ is a $\Sigma$-slice set, where we remark that the compactness assumption itself is not used in the proof of loc.\ cit., but merely its implication that the hypersurface is properly embedded.

        In particular, if $\Sigma$ is not necessarily properly embedded, $V$ is the domain of a slice chart for $\Sigma$ and $U\subseteq V$ is a simply connected open set such that $\Sigma_U$ is connected, then $U$ is a $\Sigma$-slice set.
    \end{enumerate}
\end{remark}

We fix a radial angle function domain $U$ with origin $p\in U\setminus\Sigma$. Note that the distance function $\dist_p$ is smooth on $U\setminus\{p\}$ and remains so after restricting it to the submanifold $\Sigma_U$. We introduce two smooth vector fields on $\phi_p^{-1}((0,\pi)) = \{q\in\Sigma_U\,|\,\grad(\dist_p)_q\neq \pm\nu_q\}$. For $q\in\phi_p^{-1}((0,\pi))$ let
\begin{align*}
	X_q \coloneqq \frac{\grad^{\Sigma}(\dist_p)_q}{|\grad^{\Sigma}(\dist_p)_q|}\,,
\end{align*}
i.e.\ $X$ is the normalised projection of $\grad\dist_p$ to $T\Sigma_U$. The second vector field, $Y$, is defined by
\begin{enumerate}[(i)]
	\item $Y\in\mathrm{span}\{\nu, X\}$,
	\item $Y\perp \grad^M\dd_p$,
	\item $\sphericalangle(Y,X) = \phi_p$ and $|Y|=1$.
\end{enumerate}
Note that (ii) makes sense since $X$ is, up to normalisation, just the projection of $\grad^M\dist_p$ to $T\Sigma_U$ and that then (iii) uniquely determines $Y$ pointwise on the set $\{X\neq \grad^M\dist_p\}=\{0 < \phi_p < \pi \text{ and } \phi_p\neq \pi/2\}$. Finally, the ambiguity of (iii) on $\{X= \grad^M\dist_p\}=\{\phi_p = \pi/2\}$ is overcome by defining $Y$ to be equal to $\nu$, which turns $Y$ into a smooth vector field.\\

For any $q$ as above, i.e.\ such that $0<\phi_p(q)<\pi$, denote by $\alpha:J_q\to\Sigma_U$ the maximal integral trajectory of $X$ through $q$. Let $\varphi:J_q\to I_q$ be defined by
\[
	\varphi(s)\coloneq \dist_p(\alpha(s))\,,
\]
and set $\gamma=\gamma_q\coloneqq \alpha\circ\varphi^{-1}:I_q\to\Sigma_U$. The curve $\gamma$ is the \emph{unique maximal and by distance to $p$ parametrised integral trajectory} of the gradient vector field of ${\dist_p}{|\Sigma_r}$ through $q$.\\

\begin{figure}
    \centering
    \begin{tikzpicture}
    	\node at (0,1.5) {$U$};
    	\draw[smooth,tension=.7] plot coordinates {(0,0) (2,.7) (4,1) (6,1) (8,1.5) (10,2)} node at (0,.5) {$\Sigma_U$};
        \begin{scope}
            \clip (6,1) circle (2.2);
            \draw[smooth, thick, red,tension=.7] plot coordinates {(0,0) (2,.7) (4,1) (6,1) (8,1.5) (10,2)};
        \end{scope}
        \node[red] at (7.3,1.6) {$\gamma$};
    			
    	\filldraw[black] (2,.7) circle (1pt);
    	\draw[->] (2,.7) -- (2.17, 0) node[left] {$\nu$};
    			
    	\filldraw[black] (4,-1) circle (1pt) node[left] {$p$};
    	\path (4,-1) edge [bend right=10] node[below right] {$\gamma_{pq}$} (6,1);

        \draw [domain=15:95] plot ({4.58286+2.82843*cos(\x)}, {-1.4478+2.82843*sin(\x)});
        \node at (7.9,-0.5) {$\Ssf_t(p)$};

    	\filldraw[black] (6,1) circle (1pt) node[below] {$q$};
    			
    	\draw[->] (6,1) -- (5.9, 2) node[left] {$-\nu$};
    	\draw[->] (6,1) -- (7.056,1.1088);
        \node at (7.3, 1.05) {$X$};
    	\draw[->] (6,1) -- (6.55,1.95) node[above right] {$\dot{\gamma}_{pq}$};
    	\draw[->] (6,1) -- (6.9,.55) node[right] {$Y$};
    			
    	\draw [domain=60:95] plot ({6+0.7*cos(\x)}, {1+0.7*sin(\x)});
    			
    	\path[->] (3.5,2) edge [bend left=20] (6.1,1.4) node[left] {$\phi_p(q)$};
        \end{tikzpicture}
        \caption{The situation of Proposition~\ref{prop:distance-trajectories-of-X}}
    \label{fig:integral-traj-Liouville}
\end{figure}
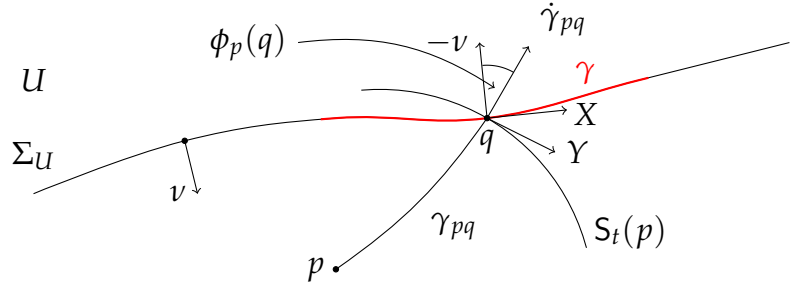

The next proposition shows that $\gamma$ satisfies an ordinary differential equation which ties together the radial angle function, the curvature of $\Sigma$ and the curvature of the geodesic sphere centred at $p$ and intersecting $\Sigma$ in $\gamma(t)$, see Figure~\ref{fig:integral-traj-Liouville}. This can be seen as a generalisation of Liouville's formula for curves in surfaces, see, e.g. \cite[Chapter~4.4, Proposition~4]{doC16}.

\begin{proposition}\label{prop:distance-trajectories-of-X}
	Suppose $U$ is a radial angle function domain with origin $p\in U\setminus\Sigma$ and radial angle function $\phi_p$. Then, for any maximal and by distance to $p$ parametrised integral trajectory $\gamma$, one has
	\begin{align*}
		\scalsecfund^{\Sigma,\nu}(X_{\gamma(t)}) = \scalsecfund^{\Ssf_t(p),-\grad^M\dd_p}(Y_{\gamma(t)})\cdot\cos\phi_p(\gamma(t)) + \frac{\dd}{\dd t}\cos\phi_p(\gamma(t))\,.	
	\end{align*}	
\end{proposition}
\begin{proof}
    This is a localised version of \cite[Proposition~3.1]{Dra24}.
\end{proof}

To set up the radial angle comparison theorem, we need to introduce the objects to which we compare the radial angle function with respect to an arbitrary origin in a given $\Sigma$-slice ball. To this end, denote for $c\in\R$ by $\mathbb{M}^n(c)$ an $n$-dimensional complete simply connected model space of constant sectional curvature $c$ and for $\lambda\in\R$, let $\Ssf^c_\lambda\subseteq \mathbb{M}^n(c)$ be a connected and properly embedded totally umbilical hypersurface of constant principal curvature $\lambda$, i.e.\ a connected properly embedded hypersurface with $\scalsecfund = \lambda$. The following list exhaustively describes these, see e.g.\ \cite[p.\ 177 ff]{doC92} for the case $c<0$.

\begin{itemize}
    \item $c>0$: in this case, $\mathbb M^n(c)$ is (isometric to) a sphere of radius $1/\sqrt{c}$ in euclidean space and each $\Ssf^c_\lambda$ is a geodesic sphere with the normal vector field pointing towards its centre if $\lambda > 0$ and towards the component not containing the centre if $\lambda < 0$, whereas for $\lambda = 0$, $\Ssf^c_\lambda$ is an equator of $\mathbb M(c)$ and either choice of a unit normal vector field is geometrically equivalent to the other.
    \item $c = 0$: in this case, $\mathbb M^n(c)$ is isometric to euclidean space and there exist two types of totally umbilical hypersurfaces.
        \begin{itemize}
            \item[$\circ$] $\lambda\neq 0$: $\Ssf^c_\lambda$ is again a geodesic sphere with its unit normal vector field analogously to the case $c>0$.
            \item[$\circ$] $\lambda = 0$: $\Ssf^c_0$ is an euclidean hyperplane and either choice of unit normal vector field is geometrically equivalent to the other.
        \end{itemize}
    \item $c<0$: $\mathbb M^n(c)$ is isometric to a rescaled hyperbolic space and there are four different types of properly embedded totally umbilical hypersurfaces.
    \begin{itemize}
        \item[$\circ$] $|\lambda|>\sqrt{-c}$: analogously to the previous two cases, $\Ssf^c_\lambda$ is a geodesic sphere.
        \item[$\circ$] $|\lambda|=\sqrt{-c}$: $\Ssf^c_\lambda$ is a horosphere, i.e.\ the limit of a sequence of geodesic spheres with increasing unbounded radii which share a common tangent plane at a fixed point.
        \item[$\circ$] $\sqrt{-c}>|\lambda|>0$: in this case, $\Ssf^c_\lambda$ is one component of an equidistant hypersurface.
        \item[$\circ$] $\lambda = 0$:  $\Ssf^c_\lambda$ is totally geodesic and thus isometric to $\mathbb M^{n-1}(c)$.
    \end{itemize}
\end{itemize}
Note that in each case $\Ssf^c_\lambda$ separates $\mathbb{M}^n(c)$ into two components. We call the one into which its normal points \emph{the inside of $\Ssf^c_\lambda$}.

For $c\leq 0$, $U=\mathbb{M}^n(c)$ is strongly convex so that it is a radial angle function domain with origin $p_\lambda$ for any $p_\lambda\in\mathbb{M}^n(c)\setminus\Ssf^c_\lambda$. For $c>0$ the situation is more complicated, but we will restrict to hypersurfaces $\Ssf^c_\lambda$ with $\lambda\geq0$ and points $p_\lambda$ on the inside of $\Ssf^c_\lambda$, in which case $U=\mathbb{M}^n(c)\setminus\{-p_\lambda\}$ is a radial angle function domain with origin $p_\lambda$ in which $\Ssf^c_\lambda$ is completely contained and where $-p_\lambda$ denotes the antipodal point to $p_\lambda$ in $\mathbb M^n(c)$.\\

If $\Ssf^c_\lambda$ happens to be a geodesic sphere, we denote by $\compR^c_\lambda$ its radius and set $\compR^c_\lambda\coloneqq\infty$ in all other cases. Recall the generalised sine function
\begin{equation}\label{eqn:snc}
	\mathrm{sn}_c(t) = \begin{cases}
		\tfrac{1}{\sqrt{c}}\sin\sqrt{c}t, & c  > 0\,,\\
		t, & c = 0\,,\\
		\tfrac{1}{\sqrt{-c}}\sinh\sqrt{-c}t, & c < 0\,,\\
	\end{cases}
\end{equation}
and the associated generalised cotangent
\begin{equation}\label{eqn:ctc}
	\mathrm{ct}_c(t)   = \frac{\mathrm{sn}_c'(t)}{\mathrm{sn}_c(t)} = \begin{cases}
				\sqrt{c}\cot\sqrt{c}t, & c  > 0\,,\\
		\frac 1t, & c = 0\,,\\
		\sqrt{-c}\coth\sqrt{-c}t, & c < 0\,,\\
	\end{cases}
\end{equation}
which connects principal curvatures and radii of geodesic spheres in constant curvature spaces \cite[Proposition~11.3]{lee},
\[
	\mathrm{ct}_c(\compR_\lambda^c) = |\lambda|\,.
\]
We are now set up for the main theorem of this section.
\begin{theorem}[Local two-sided BD-radial angle comparison]\label{thm:rad-angle-comp}
Let $(M^n,g)$ be a Riemannian manifold, $\Sigma\subseteq M$ an embedded hypersurface and $\Bsf_r(x)$ a precompact, strongly convex $\Sigma$-slice ball about $x\in \Sigma$ such that $\Bsf_{r/3}(x)$ is a $\Sigma$-slice ball too. Let $p\in \Bsf_{r/3}(x)\setminus \Sigma_{{r/3}}$ and endow $\Sigma_r$ with the unit normal vector field $\nu$ pointing towards the component of $\Bsf_r(x)\setminus \Sigma_r$ in which $p$ lies. Let $c,\lambda \leq 0$ and  $C,\Lambda\geq 0$ be real numbers such that
\begin{equation}\label{eqn:loc-2-side-angle-comp-0}
	\phantom{\quad\sigma\in \mathrm{Gr}_2(T_yM), y\in B_r(x)}c\leq \seccurv(\sigma) \leq C\quad \text{for all}\quad \sigma\in \mathrm{Gr}_2(T_yM), y\in \Bsf_r(x)
\end{equation}
and 
\begin{equation}\label{eqn:loc-2-side-angle-comp-1}
	\phantom{\quad\text{for all}\quad X\in Ty\Sigma, y\in \Sigma_r\,.}\lambda\leq \scalsecfund^{\Sigma_r,\nu}(X) \leq \Lambda \quad\text{for all}\quad X\in T_y\Sigma_r, y\in \Sigma_r\,.
\end{equation}
Assume
\begin{gather}
    \label{eqn:loc-2-side-angle-comp-2}
	\phantom{\text{ if }\quad \Lambda \leq 0\,,}
	\begin{aligned}
		\frac 56 \cdot r \leq \compR_\Lambda^C\,,\\
		\frac 23\cdot r\leq \compR_\lambda^c\,.
	\end{aligned}
\end{gather}

Furthermore, assume there exists a point $s\in\Sigma_{r/3}$ with $d_0\coloneqq\dist_p(s) = \dist_p(\Sigma_{r})$ and such that $s$ is an isolated critical point of $\dist_p|\Sigma_r$. Let $p_\lambda \in \mathbb{M}^2(c)$ be a point on the inside of $\Ssf_\lambda^c$ with $\dist(p_\lambda,\Ssf_\lambda^c) = d_0$ and $\phi_\lambda$ be the associated radial angle function with origin $p_\lambda$. Likewise, let $p_\Lambda\in \mathbb{M}^2(C)$ be a point on the inside of $\Ssf_\Lambda^C$ with $\dist(p_\Lambda,\Ssf_\Lambda^C)=d_0$ and $\phi_\Lambda$ be the associated radial angle function with origin $p_\Lambda$.

Then, the radial angle function $\phi_p:\Sigma_{r}\to[0,\pi]$ with respect to the origin $p$ satisfies for all $q\in\Sigma_{r/3}$, $q_\lambda\in\Scal_\lambda^c$ and $q_\Lambda\in\Scal_\Lambda^C$ such that $\dist(p,q) = \dist(p_\lambda, q_\lambda) = \dist(p_\Lambda, q_\Lambda)$:
\begin{equation}\label{eqn:radial-angle-comp-phi-inequ}
		\phi_\Lambda(q_\Lambda)\leq \phi_p(q)\leq \phi_\lambda(q_\lambda)\,.
\end{equation}
\end{theorem}
\begin{proof}
	Assumption \eqref{eqn:loc-2-side-angle-comp-2} establishes that for every $q\in \Sigma_r$ there exist points $q_\Lambda \in \Ssf_\Lambda^C$ and $q_\lambda\in\Ssf_\lambda^c$ such that $\dist(p,q)=\dist(p_\Lambda, q_\Lambda) = \dist(p_\lambda, q_\lambda)$. Indeed, for every $q\in\Sigma_r$ we have 
    \begin{equation}\label{eqn:distance-0}
        \dist(p,q)\leq d_2\coloneqq \sup_{q\in\Sigma_r}\dist(p,q) < \frac 43\cdot r\,.
    \end{equation}
    Assuming that $\Ssf^C_\Lambda$ is a geodesic sphere, i.e.\ $C>0$, and letting $s_\Lambda\in\Ssf^C_\Lambda$ be the unique point with $\dist(p_\Lambda,s_\Lambda)=\dist(p_\Lambda,\Ssf^C_\Lambda)$, the point in $\Ssf^C_\Lambda$ which maximises distance to $p_\Lambda$ is $-s_\Lambda$, the antipodal to $s_\Lambda$ in $\Ssf^C_\Lambda$, for which we have 
    \begin{equation}\label{eqn:distance-1}
        \dist(p_\Lambda,-s_\Lambda)=2\compR^C_\Lambda-\dist(p_\Lambda,s_\Lambda) = 2\compR^C_\Lambda-\dist(p,s)> \frac 53 r-\frac 13 r = \frac 43r\,,
    \end{equation} 
    and analogously for $\Ssf^c_\lambda$.

    There are at most two points $q_\Lambda$ in $\Ssf^C_\Lambda$ for which $\phi_\Lambda(q_\lambda)\in \{0,\pi\}$, namely $q_\Lambda = s_\Lambda$ with $\phi_\Lambda(s_\lambda)=0$ and for $C>0$ the point $q_\Lambda = -s_\Lambda$ with $\phi_\Lambda(q_\Lambda)=0$. Similarly, we have $\phi_\lambda(s_\lambda)=0$ and in case $\Ssf^c_\lambda$ is a geodesic sphere, $\phi_\lambda(-s_\lambda)=\pi$, while all points $q_\lambda\in\Ssf^c_\lambda$ with $q_\lambda\neq s_\lambda,-s_\lambda$ satisfy $0<\phi_\lambda(q_\lambda)<\pi$.\\
    
 
    Clearly, we have \eqref{eqn:radial-angle-comp-phi-inequ} for $q=s$, $q_\lambda=s_\lambda$ and $q_\Lambda=s_\Lambda$, so that in the following, we can assume $q\neq s$.\\

    Since $s$ is a global minimum and an isolated critical point of $\dist_p|\Sigma_r$, there exists a neighbourhood $U_0$ of $s$ in $\Sigma_r$ such that for each $q\in U_0\setminus\{s\}$ the maximal and by distance to $p$ parametrised integral trajectory $\gamma=\gamma_q:I_q=(d_0,d_1)\to\Sigma_r$ of $\dist_p|\Sigma_r$ satisfies $\lim_{t\searrow d_0}\gamma(t) = s$. Fix such a point $q\in U_0$, choose points $q_\lambda\in\Ssf^c_\lambda$ and $q_\Lambda\in \Ssf^C_\Lambda$ such that $\dist(p,q)=\dist(p_\lambda,q_\lambda)=\dist(p_\Lambda,q_\Lambda)$ and denote by $\gamma_\lambda=\gamma_{p_\lambda}$ and $\gamma_\Lambda=\gamma_{p_\Lambda}$ the by distance to $p_\lambda$ resp.\ $p_\Lambda$ parametrised integral trajectories of $\dist_{p_\lambda}|\Ssf^c_\lambda$ resp.\ $\dist_{p_\Lambda}|\Ssf^c_\Lambda$. Note that by \eqref{eqn:distance-0} and \eqref{eqn:distance-1} the domains of definition of $\gamma_\lambda$ and $\gamma_\Lambda$ contain $I_q$.

	By Proposition~\ref{prop:distance-trajectories-of-X}, along these curves we have
	\begin{align}
		\scalsecfund^{\Sigma_r,\nu}\left(X_{\gamma(t)}\right) &= \scalsecfund^{ \Ssf_t(p)}\left(Y_{\gamma(t)}\right)\cdot\cos\phi_p(\gamma(t)) + \frac{\dd}{\dd t}\cos\phi_p\left(\gamma(t)\right)\,,\label{eqn:radial-angle-comp-0}\\
		\lambda &= \scalsecfund^{  \Ssf_t(p_\lambda)}\left(Y^{\mathbb{M}(c)}_{\gamma_\lambda(t)}\right)\cdot\cos\phi_\lambda(\gamma_\lambda(t)) + \frac{\dd}{\dd t}\cos\phi_\lambda(\gamma_\lambda(t))\,,\label{eqn:radial-angle-comp-1}\\
		\Lambda &= \scalsecfund^{  \Ssf_t(p_\Lambda)}\left(Y^{\mathbb{M}(C)}_{\gamma_\Lambda(t)}\right)\cdot\cos\phi_\Lambda(\gamma_\Lambda(t)) + \frac{\dd}{\dd t}\cos\phi_\Lambda(\gamma_\Lambda(t))\,,\label{eqn:radial-angle-comp-2}
	\end{align} 
	Using the principal curvature comparison theorem for level sets of distance functions \cite[Corollary~11.8]{lee},
	\begin{align*}
		 \mathrm{ct}_C(t) = \scalsecfund^{  \Ssf_t(p_\Lambda)}\left(Y^{\mathbb{M}(C)}_{\gamma_\Lambda(t)}\right) \leq \scalsecfund^{ \Ssf_t(p)}\left(Y_{\gamma(t)}\right)\leq \scalsecfund^{\Ssf_t(p_\lambda)}\left(Y^{\mathbb{M}(c)}_{\gamma_\lambda(t)}\right) = \mathrm{ct}_c(t)\,,
	\end{align*}
	we subtract \eqref{eqn:radial-angle-comp-1} from \eqref{eqn:radial-angle-comp-0} and \eqref{eqn:radial-angle-comp-0} from \eqref{eqn:radial-angle-comp-2} to obtain
	
	\begin{align*}
		0 &\leq \scalsecfund^{\Sigma_r,\nu}\left(X_{\gamma(t)}\right) - \lambda \leq (\cos\phi_p(t)-\cos\phi_\lambda(t))\cdot\ct_c(t) + \frac{\dd}{\dd t}(\cos\phi_p(t)-\cos\phi_\lambda(t))\,,\\
		0 &\leq \Lambda - \scalsecfund^{\Sigma_r,\nu}\left(X_{\gamma(t)}\right) \leq (\cos\phi_\Lambda(t)-\cos\phi_p(t))\cdot\ct_C(t) + \frac{\dd}{\dd t}(\cos\phi_\Lambda(t)-\cos\phi_p(t))\,,	
	\end{align*}
	where we have used the short-hand notation $\phi_p(t)\coloneqq \phi_p(\gamma(t))$ and likewise for all other involved functions.

	Setting $f(t)\coloneqq \cos\phi_p(t)-\cos\phi_\lambda(t)$ and $h(t)\coloneqq\cos\phi_\Lambda(t)-\cos\phi_p(t)$, the above inequalities imply
	\begin{align*}
		0&\leq f(t)\cdot\ct_c(t) + f'(t) = \frac{1}{\mathrm{sn}_c(t)}\left(f(t)\mathrm{sn}'_c(t)+f'(t)\mathrm{sn}_c(t)\right)\,,\\
		0 & \leq h(t)\cdot\ct_C(t) + h'(t)= \frac{1}{\mathrm{sn}_C(t)}\left(h(t)\mathrm{sn}'_C(t)+h'(t)\mathrm{sn}_C(t)\right)\,,
	\end{align*}
	which is equivalent to
	\begin{align*}
		0&\leq \frac{\dd}{\dd t}(f(t)\cdot\sn_c(t))\,,\\
		0 & \leq  \frac{\dd}{\dd t}(h(t)\cdot\sn_C(t))\,.
	\end{align*}
	Hence, the functions $f\cdot\sn_c$ and $h\cdot\sn_C$ are monotonically nondecreasing.

        By choice of the point $q$ we have $\lim_{t\searrow d_0}f(t)\cdot\mathrm{sn}_c(t)=\lim_{t\searrow d_0}h(t)\cdot\mathrm{sn}_C(t)=0$ so that $f$ and $h$ are nonnegative functions, i.e.\ we have established \eqref{eqn:radial-angle-comp-phi-inequ} for all points in $\operatorname{img}\gamma$ with appropriate comparison points.
 
        Let $\delta\in I_q=(d_0,d_1)$ and define $\underline{\phi}\coloneqq \inf_{t\in [\delta,d_1)}\phi_\Lambda(t)$ and $\overline{\phi}\coloneqq \sup_{t\in[\delta,d_1)}\phi_\lambda(t)$. Since $\delta > d_0$ and $d_1\leq d_2$, we have by \eqref{eqn:distance-0} and \eqref{eqn:distance-1}
        \begin{equation}\label{eqn:distance-1.5}
            0<\underline{\phi}\leq \phi_\Lambda(t)\leq \phi_p(t)\leq\phi_\lambda(t)\leq\overline{\phi} < \pi
        \end{equation}
        for all $t\in [\delta,d_1)$. Hence, there exists $\epsilon > 0$ such that 
        \begin{equation}\label{eqn:distance-2}
            \epsilon < \pi/2 - \sphericalangle(\grad(\dist_p)_{\gamma(t)},\dot\gamma(t))
        \end{equation}
        for all $t\in [\delta,d_1)$. The curve $\gamma$ is parametrised by distance to $t$, i.e.\ $t = \dist_p(\gamma(t))$. Differentiating this, we obtain
        \[
            1 = \dd(\dist_p)_{\gamma(t)}\dot\gamma(t) = (\grad (\dist_p)_{\gamma(t)},\dot\gamma(t)) = |\grad(\dist_p)_{\gamma(t)}|\cdot|\dot\gamma(t)|\cdot\cos\sphericalangle(\grad(\dist_p)_{\gamma(t)},\dot\gamma(t))\,,
        \]
        which, together with \eqref{eqn:distance-2} and the fact that $|\grad(\dist_p)|=1$, shows that $\sup_{t\in[\delta,d_1)}|\dot\gamma(t)|<\infty$. In particular, $\gamma|[\delta,d_1)$ has finite length so that $o\coloneqq \lim_{t\nearrow d_1}\gamma(t)$ exists by compactness of $\overline{\Bsf}_r(x)$. Assuming that $o\in\Sigma_r$ we have $0<\phi_p(o)<\pi$ by \eqref{eqn:distance-1.5} and continuity, i.e.\ $o$ is a regular point of $\dist_p$ which contradicts the maximality of $\gamma$. Hence, $o\in\partial\Sigma_r\subseteq \Ssf_r(x)$, meaning that $\gamma$ is defined until it hits the boundary of $\Bsf_r(x)$.\\

        Let $U\subseteq \Sigma_r$ be the set consisting of the point $s$ and all points $q$ such that the maximal trajectory through $q$ emanates from $s$, i.e.\ the basin of attraction of the attractor $s$ of the (negative) gradient system of $\dist_p|\Sigma_r$. This is an open set which is invariant under the flow $\Phi$ generated by the vector field $\grad(\dist_p)$ and on which we have shown \eqref{eqn:radial-angle-comp-phi-inequ} with appropriate comparison points. We wish to show that $\Sigma_{r/3}\subseteq U$.
        
        Since $\partial U\subseteq \overline{\Bsf}_r(x)$ is compact, $\min_{\partial U}\dist_p$ is attained at a point $o\in\partial U$. Assuming $o\in\Sigma_r$, we have $\phi_\Lambda(o_\Lambda)\leq\phi_p(o)\leq \phi_\lambda(o_\lambda)$ by continuity for appropriate comparison points $o_\Lambda$ and $o_\lambda$. The set $U$ is a neighbourhood of $s$ so that $\dist(p,o) > d_0$. It follows from this and \eqref{eqn:distance-0} and \eqref{eqn:distance-1} that
        \begin{equation*}
            0 < \phi_\Lambda(o_\Lambda)\leq\phi_p(o)\leq \phi_\lambda(o_\lambda) < \pi\,,
        \end{equation*}
        i.e.\ $o$ is a regular point of $\dist_p|\Sigma_r$. Hence, there exists $\varepsilon > 0$ such that $o'\coloneqq \Phi_{-\varepsilon}(o)$ is defined and not equal to $o$. The point $o'$ does not belong to $U$, since $U$ is a $\Phi$-invariant set. Using $\Phi$, one easily sees that every neighbourhood of $o'$ contains points of $U$, so that $o'\in\partial U$. Since $\dist_p(o')<\dist_p(o)$, this is a contradiction to the definition of $o$. Hence, $\min_{\partial U}\dist_p$ is attained at points in $\Ssf_r(x)$ and at no points in $\partial U\cap\Sigma_r$.

        Suppose now that $\Sigma_{r/3}\not\subseteq U$. Then there exists a point $q\in\Sigma_{r/3}$ with $q\not\in U$. By assumption, $\Sigma_{r/3}$ is connected. Choose a continuous path $\sigma:[0,1]\to\Sigma_{r/3}$ with $\sigma(0)=s$ and $\sigma(1)=q$. Let $t_0\in (0,1)$ be the first time such that $\sigma(t_0)\in\partial U$. By the above, 
        \[
            \dist(p,\sigma(t_0)) \geq \dist\left(p,\Ssf_r(x)\right) > \frac 23r\,,
        \]
        which contradicts the fact that $\sigma$ is a path in $\Sigma_{r/3}$. Hence, $\Sigma_{r/3}\subseteq U$ as claimed.
\end{proof}

A few remarks are in order. Obviously, one could relax the assumptions of the theorem to allow for $c$ and $C$ to have the same sign and likewise for $\lambda$ and $\Lambda$. While there are situations in which these would give better estimates, we will only need the theorem as recorded above.

Another rather obvious generalisation, having the Rauch comparison theorem in mind, is to compare the radial angle function of $\Sigma$ not to ones of totally umbilical hypersurfaces in constant curvature spaces but to such of abstract hypersurfaces in arbitrary manifolds satisfying appropriate curvature assumptions in place of \eqref{eqn:loc-2-side-angle-comp-0} and \eqref{eqn:loc-2-side-angle-comp-1}.

It is noteworthy that in the following we will only need the left-hand inequality in \eqref{eqn:radial-angle-comp-phi-inequ}, but that the proof of Theorem~\ref{thm:rad-angle-comp} works by establishing the two inequalities simultaneously.\\

\begin{figure}
    \centering
    \begin{tikzpicture}
        \draw[smooth cycle,tension=.7] plot coordinates{(-0.707,0.707) (-1.5,.3) (-2,-1) (-1.5,-1.3) (-0.8,-0.4) (0,-1) (0.707,-0.707) (1,0) (0.707,0.707) (0,1)} node at (1.5, 0) {$\Sigma_r$};
        \filldraw[black] (0,0) circle (1pt) node at (0.4, -0.15) {$s$};
        \draw (0,0) -- (-0.707,0.707);
        \draw (0,0) -- (0,1);
        \draw (0,0) -- (0.707,0.707);
        \draw (0,0) -- (1,0);
        \draw (0,0) -- (0.707,-0.707);
        \draw (0,0) -- (0,-1);
        \draw (0,0) -- (-1.7,0);
        \draw (0,0) -- (-1.85,-0.4);
        \draw[red] (0,0) -- (-2,-1);
        \draw (0,0) -- (-1.5,.3);
    \end{tikzpicture}
        \caption{A gradient line of ${\dd_p}{|\Sigma_r}$ that is tangent to $\partial \Bsf_r(x)$.}
        \label{fig:gradient-lines-tangent-to-boundary}
\end{figure}
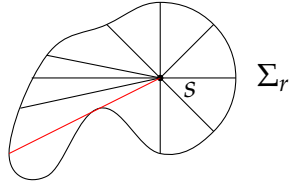

The failure of \eqref{eqn:radial-angle-comp-phi-inequ} to hold for all $q\in\Sigma_r(x)$ stems from the topology and geometry of $\Sigma_r$ which, in turn, is related to gradient lines of $\dd_p|\Sigma_r$ hitting $\partial\Bsf_r(x)$ tangentially. This leads to subsets of $\Sigma_r$ that cannot be seen by the gradient lines of $\dd_p|\Sigma_r$ emanating from $s$, see Figure~\ref{fig:gradient-lines-tangent-to-boundary}.\\

\section{A local Blaschke ball rolling theorem and the normal injectivity radius}\label{sec:local-blaschke-and-klingenberg}

In this section, we prove the local two-sided inner Blaschke ball rolling Theorem~\ref{thm:local-Blaschke} and the pointwise lower estimate on the normal injectivity radius of an embedded hypersurface, Theorem~\ref{thm:Klingenberg-hypersurface}. We continue with our setup of an embedded hypersurace $\Sigma$ in an arbitrary Riemannian manifold $(M^n,g)$ and begin with a classical result in comparision geometry.

\begin{remark}[Hinge comparison]\label{rem:Rauch-triangle-comp}
    \begin{enumerate}[(i)]
        \item The following is a corollary to the Rauch comparison theorem and can be easily deduced from \cite[Corollary~1.35]{CE08}. Let $\Bsf_r(y)\subseteq M$ be a geodesic ball and $c\in\R$ such that
	\begin{equation*}
	\phantom{\quad \text{for all}\quad \sigma\in \mathrm{Gr}_2(TpM), p\in B_r(y)} \seccurv(\sigma) \lesseqgtr c\quad \text{for all}\quad \sigma\in \mathrm{Gr}_2(T_pM), p\in B_r(y)\,.
	\end{equation*}
        Let $x,z\in \Bsf_r(y)$ and choose points $x_c,y_c,z_c\in\mathbb{M}^n(c)$ such that
	\begin{align*}
		\dist(x,y) = \dist(x_c,y_c)\,, & &\dist(y,z) = \dist(y_c,z_c)\,, && \sphericalangle x_cy_cz_c = \sphericalangle xyz\,,
	\end{align*}   	
        where $\sphericalangle xyz$ denotes the angle at $y$ spanned by $\dot\gamma_{yx}(0)$ and $\dot\gamma_{yz}(0)$.
        If $c>0$, assume $r\leq \pi/\sqrt{c}$ so that $\exp^{\mathbb{M}^n(c)}|\Bsf_r(0)$ is nonsingular. Then,
	\[
		\dist(x_c,z_c) \lesseqgtr \dist(x,z)\,.
	\]
        \item  It follows from the law of cosines \cite[Proposition~12.2.3]{Pe16} that in a constant curvature space $\mathbb M^n(c)$, the length $\dist(x_c,z_c)$ of the closing side of a hinge is monotone in the angle $\sphericalangle x_cy_cz_c$. Hence, the assumption $\sphericalangle x_cy_cz_c = \sphericalangle xyz$ in (i) can be replaced by \[\sphericalangle x_cy_cz_c \lesseqgtr \sphericalangle xyz\,.\]
        \item It is well-known that in a constant curvature space, three points always span a geodesic triangle which lies in a totally geodesic surface. Hence, in (i) we can always choose our comparison points to be in $\mathbb M^2(c)$.
    \end{enumerate}
\end{remark}

In the following, we will use the norm $|\secfund^\Sigma| = \max_{|X|=1}|\secfund(X,X)|$, so that for all $\lambda\geq 0$ the inequality $|\secfund|\leq \lambda$ is true if and only if  for all $X$ one has $-\lambda\leq \scalsecfund^{\Sigma,\nu}(X)\leq \lambda$ for any choice of local unit normal $\nu$.

\begin{theorem}[Local two-sided inner Blaschke ball rolling theorem]\label{thm:local-Blaschke}
    Let $(M,g)$ be a Riemannian manifold and $\Sigma\subseteq M$ an embedded hypersurface. Let $\Bsf_r(x)$ be a precompact, strongly convex $\Sigma$-slice ball and assume $\Bsf_{r/3}(x)$ is a $\Sigma$-slice ball too. Let $c,\lambda \geq 0$ be real numbers such that
    \begin{itemize}
        \item $|\seccurv(\sigma)| \leq c$ for all $\sigma\in \mathrm{Gr}_2(T_yM), y\in \Bsf_r(x)$,
        \item $|\secfund^\Sigma|\leq \lambda$ on $\Sigma_r$,
        \item $\frac 56 \cdot r \leq \compR_\lambda^c$.
    \end{itemize}
    Then, any ball of radius $<r/3$ can roll freely on either side of $\Sigma_{r/3}$ inside $\Bsf_{r/3}(x)$, i.e.\ if $\rho<r/3$, $o = \exp_q(v)\in \Bsf_{r/3}(x)$ with $q\in\Sigma_{r/3}$ and $v\in N_q\Sigma$ with $|v| = \rho$, then $\overline\Bsf_\rho(o)\cap \Sigma_{r/3} = \{q\}$.
	\end{theorem}
\begin{proof}
    Endow $\Sigma_r$ with the unit normal $\nu$ that points into the component in which $o$ lies. Let $s\in\Sigma_{r/3}\setminus\{q\}$. We will show that $\dd(o,s)>\rho$.
    
    By choosing some tubular neighbourhood of $\Sigma_r$, we see that there exists an $\varepsilon_0>0$ such that with $p:[0,\varepsilon_0]\to\mathbb{M}^2(c)$ defined by $p(\varepsilon)\coloneqq \gamma_{s \nu_s}(\varepsilon)$, where $\gamma_{s\nu_s}$ is the geodesic with $\gamma_{s\nu_s}(0)=s$ and $\gamma_{s\nu_s}'(0)=\nu_s$, the point $s$ is an isolated critical point of $\dist_{p(\varepsilon)}|\Sigma_r$ and $\dist(s,p(\varepsilon))=\varepsilon$ for all $\varepsilon\in (0,\varepsilon_0]$.
    
    Let $\phi_{p(\varepsilon)}$ be the associated radial angle function. We fix a point $s_\lambda\in\Ssf_\lambda^c\subseteq\mathbb{M}^2(c)$ and choose a continuous curve $p_\lambda:[0,\varepsilon_0]\to\mathbb{M}^2(c)$ such that $p_\lambda(\varepsilon)$ lies on the inside of $\Ssf^c_\lambda$ and such that $\dist(p_\lambda(\varepsilon), s_\lambda) = \dist(p_\lambda(\varepsilon),\Ssf_\lambda^c)=\varepsilon$. Furthermore, we choose continuous curves $q_\lambda:[0,\varepsilon_0]\to\Ssf^c_\lambda$ such that $\dist(p(\varepsilon),q)=\dist(p_\lambda(\varepsilon), q_\lambda(\varepsilon))$ and $o_\lambda:[0,\varepsilon]\to\mathbb{M}^2(c)$ with $o_\lambda(\varepsilon)$ in the inside of $\Ssf_\lambda^c$ such that $\dist(o_\lambda(\varepsilon),q_\lambda(\varepsilon))=\dist(o_\lambda(\varepsilon),\Ssf_\lambda^c)=\rho$. Note that, since $\rho < r/3 < \compR^c_\lambda$ and by virtue of $\varepsilon_0$ being small enough, there exists $\delta>0$ such that $\dist(o_\lambda(\varepsilon),p_\lambda(\varepsilon))\geq\rho+\delta$ for all $\varepsilon\in [0,\varepsilon_0]$. Let $\phi_\lambda = \phi_{p_\lambda(\varepsilon)}$ be the associated radial angle function. Then, noting that $\compR^c_\lambda\leq \compR^{-c}_{-\lambda}$, we have by construction and Theorem~\ref{thm:rad-angle-comp},
	\begin{equation*}
		\phantom{\qquad \text{ for all } \qquad 0<\varepsilon\leq \varepsilon_0}\phi_\lambda(q_\lambda(\varepsilon))\leq \phi_{p(\varepsilon)}(q) \qquad \text{ for all } \qquad 0<\varepsilon\leq \varepsilon_0\,.
	\end{equation*}
	
	The points $o,q,p(\varepsilon)$ and $o_\lambda(\varepsilon), q_\lambda(\varepsilon), p_\lambda(\varepsilon)$ thus satisfy the assumptions of the hinge comparison with upper curvature bounds, Remark~\ref{rem:Rauch-triangle-comp}$(\leq)$, and so we obtain \[\dist(o, p(\varepsilon))\geq \dist(o_\lambda(\varepsilon), p_\lambda(\varepsilon))>\rho+\delta\,.\] By joint continuity of the angle functions with respect to their arguments and origins, this holds true for every $0<\varepsilon \leq\varepsilon_0$, so that, as $\epsilon\to 0$, we obtain $\dist(o,s)\geq \dist(o_\lambda(0), s_\lambda)\geq\rho+\delta>\rho$.
\end{proof}

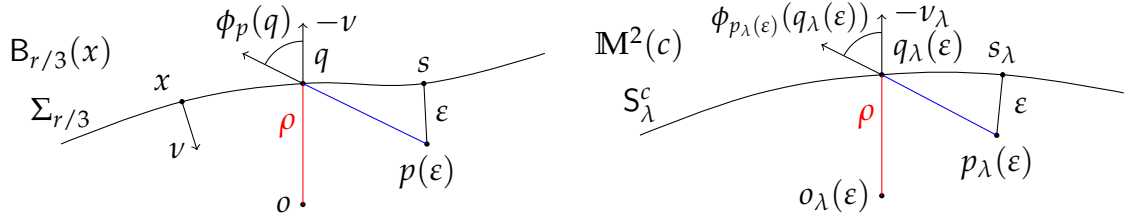
\begin{figure}\label{fig:proof-of-Blaschke}
    \begin{subfigure}[b]{0.45\textwidth}
    		\begin{tikzpicture}[scale=0.8]
			\node at (0,1.5) {$\Bsf_{r/3}(x)$};
			\draw[smooth,tension=.7] plot coordinates {(0,0) (2,.7) (4,1) (6,1) (8,1.5)} node at (0,.5) {$\Sigma_{r/3}$};
			
			\filldraw[black] (2,.7) circle (1pt) node[above left] {$x$};
            \draw[->] (2,.7)--(2.25,-0.1) node[left] {$\nu$};
   
			\filldraw[black] (4,1) circle (1pt) node[above right] {$q$};
			\draw[->] (4,1)--(4,2) node[right] {$-\nu$};
			
			\filldraw[black] (4,-1) circle (1pt) node[left] {$o$};
			\draw[red] (4,1) -- (4,-1) node at (3.75,0.3) {$\rho$};
			
			\filldraw[black] (6,1) circle (1pt) node[above] {$s$};
			
			\filldraw[black] (6.05,0) circle (1pt) node[below] {$p(\varepsilon)$};
			\draw (6,1) -- (6.05,0) node at (6.3,0.5) {$\varepsilon$};
			
                \draw[blue] (6.05,0) -- (4,1);

                \draw[->] (4,1) -- (3.0,1.5);
   			\draw [domain=90:153] plot ({4+0.7*cos(\x)}, {1+0.7*sin(\x)});
			
			\node at (3.2,2) {$\phi_p(q)$};			
    	\end{tikzpicture}
     \end{subfigure}
     \begin{subfigure}[b]{0.45\textwidth}
		\begin{tikzpicture}[scale=0.8]
			\node at (0,1.5) {$\mathbb{M}^2(c)$};
			\draw[smooth,tension=.7] plot coordinates {(0,0) (2,.7) (4,1) (6,1) (8,.7)} node at (0,.5) {$\Ssf^c_\lambda$};
			
			\filldraw[black] (4,1) circle (1pt) node[above right] {$q_\lambda(\varepsilon)$};
			\draw[->] (4,1)--(4,2) node[right] {$-\nu_\lambda$};
			
			\filldraw[black] (4,-1) circle (1pt) node[left] {$o_\lambda(\varepsilon)$};
			\draw[red] (4,1) -- (4,-1) node at (3.75,0.3) {$\rho$};
			
			\filldraw[black] (6,1) circle (1pt) node[above] {$s_\lambda$};
			
			\filldraw[black] (5.9,0) circle (1pt) node[below] {$p_\lambda(\varepsilon)$};
			\draw (6,1) -- (5.9,0) node at (6.3,0.5) {$\varepsilon$};

                \draw[blue] (5.9,0) -- (4,1);
   
			\draw[->] (4,1) -- (3.0,1.5);		
			\draw [domain=90:153] plot ({4+0.7*cos(\x)}, {1+0.7*sin(\x)});

			\node at (2.5,1.9) {\small $\phi_{p_\lambda(\varepsilon)}(q_\lambda(\varepsilon))$};
			
		\end{tikzpicture}
  \end{subfigure}

    \caption{Proof of Theorem~\ref{thm:local-Blaschke}.}
\end{figure}

\begin{definition} For $x\in\Sigma$ and $\delta\in(0,\infty]$ we let
	\[
	U_\delta(x)\coloneqq \left\{ (y,v)\in N\Sigma\,|\, \dist(x,y)< \delta \text{ and } |v|<\delta \right\}\,.
	\]
	With $\mathcal D\subseteq TM$ the domain of the exponential map, we define the \emph{normal injectivity radius of $\Sigma$ at $x$} by
	\[
	\ninj^{\Sigma}(x)\coloneqq \sup\{ \delta > 0\,|\, U_\delta(x)\subseteq \mathcal D \text{ and } \exp_{|U_\delta(x)} \text{ is a diffeomorphism onto its image} \}\,.
	\]
\end{definition}

Our main result depends on the convexity radius of a not necessarily complete Riemannian manifold $(M,g)$,  see \cite[p.~158]{nardmann}:
\[
	\conv(x) = \sup\left\{\varrho\in[0,\inj(x)]\,|\,\forall r\in [0,\varrho): \Bsf_r(x) \text{ is strongly convex}\right\}\,.
\]
Note that $\mathrm{conv}_g(x)>0$ for all $x\in M$ by \cite[Corollary~3.5]{nardmann}.\\

We also recall the following definitions. Let $q\in\Sigma$, $v\in N_q\Sigma$ with $|v|=1$ and denote by $\gamma = \gamma_{qv}$ a corresponding unit-speed geodesic segment. A vector field $J$ along $\gamma$ is called a \emph{$\Sigma$-Jacobi field} if it is pointwise orthogonal to $\gamma'$ and satisfies the Jacobi equation
\[
	\frac{\nabla^2}{\dd t^2}J + \mathrm{Riem}(J,\gamma')\gamma' = 0
\]
along with the boundary conditions
\begin{align*}
	J(0)\in T_q\Sigma \qquad \text{ and }\qquad \frac{\nabla}{\dd t}J(0) + W_{\gamma'(0)}J(0) = 0\,,
\end{align*}
where $W_v$ is the Weingarten map associated with $\Sigma$, i.e.\  $(W_vX,Y) = (\secfund(X,Y),v) = \scalsecfund^v(X,Y)$. A simple calculation shows that the $\Sigma$-Jacobi fields are precisely the variation vector fields of variations of $\gamma$ through geodesics which emanate perpendicularly from $\Sigma$. A point $\gamma(b)$ is called a \emph{focal point of $\Sigma$} if there exists a nontrivial $\Sigma$-Jacobi field along $\gamma$ with $J(b) = 0$. The focal points are precisely the critical values of the normal exponential map $\exp|\mathcal D\cap N\Sigma$. The \emph{focal radius of $\Sigma$ at $q$} is the number
\[
	\inf\left\{b>0\,|\, \gamma_{qv}(b) \text{ is a focal point}, v\in N_q\Sigma,|v|_g=1 \right\}\,.
\]
The focal radius of $\Ssf_\lambda^c\subseteq \mathbb{M}^n(c)$ is precisely $\compR^c_\lambda$, cf.\ \cite[p.\ 391]{Cha06}.

\begin{theorem}[Klingenberg's lemma for embedded hypersurfaces]\label{thm:Klingenberg-hypersurface}
	Let $(M,g)$ be a Riemannian manifold, $\Sigma\subseteq M$ an embedded hypersurface and $x\in\Sigma$. Let $c,\lambda \geq 0$ and $s > 0$ such that
	\begin{itemize}
		\item $\Bsf_s(x)$ is a precompact, regular $\Sigma$-slice ball,
		\item $| \seccurv(\sigma) | \leq c $ for all $\sigma\in \mathrm{Gr}_2(T_yM), y\in \Bsf_s(x)$,
		\item $|\secfund^\Sigma| \leq \lambda$ on $\Sigma_s(x)$.
	\end{itemize}
	Then
	\begin{equation}\label{eqn:pointwise-lower-estimate-ninj}
		\ninj^{\Sigma}(x)\geq \frac{1}{6}\min\left\{s, \conv(x), \tfrac 65\compR_\lambda^c \right\}\,.
	\end{equation}
\end{theorem}
\begin{proof}
	Set $r\coloneqq \min \{s, \conv(x), \tfrac 65\compR_\lambda^c\}$. Then $\Bsf_{r}(x)$ is a precompact,  strongly convex regular $\Sigma$-slice ball with $|\seccurv|\leq c$ on $\Bsf_{r}(x)$, $|\secfund| \leq \lambda$ on $\Sigma_r$ and $5/6\cdot r\leq \compR^c_\lambda\leq \compR_{-\lambda}^{-c}$. For each $y\in \Bsf_{r/6}(x)$, $\overline\Bsf_{r/6}(y)\subseteq \Bsf_r(x)$ is compact, so that $U_{r/6}(x)$ is in the domain of $\exp$. It remains to show that $\exp$ is regular and injective on $U_{r/6}(x)$.\\
	
	Since $|\secfund|\leq \lambda$ on $\Sigma_{r}$, we have $\scalsecfund^{\Sigma_r,\nu}\leq \lambda$ for any choice of unit normal $\nu$ along $\Sigma_{r/6}$. Along with the assumption $\seccurv \leq c$ on $\exp(U_{r/6}(x))$, we deduce from \cite[Corollary~4.2(a)]{warn} that the focal radius of $\Sigma_{r/3}(x)$ is $\geq \compR_\lambda^c$, where we note that Warner defines the Weingarten map as minus our Weingarten map, which implies that his condition $S_{\gamma'(0)}\geq \delta\coloneqq -\lambda$ has to be read as $W_{\gamma'(0)}\leq \lambda$.\\
	
	To show that $\exp|U_{r/6}(x)$ is injective, let $(o,v),(q,w)\in U_{r/6}(x)$ be such that $\exp_o(v) = \exp_q(w)$. First assume that $v$ and $w$ point to the same side of $\Sigma_{r}$. Then by Theorem~\ref{thm:local-Blaschke} we have $(o,v)=(q,w)$. Next, assume that $v$ and $w$ point to different sides of $\Sigma_{r/6}$. In order for $\gamma_{o v}(1) = \gamma_{q w}(1)$ to hold, one of the two geodesics has to leave the component of $\Bsf_{r/3}\setminus\Sigma_{r/3}$ its initial velocity points to, say $\gamma_{qw}$, so that there exists $t_0\in (0,1]$ such that $\gamma_{qw}(t_0)\in\Sigma_{r/3}$. Again, this contradicts Theorem~\ref{thm:local-Blaschke}. 
\end{proof}

\section{Intersections of hypersurfaces and convex bodies}\label{sec:intersection-hypersurfaces-convex-bodies}

In this section, we first recall recent work of Drach, \cite{Dra24}, for $\lambda$-convex domains, see the Definition below. One of his main results is an \emph{outer Blaschke ball rolling theorem}, see Theorem~\ref{thm:blaschke-outer-ball-rolling-theorem}, for such domains. Using a variant of his proof, we then establish what we call, for lack of a better name, the \emph{Blaschke intersection ball Theorem}~\ref{thm:blaschke-intersection-ball-theorem}. At last, we prove the  main theorem of this section, Theorem~\ref{thm:normal-radius-est-sigma-intersection-convex-domain}, in which we show that if a properly embedded hypersurface intersects a $\lambda$-convex domain in two components, then there exists a geodesic sphere of controlled radius touching both of them.

\begin{definition}
    Let $\lambda>0$ and $D\subseteq M$ be a domain, i.e., a connected open set. Assume that $\overline{D}$ is compact, that the boundary $\partial D = \Gamma\neq\emptyset$ is smooth and bounds $D$. Endow $\Gamma$ with the inward, i.e. towards $D$, pointing unit normal vector field. Then we call $\Gamma$, $D$ and $\overline{D}$ \emph{$\lambda$-convex}, if the scalar second fundamental form of $\Gamma$ satisfies
    \[
        \lambda\leq \scalsecfund^\Gamma\,.
    \]
\end{definition}

\begin{remark}
    \begin{enumerate}[(i)]
        \item By \cite[Theorem~1.3]{bcgs11} a $\lambda$-convex domain $D$ is (geodesically) convex. In particular, if $D\ssubset U$ for a strongly convex open set $U$, then $D$ is strongly convex and its boundary $\partial D = \Gamma$ is diffeomorphic to $S^n$. Hence, $U$ is a $\Gamma$-slice set and a radial angle function domain for every $p\in D$.
        \item By (i), given $p\in D\ssubset U$, the corresponding radial angel function $\phi_p:\Gamma_U=\Gamma\to [0,\pi]$ has image in $[0,\pi/2)$.
    \end{enumerate}      
\end{remark}

Similar to the strategy we followed in sections~\ref{sec:radial-angle-comp} and \ref{sec:local-blaschke-and-klingenberg}, the following one-sided radial angle comparison theorem is the main technical ingredient in the proof of the outer Blaschke ball rolling theorem below. Note, however, that no localization is necessary, i.e.\ the radial angle comparison holds for all points in $\Gamma$.

\begin{theorem}[upper radial angle comparison for $\lambda$-convex domains]\label{thm:radial-angle-comp-lambda-convex-domains}
    Let $\Bsf_r(x)\subseteq M$ be a strongly convex ball and $\overline D\subseteq \overline\Bsf_{r/2}(x)$ a $\lambda$-convex domain. Assume that there exists $c\in\R$ such that
    \begin{equation}\label{eqn:convex-ball-lower-curvature}
        \phantom{\text{for all} \sigma\in \mathrm{Gr}_2(T_yM), y\in \Bsf_r(x)} c\leq\sec(\sigma)\quad\text{ for all }\quad \sigma\in \mathrm{Gr}_2(T_yM), y\in \Bsf_r(x)
    \end{equation}
    and
    \begin{equation}\label{eqn:comp-radius-less-equal-r4}
        \compR^c_\lambda \leq r/4\,. 
    \end{equation}    
    Let $s\in\Gamma = \partial D$ and $p\in D$ such that $\dist(p,s) = \dist(p,\Sigma)$. Let $p_\lambda\in\mathbb{M}^2(c)$ be on the inside of $\Ssf^c_\lambda\subseteq \mathbb{M}^2(c)$ and such that $\dist(p_\lambda,\Ssf^c_\lambda) = \dist(p,\Gamma)$ and denote the associated radial angle function with $\phi_\lambda = \phi_{p_\lambda}:\Ssf^c_\lambda\to[0,\pi]$. Then, the radial angle function $\phi_p:\Gamma\to[0,\pi]$ with respect to the origin $p$ satisfies for all $q\in\Gamma$, $q_\lambda\in\Ssf^c_\lambda$ such that $\dist(p,q)=\dist(p_\lambda,q_\lambda)$:
    \begin{equation}
        \phi_p(q)\leq \phi_\lambda(q_\lambda)\,.
    \end{equation}
\end{theorem}
\begin{proof}
    This is (the main part of) Theorem~3.1 in \cite{Dra24} applied to $D$ with $\Gamma=\partial D$. We note that the global assumption, in their notation, $\inj(M) > 2 R_\lambda$ is not used. Rather, they use $\inj_p(M) > 2R_\lambda=2\compR^c_\lambda$. This follows from the assumptions of the theorem. 
\end{proof}

\begin{theorem}[Blaschke outer ball rolling theorem]\label{thm:blaschke-outer-ball-rolling-theorem}
    Let $\Bsf_r(x)\subseteq M$ be a strongly convex ball and $\overline D\subseteq \overline\Bsf_{r/2}(x)$ a $\lambda$-convex domain. Assume that there exists $c\in\R$ such that \eqref{eqn:convex-ball-lower-curvature} and \eqref{eqn:comp-radius-less-equal-r4} hold. Then, for each $q\in\Gamma=\partial D$ we have
    \[
        \overline{D}\subseteq \overline{\Bsf}_{\compR^c_\lambda}(o)\,,
    \]
    where $o\coloneqq \gamma_{qv}(\compR^c_\lambda)$, i.e.\ the $\lambda$-convex domain $D$ can roll freely inside a ball of radius $\compR^c_\lambda$.
\end{theorem}
\begin{proof}
    This is Theorem~A part I in \cite{Dra24}, where we note that the use of Toponogov's hinge comparison in its proof, for which one presupposes completeness of $(M,g)$, can be replaced by the local hinge comparison with lower curvature bounds in Remark~\ref{rem:Rauch-triangle-comp}$(\geq)$.
\end{proof}

\begin{corollary}[diameter estimate for $\lambda$-convex domains]\label{cor:diamter-lambda-convex-domains}
    In the setting of Theorem~\ref{thm:blaschke-outer-ball-rolling-theorem}, one has
    \[
        \operatorname{diam}(D)\leq 2\compR^c_\lambda\,.
    \]
\end{corollary}

Note that, in the following theorem, if one had equality in \eqref{eqn:angle-inequality}, one would be in the situation of Theorem~\ref{thm:blaschke-outer-ball-rolling-theorem}. The theorem does not cover this case as it actually depends on it in the form of the last corollary. 

\begin{theorem}[Blaschke intersection ball theorem for $\lambda$-convex domains]\label{thm:blaschke-intersection-ball-theorem}
    Let $\Bsf_r(x)\subseteq M$ be a strongly convex ball and $\overline D\subseteq \overline\Bsf_{r/2}(x)$ a $\lambda$-convex domain. Assume that there exists $c\in\R$ such that \eqref{eqn:convex-ball-lower-curvature} and \eqref{eqn:comp-radius-less-equal-r4} hold. Let $q\in\Gamma$ and $\mu\in T_q\overline{D}$ be a unit vector pointing into $D$, i.e. $\sphericalangle(\nu_q,\mu)<\pi/2$, where $\nu$ is the inward pointing unit normal of $\Gamma$. Then, with $o\coloneqq \gamma_{q\mu}(\compR^c_\lambda)$, for each $s\in\Gamma$ such that
    \begin{equation}\label{eqn:angle-inequality}
        {\sphericalangle}(\mu,\dot\gamma_{qs}(0)) < \sphericalangle(\nu_q,\dot\gamma_{qs}(0))\,,
    \end{equation}
    one has $s\in\overline{\Bsf}_{\compR^c_\lambda}(o)$.
\end{theorem}
\begin{proof}
    Let $s\in\Gamma$ such that \eqref{eqn:angle-inequality} holds. Let $\varepsilon>0$ be so small that $p\coloneqq \gamma_{s\nu_s}(\varepsilon)$ satisfies $\dist(p,s) = \dist(p,\Gamma) = \varepsilon$. By decreasing $\varepsilon$ if necessary, we can assume without loss of generality that the radial angle function with origin $p$ satisfies $\phi_{p}(q)=\sphericalangle(\dot\gamma_{pq}(1),-\nu_q)=\sphericalangle(\dot\gamma_{qp}(0),\nu_q)>\sphericalangle(\dot\gamma_{qp}(0),\mu)$.

    Let $\Ssf_\lambda^c\subseteq \mathbb{M}^2(c)$ be a geodesic sphere of radius $\compR^c_\lambda$ and denote its centre by $o_\lambda$. Pick a point $s_\lambda\in \Ssf^c_\lambda$ and let $p_\lambda$ be a point on the inside of $\Ssf^c_\lambda$ such that $\dist(p_\lambda,\Ssf^c_\lambda) = \dist(p_\lambda,s_\lambda) = \varepsilon$. We also pick a point $q_\lambda\in\Ssf^c_\lambda$ such that $\dist(p_\lambda,q_\lambda) = \dist(p,q)$, noting that, according to the last corollary, $\dist(p,q)< \operatorname{diam}(D)\leq 2\compR^c_\lambda$.

    The hinge formed by the three points $p, q, o$ has side lengths identical to the ones of the hinge formed by $p_\lambda, q_\lambda, o_\lambda$ and, by assumption and Theorem~\ref{thm:radial-angle-comp-lambda-convex-domains}, $\sphericalangle p q o < \sphericalangle p_\lambda q_\lambda o_\lambda$, so that that the hinge comparison with lower curvature bounds, Remark~\ref{rem:Rauch-triangle-comp}$(\geq)$, yields $\dist(p,o)\leq\dist(p_\lambda,o_\lambda) = \compR^c_\lambda - \varepsilon$. Since this is true for every $\varepsilon > 0$ subject to the above conditions, we have $\dist(s,o)\leq \compR^c_\lambda$ as claimed.    
\end{proof}

\begin{figure}\label{fig:outer-intersection-blaschke}
    \begin{subfigure}[b]{0.45\textwidth}
    		\begin{tikzpicture}
                \draw[black] (0,0) circle (60pt) node at (2,2) {$\overline{\Bsf}_{\compR^c_\lambda}(o)$};
                \begin{scope}[rotate = 30]
                    \draw[black!60!green, fill=green, fill opacity=0.1] (0.175,0) ellipse (55pt and 40pt);   
                \end{scope}
                \draw[black!60!green] node at (0.5,-1.7) {$\Gamma$};
                \draw[black!60!green] node at (1.3,0) {$D$};
                \filldraw[black] (1.85,1) circle (1pt) node[right] {$q$};
                \draw (1.85,1)  to [out=200, in=40] (0,0);
                \filldraw[black] (0,0) circle (1pt) node[left] {$o$};
    	\end{tikzpicture}
     \end{subfigure}
     \begin{subfigure}[b]{0.45\textwidth}
		\begin{tikzpicture}
		\draw[black!60!green, fill=green, fill opacity=0.1] (0,0) ellipse (70pt and 50pt);   
            \draw[black!60!green] node at (-2,1.5) {$\Gamma$};
            \draw[black!60!green] node at (0,1) {$D$};
            \draw[smooth,tension=.7, blue] plot coordinates {(-4,1.3) (-3,1) (0,-2)  (2.5,1.5)  (3.5,2)} node at (-4,1) {$\Sigma$};
            \filldraw[black] (-0.8,-1.655) circle (1pt) node[below] {$q$};
            \draw (-0.8,-1.655) to [out=47,in=220] (0.5,-0.1);
            \filldraw[black] (0.5,-0.1) circle (1pt) node[right] {$o$};
            \draw[->,orange] (-0.8,-1.655) to (-0.5,-1.3) node[right] {$\mu$};
            \draw[->] (-0.8,-1.655) to (-0.7,-1.15) node[left] {$\nu$};
            \filldraw[black] (0.555,-1.72) circle (1pt) node[below] {$s$};
            \draw (-0.8,-1.655) to (0.555,-1.72);
            \draw[black] (0.5,-0.1) circle (58pt) node at (2,2) {$\overline{\Bsf}_{\compR^c_\lambda}(o)$};
    \end{tikzpicture}
  \end{subfigure}

    \caption{The left-hand figure illustrates Theorem~\ref{thm:blaschke-outer-ball-rolling-theorem}. The right-hand figure illustrates the application of Theorem~\ref{thm:blaschke-intersection-ball-theorem} as in case 3 in the proof of Theorem~\ref{thm:normal-radius-est-sigma-intersection-convex-domain}.}
\end{figure}

\begin{lemma}
    Let $\Sigma\subseteq M$ be a properly embedded hypersurface, $x\in M$ and $r>0$ such that $r<\inj(x)$. Assume $\Sigma\cap\Bsf_r(x)\neq\emptyset$. Then there exists $y\in \Sigma\cap\Bsf_r(x)$ such that $\dist(x,\Sigma) = \dist(x,y)$.
\end{lemma}
\begin{proof}
    By assumption we have $\dist(x,\Sigma)<r$. Let $(y_n)_{n\in\N}$ be a sequence in $\Sigma\cap\Bsf_r(x)$ with $\lim_{n}\dist(x,y_n)=\dist(x,\Sigma)$. By compactness of $\overline{\Bsf}_r(x)$ we can extract a subsequence, still denoted by $(y_n)_n$, which converges to a $y$ in $\Bsf_r(x)$. Since $\Sigma$ is properly embedded, $y\in\Sigma$.
\end{proof}

\begin{lemma}
Let $\Sigma\subseteq M$ be a properly embedded hypersurface and $p\in\Sigma$. Then there exists an $\varepsilon_0>0$ such that for all $0<\varepsilon<\varepsilon_0$, any geodesic sphere $\Ssf_\varepsilon(q)$ touching $\Sigma$ in $p$ satisfies $\Ssf_\varepsilon(q)\cap \Sigma = \{p\}$.
\end{lemma}
\begin{proof}
    Let $U$ be a tubular neighbourhood of $\Sigma$ of size $\delta:\Sigma\to(0,\infty)$ and define
    \begin{equation*}
        \varepsilon_0\coloneqq \frac12\min\{\conv(p),\sup\{\delta>0\,|\,\exp(U_\delta(p))\subseteq U\}\}\,.
    \end{equation*}
    Let $0<\varepsilon<\varepsilon_0$ and $\Ssf_\varepsilon(q)$ be a geodesic sphere touching $\Sigma$ at $p$. Note that $\Ssf_\varepsilon(q)$ is contained in $U$. Assume there exists $s\in\Sigma\cap\Ssf_\varepsilon(q)\setminus\{p\}$. If $\Ssf_\varepsilon(q)$ touches $\Sigma$ at $s$, then there are two distinct normal geodesic segments of length $\varepsilon$ emanating from $\Sigma$ and meeting at $q$, which is a contradiction to $\Ssf_\varepsilon(q)\subseteq U$. If there is no point in $\Sigma\cap\Ssf_\varepsilon(q)\setminus\{p\}$ at which $\Sigma$ and $\Ssf_\varepsilon(q)$ touch, then $\Bsf_\varepsilon(q)\cap\Sigma\neq\emptyset$. By the preceeding lemma, we obtain two distinct normal geodesic segments of length $<\varepsilon_0$ which emanate from $\Sigma$ and meet at $q$, which is, again, a contradiction to $\Ssf_\varepsilon(q)\subseteq U$.
\end{proof}

The following proposition represents what we call the \emph{moving ball technique}. Given a closed geodesic ball $\overline{\Bsf}_R(q)$ which touches the hypersurface $\Sigma$ in a point $p$ and which has at least one more point in common with $\Sigma$, it produces a geodesic sphere $\Ssf_{r_0}(o)$ which touches $\Sigma$ in $p$ and a second point $s$. The proof works by moving the centre $q$ of $\overline{\Bsf}_R(q)$ towards $p$, thereby decreasing the radius, until it cannot be decreased anymore without sharing any points with $\Sigma$ other than $p$.

\begin{proposition}\label{prop:geodesic-sphere-touching-Sigma-two-points}
    Let $\Sigma\subseteq M$ be a properly embedded hypersurface and $p\in\Sigma$. Let $\Ssf_R(q)$ be a geodesic sphere touching $\Sigma$ at $p$ with $R<\conv(q)$ and satisfying $\overline\Bsf_R(q)\cap\Sigma\setminus\{p\}\neq\emptyset$. Then there exists $r_0\in(0,R]$ and a geodesic sphere $\Ssf_{r_0}(o)$ touching $\Sigma$ at $p$ and a point $s\in\Sigma\setminus\{p\}$.
\end{proposition}
\begin{proof}
    Let $\sigma:[0,R]\to M$ be the unit-speed geodesic segment with $\sigma(0)=p$ and $\sigma(R)=q$. By the preceeding lemma and the assumption $\overline\Bsf_R(q)\cap\Sigma\setminus\{p\}\neq\emptyset$, we have
    \[
        r_0\coloneqq \inf\{ t\in (0,R]\,|\, \overline\Bsf_t(\sigma(t))\cap\Sigma\setminus\{p\} \neq \emptyset  \} > 0\,.
    \]
    Thus, there exists a point $s\in\Ssf_{r_0}(\sigma(r_0))\cap\Sigma\setminus\{p\}$. We claim that $\Ssf_{r_0}(\sigma(r_0))$ and $\Sigma$ must necessarily touch at $s$. Suppose this was not the case. Then $\dd(\dist_{\sigma(r_0)})_s\neq 0$ and $\dd(\dist_s)_{\sigma(r_0)}(\dot\sigma(r_0))<1$, so that there exist $\varepsilon > 0$ and an open neighbourhood $U\subseteq \Sigma$ of $s$ such that for all $t\in (r_0-\varepsilon,r_0+\varepsilon)$ and $x\in U$ one has $\omega^t_x\coloneqq \dd(\dist_{\sigma(t)})_x\neq 0$ and $f^t(x)\coloneqq \dd(\dist_x)_{\sigma(t)}(\dot\sigma(t)) - 1 \neq  0$. Define a time-dependent vector field $X$ by
    \[
        \phantom{\qquad t\in (t_0-\varepsilon,t_0+\varepsilon), x\in U\,.}   X^t_x\coloneqq -\frac{f^t(x)}{|\omega^t_x|^2}(\omega^t_x)^\sharp\,,\qquad t\in (r_0-\varepsilon,r_0+\varepsilon), x\in U\,,
    \]
    and let $\gamma:I\to \Sigma$ be a solution to the nonautonomous ODE $\dot\gamma(t) = X^t_{\gamma(t)}$ with initial condition $\gamma(r_0) = s$, where $I$ is a neighbourhood of $t_0$. Then, for $t\in I$, we have
    \begin{align*}
        \dist(\sigma(t),\gamma(t)) &= r_0+\int_{r_0}^t\!\frac{\dd}{\dd s}\dist(\sigma(s),\gamma(s))\,\dd s = r_0+\int_{r_0}^t\!\dd\dist_{(\sigma(s),\gamma(s))}(\dot\sigma(s),\dot\gamma(s))\,\dd s\\
        &=r_0+\int_{r_0}^t\!\left(\dd (\dist_{\gamma(s)})_{\sigma(s)}(\dot\sigma(s))+ \dd(\dist_{\sigma(s)})_{\gamma(s)}(\dot\gamma(s))\right)\,\dd s\\
        &= r_0+\int_{r_0}^t\!\left(1+f^s(\gamma(s))+\omega^s_{\gamma(s)}(\dot\gamma(s))\right)\,\dd s = t\,,        
    \end{align*}
    i.e. $\gamma(t)\in\Sigma\cap \Ssf_t(\sigma(t))$, which is a contradiction to the definition of $r_0$. Hence, $\Ssf_{r_0}(\sigma(r_0))$ and $\Sigma$ do indeed touch at $s$.    
\end{proof}

\begin{theorem}\label{thm:normal-radius-est-sigma-intersection-convex-domain}
    Let $\Sigma$ be a properly embedded hypersurface in $(M,g)$, $x\in M$, $r \leq \conv(x)$ and $\overline D\subseteq \overline\Bsf_{r/2}(x)$ be a $\lambda$-convex domain. Assume that there exists $c\in\R$ such that \eqref{eqn:convex-ball-lower-curvature} and \eqref{eqn:comp-radius-less-equal-r4} hold. Then, if $\Sigma\cap D\neq \emptyset$ and  $\Sigma$ touches $\Gamma=\partial D$ or $\Sigma\cap D$ has two components, there exists a geodesic sphere $\Ssf_\rho(p)\subseteq\Bsf_r(x)$ with $\rho\leq \compR^c_\lambda$ touching $\Sigma$ at two distinct points in $\Sigma\cap\Bsf_r(x)$.
\end{theorem}
\begin{proof}
    We divide the proof into three cases.\\
    \noindent
    \emph{Case 1:} We assume $\Sigma$ touches $\Gamma$ at a point $q$. By Theorem~\ref{thm:blaschke-outer-ball-rolling-theorem}, there exists a geodesic sphere $\Ssf_{\compR^c_\lambda}(o)$ touching $\Gamma$ and thus $\Sigma$ in $q$ and such that $\overline D\subseteq \overline{\Bsf}_{\compR^c_\lambda}(o)$. By Proposition~\ref{prop:geodesic-sphere-touching-Sigma-two-points}, we find a geodesic sphere $\Ssf_\rho(o)$ with $s\leq \compR^c_\lambda$ touching $\Sigma$ in $q$ and a second point $s\neq q$.\\
    
    \noindent
    \emph{Case 2:} We assume that $\Sigma\cap D$ has two components, say $\Sigma_0$ and $\Sigma_1$, and that there exists a geodesic $\gamma_{qp}$ of length $d=\dist(\Sigma_0,\Sigma_1)$ between points $q\in\overline{\Sigma_0}$ and $p\in\overline{\Sigma_1}$ which is normal to $\Sigma$ at, say, $q$. Note that we can assume $d>0$ for otherwise we are in case 1. By Corollary~\ref{cor:diamter-lambda-convex-domains}, $d \leq 2\compR^c_\lambda$. The geodesic sphere $\Ssf_{d/2}(\gamma_{qp}(1/2))$ touches $\Sigma$ in $q$ and $p\in\overline{\Bsf}_{d/2}(\gamma_{qp}(1/2))$ so that Proposition~\ref{prop:geodesic-sphere-touching-Sigma-two-points} yields a geodesic sphere as claimed.\\
    
    \noindent
    \emph{Case 3:} Let $\Sigma_0$ and $\Sigma_1$ be components of $\Sigma\cap D$ and $d=\dist(\Sigma_0,\Sigma_1)>0$. Let $q\in\partial\Sigma_0$ and $s\in\partial\Sigma_1$ such that $\dist(q,s) = d$ and assume that the corresponding geodesic $\gamma_{qs}$ is not normal to $\Sigma$ at either points.
    
    By case 1, we can assume that $\Sigma$ and $\Gamma$ touch at no point of their intersection, i.e.\ their intersection is transversal, so that $\partial\Sigma_0$ and $\partial\Sigma_1$ are closed embedded codimension two submanifolds of $M$. In particular, the geodesic segment $\gamma_{qs}$ realises the distance between $\partial\Sigma_0$ and $\partial\Sigma_1$ so that $\dot\gamma_{qs}(0)\in N_q\partial\Sigma_0$. Denote by $\nu$ the inward pointing unit normal vector field of $\Gamma$ and let $\mu\in N_q\Sigma$ be the unit vector pointing into $D$. Note that $\mu\in N_q\partial\Sigma_0$ and that $\mu$ lies in the open sector $\{\alpha\cdot\nu_q+\beta\cdot\dot\gamma_{qs}(0)\,|\,\alpha,\beta >0\}$. Indeed, if $\mu$ conincided with $\nu_q$ or was a multiple of $\dot\gamma_{qs}(0)$, we were in case 1 respectively case 2, and if $\mu$ was outside the closure of the sector we could find a geodesic segment between $\Sigma_0$ and $\overline{\Sigma_1}$ of length $< d$. It follows in particular that $\sphericalangle(\dot\gamma_{qs}(0),\mu)\leq \sphericalangle(\dot\gamma_{qs}(0),\nu_q)$. Theorem~\ref{thm:blaschke-intersection-ball-theorem} yields a geodesic sphere $\Ssf_{\compR^c_\lambda}(o)$ touching $\Sigma$ at $q$ and such that $s\in \overline \Bsf_{\compR^c_\lambda}(o)$. Proposition~\ref{prop:geodesic-sphere-touching-Sigma-two-points} then yields a geodesic sphere as claimed.    
\end{proof}

\begin{corollary}\label{cor:lambda-convex-body-in-Sigma-tubular-neighbourhood}
    Let $x\in\Sigma$ and $\delta>0$ such that $\exp|U_\delta(x)$ is a diffeomorphism onto. Let $r\leq \min\{\conv(x),\delta/2\}$ and $\overline D\subseteq \overline\Bsf_{r/2}(x)$ be a $\lambda$-convex domain. Assume there exists $c\in\R$ such that \eqref{eqn:convex-ball-lower-curvature} and \eqref{eqn:comp-radius-less-equal-r4} hold. Then, $D\cap\Sigma$ is connected. In particular, $D$ is a $\Sigma$-slice domain.
\end{corollary}
\begin{proof}
    The ``in particular'' part follows from Remark~\ref{rem:continuity-seperation-property-and-slice-balls}(iii). Assuming that $D\cap \Sigma$ has two components, the last theorem asserts that there exists a geodesic sphere $\Ssf_\rho(o)\subseteq \Bsf_r(x)$ touching $\Sigma$ at two distinct points. Since $\rho\leq\compR^c_\lambda< r$ and $\Bsf_r(x)\subseteq \exp(U_{\delta}(x))$, this is a contradiction to $\exp|U_{\delta}(x)$ being a diffeo onto.   
\end{proof}

\section{Hypersurfaces as graphs}

We assume for the moment that $M=I\times N$, where $I$ is an open interval containing zero. For each $t\in I$ we write $M^t = \{t\}\times N$. We further assume that the metric $g$ decomposes as $g = \dd t^2+g_t$ with $(g_t)_{t\in I}$ a family of Riemannian metrics on $N$. This implies that the $\partial_t$-integral curves are geodesics and that each such geodesic meets each $M^t$ perpendicularly. By $\scalsecfund^t$ we denote the, with respect to $\partial_t$, scalar-valued second fundamental form of $M^t$ and note in passing that $\scalsecfund^t(X,Y) = -\tfrac12 \partial_t g_t(X,Y)$.\\

Let $\Sigma\subseteq M$ be a properly embedded hypersurface. We assume that there exists a function $f:N\to\R$ such that $\Sigma = \Gamma(f) = \left\{(f(x),x)\,|\,x\in N\right\}$. The tangent spaces to $\Sigma$ are given by $T_{(f(x),x)}\Sigma = \{\dd f_x(v)\cdot \partial_t(f(x),x) + v\,|\,v\in T_xN\}$ and a unit normal to $\Sigma$ is given by
\begin{equation}\label{eqn:Sigma-normal}
	\nu_{(f(x),x)} = \frac{\partial_t(f(x),x)-\grad_{f(x)}f_x}{\sqrt{1 + |\grad_{f(x)}f_x|^2_{f(x)}}}\,,
\end{equation}
where we denote objects associated with the metric $g_t$ by $\bullet_t$, e.g.\ $\grad_{f(x)} = \grad_{g_{f(x)}}$. We will now calculate the scalar valued second fundamental form $\scalsecfund^\Sigma$ of $\Sigma$. To that end, denote for a point $x\in M^0$ by $\underline{x}$ the point $(f(x),x)$ and $\underline{X} = \dd f X+X$ for a tangent vector $X$ of $N$.

\begin{lemma}  For $\underline{X},\underline{Y}\in T_{\underline{x}}\Sigma$ one has
    {\small
    \begin{multline}\label{eqn:second-fundamental-form-Sigma}
		\scalsecfund^\Sigma(\underline{X},\underline{Y}) = \\ \frac{1}{\sqrt{1 + |\grad_{f(x)}f|^2_{f(x)}}}\left(\Hess_{f(x)}f(X,Y) + \scalsecfund^{f(x)}(X,Y) + \scalsecfund^{f(x)}\left(\grad_{f(x)}f,X(f)Y+Y(f)X\right)\right)\,.
	\end{multline}	}
\end{lemma}
\begin{proof}
	Let $X,Y$ be vector fields on $N$ which we extend to $M$ by requiring them to be constant along $\partial_t$-lines. Similarly, we extend $f$ to $M$. Then
	\begin{align*}
		&\scalsecfund^\Sigma(\underline{X},\underline{Y}) = g(\nu_{\underline{x}},\nabla_{\underline{X}}\underline{Y}) = g(\nu_{\underline{x}},\nabla_{X(f)\partial_t+X}(Y(f)\partial_t+Y)) \\
		&=g(\nu_{\underline{x}}, X(f)\nabla_{\partial_t}(Y(f)\partial_t+Y) + \nabla_X(Y(f)\partial_t+Y)) \\
		&=g(\nu_{\underline{x}}, X(f)\partial_t(Y(f))\partial_t + X(f)Y(f)\nabla_{\partial_t}\partial_t + X(f)\nabla_{\partial_t}Y   + X(Y(f))\partial_t + Y(f)\nabla_X\partial_t + \nabla_XY)\,. \\
	\end{align*}
	Note that $\partial_t(Y(f))=0$ since $f$ is constant in $\partial_t$-direction and that $\nabla_{\partial_t}\partial_t = 0$ since the $\partial_t$-integral curves are geodesics. Hence, $\scalsecfund^\Sigma(\underline{X},\underline{Y})$ is equal to
	\begin{multline*}
		g(\nu_{\underline{x}}, X(f)\nabla_{\partial_t}Y   + X(Y(f))\partial_t + Y(f)\nabla_X\partial_t + \nabla_XY)\\
		=\frac{1}{\sqrt{1 + |\grad_{f(x)}f|^2_{f(x)}}}g(\partial_t-\grad_{f(x)}f,  X(f)\nabla_{\partial_t}Y   + X(Y(f))\partial_t + Y(f)\nabla_X\partial_t + \nabla_XY)\\
		=\frac{1}{\sqrt{1 + |\grad_{f(x)}f|^2_{f(x)}}}\Big(X(Y(f)) - g(\grad_{f(x)}f,\nabla_XY) + g(\partial_t,\nabla_XY) \\
		-g(\grad_{f(x)}f,X(f)\nabla_{\partial_t}Y+Y(f)\nabla_X\partial_t \Big)\,,
	\end{multline*}
	where we have used that $g(\partial_t,\nabla_{\partial_t}Y) = g(\partial_t,\nabla_X\partial_t) = 0$. Note that the second term in parantheses is equal to $g(\grad_{f(x)}f,\nabla^{f(x)}_XY) = \dd f(\nabla^{f(x)}_XY)$ and that the third term is equal to $\scalsecfund^{f(x)}(X,Y)$, so that the above is equal to
    \begin{small}
	\begin{equation*}
		\frac{1}{\sqrt{1 + |\grad_{f(x)}f|^2_{f(x)}}}\Big(\Hess_{f(x)}f(X,Y) + \scalsecfund^{f(x)}(X,Y)-g(\grad_{f(x)}f,X(f)\nabla_{\partial_t}Y+Y(f)\nabla_X\partial_t \Big)\,.
	\end{equation*}
    \end{small}
	Using the Koszul-formula, one easily computes $g(U,\nabla_{\partial_t}V)_{f(x)} = -\scalsecfund^{f(x)}(U,V)$, which gives the claimed formula.
\end{proof}

\begin{lemma}\label{lem:gradient-estimate}
	Assume there exist constants $C^\Sigma,C^F\geq 0$ with $C^\Sigma+C^F>0$ such that $|\scalsecfund^\Sigma|\leq C^\Sigma$ and $|\scalsecfund^t|_t\leq C^F$ for all $t\in I$. Let $x\in N$ and assume that $f(x)= 0$ and $\dd f_x = 0$, i.e.\ $\Sigma$ is tangent to $M^0$ at $(0,x)$. Let $\gamma:[0,L]\to N$ be a curve with $\gamma(0)=x$ and assume there exists $C^G>0$ such that $|\dot{\gamma}|_{t}\leq C^G$ for all $t\in I$. Then
	\begin{equation}\label{eqn:gradient-estimate}
		\phantom{\qquad \text{ for all } \qquad t\in [0,\alpha) \,,	}|\grad_{f(\gamma(s))}f_{\gamma(s)}|_{f(\gamma(s))}\leq \psi^{-1}(s)\qquad \text{ for all } \qquad s\in [0,\alpha) \,,
	\end{equation}
	where $\psi:[0,\infty)\to [0,\alpha)$ is given by
	\begin{equation}\label{eqn:psi-integral}
		\psi(s)\coloneqq \frac{1}{C^G}\int\limits_0^s\!\frac{1}{C^\Sigma(1+\sigma^2)^{3/2}+C^F(1+\sigma^2)} \,\dd\sigma\,.
	\end{equation}
\end{lemma}
\begin{proof}
	As before, we extend $f$ to all of $M$ such that it is constant along $\partial_t$-lines. With $\widetilde{\gamma}(s)\coloneqq (f(\gamma(s)),\gamma(s))$ we then have
	\[
		\grad_{f(\gamma(s))}f_{\gamma(s)} = \grad_gf_{\widetilde{\gamma}(s)}\,.
	\]
	To estimate the norm of this gradient, we compute the $s$-derivative of its square:
	\begin{align*}
		\frac{\dd}{\dd s} |\grad_gf_{\widetilde{\gamma}(s)}|^2_g &= \frac{\dd}{\dd s} g(\grad_gf_{\widetilde{\gamma}(s)},\grad_gf_{\widetilde{\gamma}(s)}) = 2g\left(\frac{\nabla^{\widetilde{\gamma}}}{\dd s}\grad_gf_{\widetilde{\gamma}(s)} , \grad_gf_{\widetilde{\gamma}(s)}\right)\\
		&= 2\Hess_gf(\dot{\widetilde{\gamma}}(s),\grad_gf_{\widetilde{\gamma}(s)})\\
		&= 2\Hess_gf(\dd f_{\gamma(s)}\dot{\gamma}(s)\partial_t,\grad_gf_{\widetilde{\gamma}(s)}) + 2\Hess_gf(\dot{\gamma}(s),\grad_gf_{\widetilde{\gamma}(s)})\,.
	\end{align*}
	Since $f$ is constant in $\partial_t$-directions, 
	\[
		\Hess_gf(\dot{\gamma}(s),\grad_gf_{\widetilde{\gamma}(s)}) = \Hess_{f(\gamma(s))}f(\dot{\gamma}(s), \grad_gf_{\widetilde{\gamma}}(s))
	\]
	by \eqref{eqn:covariant-derivative-submldf-0}. For $X$ a vetor field tangent to $M^t$ and constant in $\partial_t$-direction, we have \newline $\Hess_gf(\partial_t,X) = \partial_t(X(f)) - \dd f(\nabla_{\partial_t}X)$, the first term of which vanishes, so that
	\[
		\Hess_gf(\partial_t,X) = \scalsecfund^t(X,\grad_tf)\,.
	\]
	Hence, we obtain
	\begin{multline*}
		\frac{\dd}{\dd s} |\grad_gf_{\widetilde{\gamma}(s)}|_g = \frac12 |\grad_gf_{\widetilde{\gamma}(s)}|_g^{-1}\frac{\dd}{\dd s}|\grad_gf_{\widetilde{\gamma}(s)}|^2_g\\
		=  |\grad_gf_{\widetilde{\gamma}(s)}|_g^{-1}\big(\dd f_{\gamma(s)}\dot{\gamma}(s)\scalsecfund^{f(\gamma(s))}(\grad_gf_{\widetilde{\gamma}(s)}, \grad_gf_{\widetilde{\gamma}(s)}) \\+ \Hess_{f(\gamma(s))}f(\dot{\gamma}(s), \grad_gf_{\widetilde{\gamma}(s)})\big).
	\end{multline*}
	Using \eqref{eqn:second-fundamental-form-Sigma}, we express this in terms of the second fundamental forms of all $M^t$ and $\Sigma$ as
	\begin{multline*}
		\frac{\dd}{\dd s} |\grad_gf_{\widetilde{\gamma}(s)}|_g  = |\grad_gf_{\widetilde{\gamma}(s)}|_g^{-1}\cdot \\
		\Big( \sqrt{1+|\grad_gf_{\widetilde{\gamma}(s)}|_g^2}\scalsecfund^\Sigma(\dd f_{\gamma(s)}\dot{\gamma}(s)\partial_t+\dot{\gamma}(s), \dd f_{\widetilde{\gamma}(s)}(\grad_gf_{\widetilde{\gamma}(s)})\partial_t+\grad_gf_{\widetilde{\gamma}(s)})\\
		 - (1+|\grad_gf_{\widetilde{\gamma}(s)}|_g^2)\scalsecfund^{f(\gamma(s))}(\dot{\gamma}(s),\grad_gf_{\widetilde{\gamma}(s)})\Big)\,.
	\end{multline*}
	Using the assumptions on the boundedness of $\scalsecfund^\Sigma$, $\scalsecfund^t$ and $\dot{\gamma}$, this yields
	\begin{equation*}
		\frac{\dd}{\dd s} |\grad_gf_{\widetilde{\gamma}(s)}|_g\leq C^G\cdot\left(C^\Sigma (1+|\grad_gf_{\widetilde{\gamma}(s)}|^2_g)^{3/2} + C^F(1+|\grad_gf_{\widetilde{\gamma}(s)}|^2_g)\right)\,.
	\end{equation*}
	Integration leads to
	\begin{equation*}
		|\grad_gf_{\widetilde{\gamma}(s)}|_g\leq \int_0^s\!C^G\cdot\left(C^\Sigma (1+|\grad_gf_{\widetilde{\gamma}(\sigma)}|^2_g)^{3/2} + C^F(1+|\grad_gf_{\widetilde{\gamma}(\sigma)}|^2_g)\right)\,\dd\sigma\,,
	\end{equation*}
	to which we apply the Bihari-LaSalle-inequality \cite[\S 3]{Bihari}, \cite[Lemma~1]{LaSalle}, yielding the claimed estimate.  
\end{proof}

\begin{lemma}\label{lem:curv-par-hyp}
    With $I=(-a,a)$ and $M=I\times N$, assume there exist $c,\lambda\geq0$ such that $|\seccurv|\leq c$, $|\scalsecfund^0|\leq \lambda$ and $a\leq \compR^c_{-\lambda}$. Then, for all $t\in I$ one has
    \begin{equation}\label{eqn:parallel-hypersurface-curvatures}
        -\mathrm{ct}_c(\compR^c_{-\lambda}-|t|)\leq \scalsecfund^t\leq \begin{cases}
            \ct_{-c}(|t|+\compR^{-c}_\lambda) & \text{ if } \lambda > \sqrt{c}\,,\\
            \sqrt{c} & \text{ if } \lambda = \sqrt{c}\,,\\
             \mathrm{tn}_{-c}(|t|+\mathrm{atn}_{-c}(\lambda))&\text{ if } \lambda < \sqrt{c}\,,            
        \end{cases}
    \end{equation}
    where $\mathrm{ct}_c$ is the generalised cotangent \eqref{eqn:ctc}, $\compR^c_\lambda$ the radius of a geodesic sphere of curvature $\lambda$ in $\mathbb M^n(c)$, $\mathrm{tn}_{-c}=c\cdot \sn_{-c}/\sn_{-c}'$ is the generalised tangent and $\mathrm{atn}_{-c}=\mathrm{tn}_{-c}^{-1}$ the corresponding inverse.
\end{lemma}
\begin{remark}
    The three cases on the right-hand side are the principal curvatures of, respectively, geodesic spheres ($\lambda>\sqrt{c}$), horospheres ($\lambda=\sqrt{c}$) and equidistant hypersurfaces ($\lambda < \sqrt{c}$) in $\mathbb M^n(c)$.
\end{remark}
\begin{proof}
    This follows from \cite[Theorem~3.1]{Es87} by comparing $M^0$ with $\Ssf^c_{-\lambda}\subseteq \mathbb M^n(c)$ and $\Ssf^{-c}_\lambda\subseteq \mathbb{M}^n(-c)$, where we note that \cite[Theorem~3.1]{Es87} does not need the assumed completeness of $(M,g)$, but only the well-definedness of all geodesics emanating perpendicularly from $M^0$ up to, but not including, $|t|=\compR^c_{-\lambda}$.
\end{proof}

Recall that we use the maximum norm on the second fundamental form of hypersurfaces.

\begin{proposition}\label{prop:curvature-exp-hyperdisc}
    Let $(M,g)$ be a Riemannian manifold, $p\in M$ and $r\leq \inj^M_g(p)$. Let $V\subseteq T_pM$ be a vector subspace of codimension one and define an embedded hypersurface in $M$ by
    \[
        \Gamma\coloneqq \Gamma^{x,r,V} \coloneqq \exp_p(\Bsf_r(0)\cap V)\,.
    \]
    Assume there exist constants $c,\Lambda\geq 0$ such that
   \begin{equation*}
        \phantom{\quad \text{ on } \quad \Bsf_r(p)\,.} |\seccurv|\leq c \quad \text{ and } \quad \left|\nabla\mathrm{Riem}\right|\leq \Lambda \quad \text{ on } \quad \Bsf_r(p)
   \end{equation*}
   and $r< \pi/\sqrt{c}$.
    Then, for any point $x\in \Gamma\setminus\{p\}$ with $s=\dist(p,x)$ we have
    \begin{equation}\label{eqn:2nd-fundamental-form-estimate}
        |\secfund^\Gamma|_x\leq \frac{8}{9}\cdot s^2\cdot\frac{\mathrm{sn}_{-c}(s/2)^2}{\mathrm{sn}_{c}(s)^2}\left(3\cdot\Lambda\cdot\mathrm{sn}_{-c}(s/2)^2+4\cdot c\cdot\mathrm{sn}_{-c}(s)\right)\,,
    \end{equation}
    where $\mathrm{sn_{c}}$ denotes the generalised sine \eqref{eqn:snc}.
\end{proposition}
\begin{proof}
    By definition of $\Gamma$, we have $\secfund^\Gamma_p=0$. Moreover, for every $x\in\Gamma$, we have \[\secfund^\Gamma(\grad(\dist_p)_x, \grad(\dist_p)_x)=0\,,\]so that it suffices to bound $\secfund^\Gamma(X,X)$ for every $X\in T(\Gamma\setminus\{p\})$ with $X\perp\partial_\rho$.

    Let $x\in\Gamma\setminus\{p\}$, $X\in T_x\Gamma$ with $|X|=1$ and $X\perp\grad(\dist_p)$, and let $\sigma:I\to\Gamma$ be a $\Gamma$-geodesic with $\sigma(0)=x$ and $\dot\sigma(0)=X$ so that we have
    \[
        \secfund^\Gamma(X,X) = \left(\frac{\nabla}{\dd t}\dot\sigma(t)_{|t=0}\right)^{N\Gamma}\,,
    \]
    where ${\bullet}^{N\Gamma}$ denotes the projection onto the normal bundle $N\Gamma$ of $\Gamma$.
    With $\eta:I\to V$ defined by $\eta(t)\coloneqq \exp_p^{-1}(\sigma(t))$ we thus have
    \begin{align*}
        \secfund^\Gamma(X,X) &= \left(\frac{\nabla}{\dd t}\dd (\exp_p)_{\eta(t)}(\dot\eta(t))_{|t=0}\right)^{N\Gamma} \\&= \left(\Hess (\exp_p)_{\eta(0)}(\dot\eta(0),\dot\eta(0)) + \dd(\exp_p)_{\eta(0)}(\Ddot{\eta}(0))\right)^{N\Gamma}\\
        &=\Hess (\exp_p)_{\eta(0)}(\dot\eta(0),\dot\eta(0))^{N\Gamma}\,.
    \end{align*}
    Hence,
    \[
        \left|\secfund^\Gamma(X,X)\right|= \left| \Hess (\exp_p)_{\eta(0)}(\dot\eta(0),\dot\eta(0))^{N\Gamma} \right|\,.
    \]
    With $s=\dist(p,x)$ and $v\coloneqq \eta(0)/|\eta(0)|$ we have $sv=\eta(0)$ and $\dd (\exp_p)_{sv}(\dot\eta(0))=X$.
    Let $J$ be the Jacobi field along the geodesic $\gamma(t)\coloneqq \exp_p(tv)$ with $J(0)=0$ and $\tfrac{\nabla}{\dd t}J(0)=\dot\eta(0)$, so that $J(t)=\dd(\exp_p)_{tv}(t\dot\eta(0))$. From the Jacobi field comparison theorem with upper curvature bounds, see \cite[Theorem~11.9(a)]{lee}, we have
    \[
        s = s\cdot|X|=s|\dd(\exp_p)_{sv}(\dot\eta(0))| = |J(s)|\geq \mathrm{sn}_c(s) |\dot\eta(0)|\,,
    \]
    where $\mathrm{sn}_c$ is the generalised sine \eqref{eqn:snc}, which yields the bound
    \[
        |\dot\eta(0)|\leq \frac{s}{\mathrm{sn}_c(s)}\,.
    \]
    We plug this into the Hessian estimate \cite[Theorem~4.10]{lezcano}\footnote{It seems that the factor $1/r^2$ in the definition of $\rho$ in \cite[Theorem~4.10]{lezcano} is superfluous. Indeed, the vector field $K$ in the proof of loc.\ cit.\ contains an additional coefficient $1/r^2$ compared to its definition \cite[Proposition~4.1]{lezcano}.} to obtain
    \[
    \left| \Hess(\exp_p)_{\eta(0)}(\dot\eta(0), \dot\eta(0))^{N\dot\gamma}\right|\leq \frac{8}{9}\cdot s^2\cdot\frac{\mathrm{sn}_{-c}(s/2)^2}{\mathrm{sn}_{c}(s)^2}\left(3\cdot\Lambda\cdot\mathrm{sn}_{-c}(s/2)^2+4\cdot c\cdot\mathrm{sn}_{-c}(s)\right)\,,
    \]
    where $N\dot\gamma$ denotes the normal bundle of $\mathrm{im}\dot\gamma$. Since $N_x\Gamma\subseteq N_x\dot\gamma$, the proposition is proved.    
\end{proof}

\begin{definition}
    For a properly embedded hypersuface $\Sigma$ in $(M,g)$ and $x\in M$ we define the \emph{graphical radius of $\Sigma$ at $x$} to be
    \begin{align}\label{eqn:graphical-radius}
        \mathrm{graph}^\Sigma(x)\coloneqq \sup \>&\big\{r\in (0,\inj^M_g(x)]\,\big|\,r\leq \ninj^{\Gamma^{x,\inj(x),T_x\Sigma}}(x) \text{ and there exists } f\in \Csf^\infty(\Gamma^{x,r,T_x\Sigma})\text{ s.t. }\\\nonumber
        & -r<f<r \text{ and }  (\Sigma\cap\Bsf_r(x))^0 = \Gamma(f) \big\}\,,
    \end{align}
    where $(\Sigma\cap\Bsf_r(x))^0$ denotes the connected component of $\Sigma\cap\Bsf_r(x)$ containing $x$ and $\Gamma(f) \coloneqq \left\{\exp_y(f(y)\cdot\nu_y)\,|\,y\in \Gamma^{x,r,T_xM}\right\}$ is the graph of $f$ with $\nu$ a unit normal to $\Gamma^{x,r,T_x\Sigma}$.
\end{definition}

With the above definition, we finally arrive at the main result of this section:

\begin{theorem}\label{thm:graphing-radius-estimate}
    Let $\Sigma$ be a properly embedded hypersurface in a Riemannian manifold $(M,g)$ and $x\in\Sigma$. Assume there exist $r>0$ and $c,\Lambda\geq 0$ such that
    \begin{itemize}
        \item The ball $\Bsf_r(x)$ is precompact,
        \item $|\seccurv(\sigma)|\leq c$ for all $\sigma\in\mathrm{Gr}_2(T_yM)$, $y\in\Bsf_r(x)$,
        \item $|\nabla \mathrm{\mathrm{Riem}}|\leq \Lambda$ on $\Bsf_r(x)$,
        \item $|\secfund^\Sigma|\leq\lambda$ on $\Bsf_r(x)$.
    \end{itemize}
    Then, there exists an explicitely computable constant $C = C(r,c,\Lambda,\lambda,\inj(x),\conv(x))>0$ such that
    \[
        \mathrm{graph}^\Sigma(x)\geq C\,.
    \]
\end{theorem}
\begin{proof}
By the implicit function theorem and a unique continuation argument, it suffices to show $|f(s)|\leq s$ and the finiteness of $|\grad_{f(\gamma(s))}f_{\gamma(s)}|_{f(\gamma(s))}$ along unit speed, minimising geodesics $\gamma$ in $\Gamma^{x,s,T_x\Sigma}$ emanating from $x$, where $f$ is as in \eqref{eqn:graphical-radius} and $s$ is suitably chosen.

Let $r_0\coloneq \min\{r,\inj(x)\}$ and $\Gamma_{r_0}\coloneq \Gamma^{x,r_0,T_x\Sigma}$. Denote by $\rho(s)$ the right-hand side in \eqref{eqn:2nd-fundamental-form-estimate}. Let $r_1\in (0,\min\{r_0,\pi/\sqrt(c))$ be the argument at which
\[
    s\mapsto \min\{s,\tfrac 65\compR^c_{\rho(s)}\}
\]
attains its maximum. Noting that $\Bsf_{r_1}(x)$ is a $\Gamma_{r_1}$-slice ball, by Theorem~\ref{thm:Klingenberg-hypersurface} we have
\[
    \ninj^{\Gamma_{r_0}}(x)\geq \frac{1}{6}\min\{r_1,\conv(x),\tfrac 65 \compR^c_{\rho(r_1)}\}\eqcolon r_2\,.
\]
Since the absolute value of the left-hand side of \eqref{eqn:parallel-hypersurface-curvatures} is greater than or equal to the right-hand side and using $\tfrac 15\compR^c_{\rho(r_2)}\geq \tfrac 15\compR^c_{\rho(r_1)}\geq r_2$, Lemma~\ref{lem:curv-par-hyp} implies that the hypersurfaces parallel to $\Gamma_{r_2}$ up to distance $r_2$ have their second fundamental form bounded in norm by
\[
    C^F\coloneq \ct_c\left(R^c_{-\rho(r_2)}-\tfrac 15R^c_{-\rho(r_2)}\right) = \ct_c\left(\tfrac 45 R^c_{-\rho(r_2)}\right)\,.
\]
In order to apply Lemma~\ref{lem:gradient-estimate} with $N=\Gamma_{r_2}$ and $I=(-r_2,r_2)$, we have to determine the metric distortion constant $C^G$, which we do by standard Jacobi-field comparison. Namely, with $p\in \Gamma_{r_2}$, $\nu$ a unit normal to $T_p\Gamma_{r_2}$ and $w\in T_p\Gamma_{r_2}$ with $|w|=1$, we have by \cite[Theorem~11.9(b)]{lee},
\[
    |\dd (\exp_p)_{t\nu}(w)| = \frac{J(t)}{t} \leq \frac {\sn_{-c}(t)}t\cdot\left|\frac{\nabla}{\dd t}J(0)\right| = \frac {\sn_{-c}(t)}t\,,
\]
where $J$ is the unique Jacobi-field along $t\mapsto \exp_p(t\nu)$ with $J(0)=\nu$ and $(\nabla/\dd t) J(0)=w$. By monotony of the right-hand side, we thus obtain $C^G\coloneq \sn_{-c}(r_2)/r_2$. Together with $C^\Sigma\coloneq \rho(r_2)$, Lemma~\ref{lem:gradient-estimate} yields
\[
    \alpha = \alpha(C^\Sigma, C^F,C^G) = \frac{1}{C^G}\int\limits_0^\infty\!\frac{1}{C^\Sigma(1+\sigma^2)^{3/2}+C^F(1+\sigma^2)} \,\dd\sigma\,.
\]
Set $r_3\coloneq \min\{r_2,\alpha\}$. Lastly, to ensure that any potential graphing function does not take values greater than the normal radius of $\Gamma_{r_3}$, we integrate the gradient estimate \eqref{eqn:gradient-estimate}. Hence,
\[
    C\coloneq \min\{r_3,r_4\}\,,
\]
where
\[
    r_4\coloneq C^G\int\limits_0^{r_3}\psi^{-1}(s)\,\dd s\,,
\]
is our desired constant.
\end{proof}
\begin{remark}
    Note that the integral \eqref{eqn:psi-integral} can be evaluated in closed form. Hence, $\alpha$ and $r_4$ can be easily calculated numerically.
\end{remark}


\section{Submanifolds of bounded geometry}

In this section, we let $M$ be an $n$-dimensional manifold and $\Sigma\subseteq M$ a properly embedded submanifold, not necessarily of codimension one at first.\\

We recall some concepts from \cite{nardmann}. To this end, denote for $m,d\in\N$ with $\R\mathrm{Poly}_m^d$ the $\R$-vector space of real polynomials of (total) degree $\leq d$ in $m$ variables, equipped with the usual Euclidean topology. For a Riemannian metric $g$ on $M$ and $u\in \Csf^\infty(M)$, we let $g[u]\coloneqq\e^{2u}g$ and for a vector field $X$ we also define $X[u]\coloneqq \e^{-u}X$ so that $g[u](X[u],X[u]) = g(X,X)$.

\begin{definition}\label{def:flatzoomer}
	A map $\Phi:\Csf^\infty(M,\R)\to \Csf^0(M,\R_{\geq 0})$ is a \emph{flatzoomer} if for some - and then every - Riemannian metric $\eta$ on $M$, there exist $k,d\in \N_0$, $\alpha\in\R_{> 0}$, $u_0\in \Csf^0(M,\R)$ and a polynomial-valued map $P\in \Csf^0(M,\R\mathrm{Poly}_{k+1}^d)$ such that
	\begin{equation}\label{eqn:flatzoomer-inequality}
		\Phi(u)(x) \leq \e^{-\alpha u(x)} P(x)\left(u(x), \left|\nabla^1_\eta u \right|_\eta(x) ,\ldots,\left|\nabla^k_\eta u \right|_\eta(x) \right)
	\end{equation}
	for all $u\in \Csf^{\infty}(M,\R)$ with $u>u_0$ and $x\in M$.
\end{definition}

\begin{example}\label{ex:flatzoomer} Let $\mathcal{F}$ be a foliation on $M$, $g$ a Riemannian metric on $\mathcal{F}$ and $k\in \N_0$.
	\begin{enumerate}[(a)]
		\item \cite[Example~2.5]{nardmann} The map $\Phi^{(M,g)}_{\nabla^k\mathrm{Riem}}:C^\infty(M,\R)\to C^0(M,\R_{\geq 0})$,
		\[
			\Phi^{(M,g)}_{\nabla^k\mathrm{Riem}}(u)\coloneqq \left|\nabla^k_{g[u]} \mathrm{Riem}_{g[u]}\right|_{g[u]}
		\]
		is a flatzoomer.
		\item \cite[Example~2.6]{nardmann} With $\secfund$ the second fundamental form of (the leaves of) $\mathcal F$ with respect to $g$, the map $\Phi^{(\Fcal,g)}_{\nabla^k{\secfund}}\Csf^\infty(M,\R)\to \Csf^0(M,\R_{\geq 0})$,
		\[
			\Phi^{(\Fcal,g)}_{\nabla^k{\secfund}}(u)\coloneqq \left|\nabla^k_{g[u]} \secfund_{g[u]}\right|_{g[u]}
		\]
		is a flatzoomer
	\end{enumerate}
\end{example}

As remarked on top of page 157 in \cite{nardmann}, the definition of a flatzoomer is inadequate to cover the (inverse) injectivity radius, hence the more general definition of a quasi-flatzoomer below.

Denote with $\Fct(M,\R_{\geq 0})$ the space of (rough) functions on $M$ with values in the nonnegative real numbers.

\begin{definition}\label{def:quasi-flatzoomer} Let $\Kcal=(\Kcal_i)_{i\in\N_0}$ be a compact exhaustion\footnote{\cite[Conventions~1.1]{nardmann} defines a compact exhaustion without any regularity assumptions on the individual sets.} of $M$. Set additionally $\Kcal_{-2}\coloneqq\Kcal_{-1}:=\emptyset$. A map $\Phi:\Csf^\infty(M,\R)\to\Fct(M,\R_{\geq 0})$ is a \emph{quasi-flatzoomer for} $\Kcal$ if for some - and then every - Riemannian metric $\eta$ on $M$, there exist $k,d\in\N_0$, $\alpha\in\R_{> 0}$, $u_0\in \Csf^0(M,\R)$ and $P\in \Csf^0(M,\R\mathrm{Poly}^d_{k+1})$ such that
	\[
		\Phi(u)(x)\leq \sup \left\{\e^{-\alpha u(y)}P(y)\left(u(y), \left|\nabla^1_\eta u \right|_\eta(y) ,\ldots,\left|\nabla^k_\eta u \right|_\eta(y) \right)\,\Big|\,y\in \mathcal K_{i+1}\setminus \mathcal K_{i-2}    \right\}
	\]
	for all $i\in\N_0$, $x\in \mathcal K_{i}\setminus\mathcal K_{i-1}$ and $u\in\Csf^\infty(M,\R)$ which satisfy $u>u_0$ on $\Kcal_{i+1}\setminus \Kcal_{i-2}$.
\end{definition}

\begin{example}\label{ex:quasi-flatzoomers} Let $\Fcal$ be a foliation on $M$ and $g$ a Riemannian metric on $\Fcal$.
	\begin{enumerate}[(a)]
		\item Every flatzoomer $\Phi$ is a quasi-flatzoomer for every compact exhaustion.
		\item Let $\Kcal=(\Kcal_i)_{i\in\N_0}$ be a compact exhaustion of $M$, $m\in\N$, $\Phi_i:\Csf^\infty(M,\R)\to\Fct(M,\R_{\geq 0})$, $1\leq i \leq m$, quasi-flatzoomers for $\Kcal$ and $Q\in\Csf^0(M\times (\R_{\geq 0})^m, \R_{\geq 0})$ \emph{homogeneous-polynomially bounded}, i.e.\ assume there exists $r\in\R_{> 0}$ and $c\in\Csf^0(M,\R_{\geq 0})$ such that for all $x\in M$ and $v_1,\ldots,v_m\in[0,1]$ one has
		\[
			Q(x,v_1,\ldots,v_m)\leq c(x)\cdot (v_1+\ldots+v_m)^r\,.
		\] 
		Then the map $\Phi:\Csf^\infty(M,\R)\to\Fct(M,\R_{\geq 0})$ given by
		\[
			\Phi(u)(x)\coloneq Q(x,\Phi_1(u)(x),\ldots,\Phi_m(u)(x))
		\]
		is a quasi-flatzoomer for $\Kcal$, \cite[Example~2.11]{nardmann}. In particular, sums, products, powers, minima and maxima of quasi-flatzoomers are quasi-flatzoomers.
		\item Denote by $\mathrm{conv}^\mathcal{F}_g$ and $\mathrm{inj}^\mathcal{F}_g$ the convexity respectively injectivity radius of the leaves of $(\Fcal,g)$. Then $\Phi^{(\Fcal,g)}_{\mathrm{inj}},\Phi^{(\Fcal,g)}_{\mathrm{conv}}:\Csf^\infty(M,\R)\to \Fct(M,\R_{\geq 0})$ defined by
		\begin{align*}
			\Phi^{(\Fcal,g)}_{\mathrm{inj}}(u) &\coloneqq 1/\mathrm{inj}^\mathcal{F}_{g[u]}\,,\\
			\Phi^{(\Fcal,g)}_{\mathrm{conv}}(u) & \coloneqq1/\mathrm{conv}^\mathcal{F}_{g[u]}\\
		\end{align*}
		are quasi-flatzoomers for every compact exhaustion $\mathcal K$ of $M$, \cite[Theorem~3.8]{nardmann}.
	\end{enumerate}
\end{example}

\begin{remark}\label{rem:conv-bounded-polynomial}
    An inspection of the arguments in the proof that $\Phi^{(M, g)}_{\conv}$ resp.\ $\Phi^{(M, g)}_{\inj}$ is a quasi-flatzoomer (i.e.\ with the trivial foliation $\Fcal = M$) shows that if the metric $g$ has positive injectivity radius, uniformly bounded $\left|\mathrm{Riem}_{g}\right|_g$ and $\left|\nabla_g\mathrm{Riem}_{g}\right|_g$, and one chooses the compact exhaustion $\mathcal K$ in such a way that there exists a constant $D>0$ such that for all $i\in\N_0$ one has 
    \begin{gather}\label{eqn:compact-exhaustion-boundary-spacing}
        \begin{aligned}
            \dist(\Kcal_i\setminus\Kcal_{i-1}, \partial\Kcal_{i+1}) > D\,,  \\
            \dist(\Kcal_i\setminus\Kcal_{i-1}, \partial\Kcal_{i-2}) > D
        \end{aligned}
    \end{gather}
     then the function $u_0$ and the polynomial $P$ from Definition~\ref{def:quasi-flatzoomer} can be chosen bounded respectively with bounded coefficients.

    Indeed, letting $\iota = \inj_g^M$ we choose a $\iota/3$-net $\{x_i\,|\,i\in\N\}$ in $(M,\dist^g)$ and, in the notation of the proof of Theorem~3.8 in \cite{nardmann}, let the covering $\mathcal U = (U_i)_{i\in\N_0}$ be given by $U_i\coloneqq\Bsf_{2\iota / 3}(x_i)$. Then each $U_i$ is compactly contained in an exponential chart $\varphi_i$ with domain $\Bsf_\iota(x_i)$ and by the usual Gromov packing argument, using that $\mathrm{Ric}_g$ is bounded from below by a constant, one finds that $\mathcal U$ is locally finite. Since the curvature and its derivative are bounded, there exist constants $A,C>0$ such that for all $i\in\N_0$ one can choose $A_i\coloneqq A$, $C_i\coloneqq C$ on p.\ 161 of \cite{nardmann}. Consequently, the function $H\in \Csf^{0}(M,\R_{>0})$ can be chosen constant. Next, given our assumption \eqref{eqn:compact-exhaustion-boundary-spacing} and the choice of $\mathcal U$, for all $i\in\N_0$ and all $x\in\Kcal_i\setminus\Kcal_{i-1}$ one has $\overline{\Bsf}^g_{\min\{D,\iota/4\}}(x)\subseteq U_j\cap \left(\Kcal_{i+1}\setminus\Kcal_{i-2}\right)$, so that the function $u_1\in \Csf^0(M,\R)$ on top of page 162 of \cite{nardmann} can be chosen constant. Going through the rest of the proof of Theorem~3.8 in \cite{nardmann}, the above claim is now evident.
\end{remark}

Here comes the main flatzoomer theorem.

\begin{theorem}[{\cite[Theorem~4.1]{nardmann}}]\label{thm:main-flatzoomer-thm}
    Let $\Kcal = (\Kcal_i)_{i\in\N_0}$ be a smooth compact exhaustion of $M$, let $(\Phi_i)_{i\in\N_0}$ be a sequence of quasi-flatzoomers for $\Kcal$, let $(\varepsilon_i)_{i\in\N_0}$ be a sequence in $\Csf^0(M,\R_{>0})$ and let $w\in \Csf^0(M,\R)$. Then there exists a real-analytic $u:M\to\R$ with $u>w$ such that
    \[
        \forall i\in \N_0: \Phi_i(u) < \varepsilon_i \quad\text{ holds on } \quad M\setminus \Kcal_i\,.
    \]
\end{theorem}

From the above theorem it is straight-forward to obtain on a given Riemannian manifold $(M,g)$ a conformal metric $g[u]$ such that $(M,g[u])$ has bounded geometry. Indeed, choosing an arbitrary smooth compact exhaustion $\Kcal$, $\Phi_0=\Phi^{(M,g)}_{\inj}$ and $\Phi_i=\Phi^{(M,g)}_{\nabla^{i-1}\mathrm{Riem}}$ for all $i>0$, $\varepsilon_i=C_i\in\R$ $(i\in\N_0)$ and $w=0$ we obtain a $u\in \Csf^{\infty}(M,\R_{>0})$ as desired.\\

We will now turn our attention to submanifolds. Theorem~\ref{thm:submanifold-qfz} below shows that a quasi-flatzoomer of a properly embedded submanifold can be suitably extended to a quasi-flat\-zoomer on the ambient manifold. For its proof, we need to express the $k$-th covariant derivative of function applied to tangent vectors of the submanifold in terms of the $k$-th covariant derivative of the function restricted to the submanifold and certain terms that are normal to the submanifold. To that end, denote with $\nabla_g$ the covariant derivative of $(M,g)$, with $\nabla_\Sigma$ the one of $(\Sigma,g_{|\Sigma})$ and let $\odot$ be the symmetric tensor product. For vectors $X_1,\ldots,X_k\in T_pM$ and $I\subseteq \{1,\ldots,k\}$ let $I^c\coloneqq \{1,\ldots,k\}\setminus I$ and define
\[
	X_I\coloneqq X_{i_1}\odot \ldots\odot X_{i_\ell}
\]
where $I = \{i_1<\ldots<i_\ell\}$.
\begin{lemma}\label{lem:covariant-derivative-submlfd}
	Let $u\in\Csf^\infty(M)$, $k\in\N_{\geq 2}$ and $X_1\ldots,X_k\in T_p\Sigma$. Then
	\begin{equation}\label{eqn:covariant-derivative-submldf-0}
		\nabla^k_gu(X_1,\ldots,X_k) = \nabla^k_\Sigma u(X_1,\ldots,X_k) - \sum_{\ell = 2}^{k}\sum_{\substack{I\subseteq\{1,\ldots,k\}\\ |I|=\ell}}\nabla_g^{k-\ell+1}u(X_{I^c},\nabla_g^{\ell-2}\secfund(X_I))\,.
	\end{equation}
\end{lemma}
\begin{proof}
	We set $\widetilde{\nabla}\coloneqq\nabla_g$ and $\nabla\coloneqq\nabla_\Sigma$.
	Let $k = 2$ and $X_1,X_2\in T_p\Sigma$. Extend $X_2$ to smooth local vector field tangent to $\Sigma$. Then
	\begin{align*}
		\widetilde{\nabla}^2 u(X_1,X_2) &= X_1(\widetilde{\nabla}u(X_2)) - \widetilde{\nabla}u(\widetilde{\nabla}_{X_1}X_2)\\
				&=X_1(\nabla u(X_2)) - \widetilde{\nabla}u(\nabla_{X_1}X_2+\secfund(X_1,X_2))\\
				&=X_1(\nabla u(X_2)) - \nabla u(\nabla_{X_1}X_2)+\widetilde{\nabla} u(\secfund(X_1,X_2))\\
				&=\nabla^2 u(X_1,X_2) - \widetilde{\nabla} u(\secfund(X_1,X_2))\,,
	\end{align*}
	which is \eqref{eqn:covariant-derivative-submldf-0} for $k=2$.
	
	Next, assume \eqref{eqn:covariant-derivative-submldf-0} holds for a fixed $k$ and all $X_1,\ldots,X_k\in T_p\Sigma$. Fix $X_1,\ldots,X_{k+1}\in T_p\Sigma$. Using Fermi-coordinates over normal coordinates about $p$ in $\Sigma$, we extend $X_1,\ldots,X_{k+1}$ to smooth local vector fields tangent to $\Sigma$ which satisfy $\nabla X_i = 0$ in $p$ for all $i=1,\ldots,k+1$. The Gauss equation then reads $\widetilde{\nabla}_{X_i}X_j = \secfund(X_i,X_j)$ so that we have
	\begin{align*} 
		\widetilde{\nabla}^{k+1}u(X_1,\ldots,X_k) &= X_1(\widetilde{\nabla}^ku(X_2,\ldots,X_{k+1}))\\&-\sum_{i = 2}^{k+1}\widetilde{\nabla}^ku(X_2,\ldots,X_{i-1},\hat{X}_i,X_{i+1},\ldots,X_{k+1},\widetilde{\nabla}_{X_1}X_i)\\
		& = X_1\left(\nabla^ku(X_2,\ldots,X_{k+1}) - \sum_{\ell = 2}^{k}\sum_{\substack{I\subseteq\{2,\ldots,k+1\}\\ |I|=\ell}}\widetilde{\nabla}^{k-\ell+1}u(X_{I^c},\nabla_g^{\ell-2}\secfund(X_I))\right)\\&-\sum_{i = 2}^{k+1}\widetilde{\nabla}^ku(X_2,\ldots,X_{i-1},\hat{X}_i,X_{i+1},\ldots,X_{k+1},\secfund(X_1,X_i))\\
		& = \nabla^{k+1}u(X_1,\ldots,X_{k+1})\\
		& - \sum_{\ell = 2}^{k}\sum_{\substack{I\subseteq\{2,\ldots,k+1\}\\ |I|=\ell}}\big(\widetilde{\nabla}^{k+1-\ell+1}u(X_{I^c\cup\{1\}},\nabla_g^{\ell-2}\secfund(X_I)) \\ &\mkern140mu+  \widetilde{\nabla}^{k-\ell+1}u(X_{I^c},\nabla_g^{\ell-1}\secfund(X_{I\cup\{1\}}))\big)\\
		&-\sum_{i = 2}^{k+1}\widetilde{\nabla}^ku(X_1,\ldots,X_{i-1},\hat{X}_i,X_{i+1},\ldots,X_{k+1},\secfund(X_1,X_i))\\
		& = \nabla^{k+1}u(X_1,\ldots,X_{k+1}) - \sum_{\substack{I\subseteq \{1,\ldots,k+1\} \\ |I| = 2}} \widetilde{\nabla}^{k}u(X_{I^c},\secfund(X_I))\\
		& - \sum_{\ell = 3}^{k}\sum_{\substack{I\subseteq\{2,\ldots,k+1\}\\|I|=\ell}}\widetilde{\nabla}^{k+1-\ell+1}u(X_{I^c\cup\{1\}},\nabla_g^{\ell-2}\secfund(X_I))\\
		& - \sum_{\ell = 2}^{k}\sum_{\substack{I\subseteq \{2,\ldots,k+1\} \\ |I|=\ell}} \widetilde{\nabla}^{k-\ell+1}u(X_{I^c},\nabla_g^{\ell-1}\secfund(X_{I\cup\{1\}}))\\
		& = \nabla^{k+1}u(X_1,\ldots,X_{k+1}) - \sum_{\substack{I\subseteq \{1,\ldots,k+1\} \\ |I| = 2}} \widetilde{\nabla}^{k}u(X_{I^c},\secfund(X_I))\\
		& - \sum_{\ell = 3}^{k}\sum_{\substack{I\subseteq\{2,\ldots,k+1\}\\ |I|=\ell}}\widetilde{\nabla}^{k+1-\ell+1}u(X_{I^c\cup\{1\}},\nabla_g^{\ell-2}\secfund(X_I))\\
		& - \sum_{\ell = 3}^{k+1}\sum_{\substack{I\subseteq \{2,\ldots,k+1\} \\ |I|=\ell-1}} \widetilde{\nabla}^{k+1-\ell+1}u(X_{I^c},\nabla_g^{\ell-2}\secfund(X_{I\cup\{1\}}))\\
         & = \nabla^{k+1}_\Sigma u(X_1,\ldots,X_{k+1}) - \sum_{\ell = 2}^{k+1}\sum_{\substack{I\subseteq\{1,\ldots,k+1\}\\ |I|=\ell}}\nabla_g^{k+1-\ell+1}u(X_{I^c},\nabla_g^{\ell-2}\secfund(X_I))\,.
	\end{align*}
\end{proof}

By virtue of our metric $g$ we have an isometric inclusion $(T^*\Sigma)^{\otimes k}\hookrightarrow (T^*M)^{\otimes k}$. In particular, we can view each term on the right hand side of \eqref{eqn:covariant-derivative-submldf-0} as an element of $(T_p^*M)^{\otimes k}$ for $p\in\Sigma$.
\begin{corollary}\label{cor:norm-covariant-derivative-function-submfld}
Let $u\in\Csf^\infty(M)$. Then on $\Sigma$ one has
\[
	|\nabla_\Sigma^k u|_g\leq |\nabla_g^k u|_g + \sum_{\ell = 2}^{k}\binom{k}{l} |\nabla_g^{\ell-2}\secfund|_g\cdot |\nabla_g^{k-\ell+1}u|_g\,.
\]
\end{corollary}
\begin{proof}
	This follows from Lemma~\ref{lem:covariant-derivative-submlfd} by noting that there are $\binom{k}{\ell}$ ways to choose $I\subseteq\{1,\ldots,k\}$ with $|I|=\ell$.
\end{proof}

\begin{theorem}[submanifold quasi-flatzoomers are quasi-flatzoomers]\label{thm:submanifold-qfz} Let $\Sigma$ be a properly embedded submanifold of $M$ and $\Kcal=(\Kcal_i)_{i\in \N_0}$ a compact exhaustion of $M$. Set $\Kcal':=\left(\Kcal_i\cap \Sigma\right)_{i\in\N_0}$ and let $\Phi':\Csf^\infty(\Sigma,\R)\to \Fct(\Sigma,\R_{\geq 0})$ be a quasi-flatzoomer for $\Kcal'$ on $\Sigma$. Then
	\begin{align*}
		\Phi : \Csf^\infty(M,\R)&\to \Fct(M,\R_{\geq 0})\\
		\Phi(u)(x)&\coloneqq\begin{cases}
			\Phi'(u_{|\Sigma})(x) & \qquad \text{ if } x\in \Sigma,\\
			0&\qquad\text{ else },
		\end{cases}
	\end{align*}
	is a quasi-flatzoomer for $\Kcal$ on $M$.
\end{theorem}
\begin{proof}
	By definition, there exist $\tilde{u}_0\in \Csf^{0}(\Sigma,\R),\tilde{b}_0\in \Csf^{0}(\Sigma,\R_{\geq 0})$, $\alpha > 0$ and $k,d\in\N$ such that
	\begin{equation*}
		\Phi'(\tilde{u})(x)\leq \sup\left\{\e^{-\alpha\tilde{u}(y)}\tilde{b}_0(y)\left(1+\sum_{\ell=0}^{k}\left|\nabla^\ell_\Sigma \tilde{u}\right|_{g}(y)\right)^d \,\Big|\,y\in \mathcal K'_{i+1}\setminus\mathcal K'_{i-2}      \right\}
	\end{equation*}
	for all $i\in\N_0$, $x\in \Kcal'_{i}\setminus\Kcal'_{i-1}$ and $\tilde{u}\in\Csf^\infty(\Sigma,\R)$ which satisfy $\tilde{u}>\tilde{u}_0$ on $\Kcal'_{i+1}\setminus \Kcal'_{i-2}$.
	
	By virtue of Corollary~\ref{cor:norm-covariant-derivative-function-submfld} we now find a function $\tilde{b}_1\in \Csf^{0}(\Sigma,\R_{\geq 0})$ such that
	\begin{equation*}
		\Phi'(\tilde{u})(x)\leq \sup\left\{\e^{-\alpha u(y)}\tilde{b}_1(y)\left(1+\sum_{\ell=0}^{k}\left|\nabla^\ell_g u\right|_{g}(y)\right)^d \,\Big|\,y\in \mathcal K'_{i+1}\setminus\mathcal K'_{i-2}      \right\}
	\end{equation*}
	for all $i\in\N_0$, $x\in \Kcal'_{i}\setminus\Kcal'_{i-1}$ and $u\in\Csf^\infty(M,\R)$ which satisfy $u>\tilde{u}_0$ on $\Kcal'_{i+1}\setminus \Kcal'_{i-2}$.
	
	Since $\Sigma$ is a closed set, we can extend $\tilde{u}_0$ and $\tilde{b}_1$ to functions $u_0\in \Csf^0(M)$ resp.\ $b\in \Csf^0(M,\R_{\geq 0})$.	We obtain
	\[
		\Phi(u)(x)\leq \sup\left\{\e^{-\alpha u(y)}b(y)\left(1+\sum_{\ell=0}^{k}\left|\nabla^\ell_{g}u\right|_{g}(y)\right)^d \,\Big|\,y\in \mathcal K_{i+1}\setminus\mathcal K_{i-2}\right\}
	\]
	for all $i\in\N_0$, $x\in \Kcal_{i}\setminus\Kcal_{i-1}$ and $u\in\Csf^\infty(M,\R)$ which satisfy $u>u_0$ on $\Kcal_{i+1}\setminus \Kcal_{i-2}$.
\end{proof}

\begin{remark}\label{rem:submanifold-qfz-bounded-polynomial}
    The proof of Theorem~\ref{thm:submanifold-qfz} above shows that if 
    $\Phi':\Csf^\infty(\Sigma,\R)\to \Fct(\Sigma,\R_{\geq 0})$ is such that
    $u_0\in\Csf^0(\Sigma,\R)$ and $P\in \Csf^0(\Sigma,\R\mathrm{Poly}_{k+1}^d)$ corresponding to the Riemannian metric $g{|\Sigma}$ as in Definition~\ref{def:quasi-flatzoomer} can be chosen bounded resp.\ with bounded coefficients, then the same can be done for the analogous data of $\Phi:\Csf^\infty(M,\R)\to \Fct(M,\R_{\geq 0})$ if $|\nabla^i\secfund^\Sigma_g|_g$ is bounded for all $i=0,\ldots,k-2$. 
\end{remark}

\begin{example}\label{ex:submanifold-quasi-flatzoomer}
	Let $\Sigma$ be a properly embedded submanifold of $M$ and $g$ a Riemannian metric on $M$. Then the following maps $\Phi_*:\Csf^\infty(M,\R)\to \mathsf{Fct}(M,\R_{\geq 0})$ are quasi-flatzoomers on $M$ for every compact exhaustion:
	\begin{enumerate}[(i)]
		\item 	\begin{align*}
					\Phi_\mathrm{conv}^{(\Sigma,g)}(u)(x):=\begin{cases}
						\Phi^{(\Sigma,g|\Sigma)}_{\mathrm{conv}}(u_{|\Sigma})(x) & \text{ if } x\in \Sigma\,,\\
						0&\text{ else}\,,
					\end{cases}
				\end{align*}
		
		\item	\begin{align*}
					\Phi_\mathrm{inj}^{(\Sigma,g)}(u)(x):=\begin{cases}
						\Phi^{(\Sigma,g|\Sigma)}_{\mathrm{inj}}(u_{|\Sigma})(x) & \text{ if } x\in \Sigma\,,\\
						0&\text{ else}\,,
					\end{cases}
				\end{align*}
		\item 	\begin{align*}
					\Phi^{(\Sigma,g)}_{\nabla^k\mathrm{Riem}}(u)(x) := \begin{cases}
						\Phi^{(\Sigma,g|\Sigma)}_{\nabla^k\mathrm{Riem}}(u_{|\Sigma})(x)&\text{ if } x\in \Sigma\,,\\
						0& \text{ else}\,.
					\end{cases}
				\end{align*}		
	\end{enumerate}
\end{example}
The second fundamental form of a submanifold cannot shown to be a quasi-flatzoomer with a combination of Examples~\ref{ex:flatzoomer}, \ref{ex:quasi-flatzoomers} and Theorem~\ref{thm:submanifold-qfz}. Indeed, the second fundamental form of a submanifold is an extrinsic datum and is therefore not expressible purely by data of $(\Sigma,g_{|\Sigma})$. However, we can still leverage Example~\ref{ex:flatzoomer}(b) by using the foliation furnished by a tubular neighbourhood of $\Sigma$ in $M$.
\begin{proposition}\label{prop:second-fund-form-qfz}
	Let $\Sigma\subseteq M$ be a properly  embedded submanifold of codimension $n-\ell < n=\dim M$ and with trivialisable normal bundle $N\Sigma$ and let $g$ be a Riemannian metric on $M$. Then the map $\Phi^{(\Sigma,g)}_{\nabla^k{\secfund}}\Csf^\infty(M,\R)\to \Fct(M,\R_{\geq 0})$,
	\begin{equation*}
		\Phi^{(\Sigma,g)}_{\nabla^k{\secfund}}(u)\coloneqq \begin{cases}
			\left|\nabla^k_{g[u]} \secfund^{\Sigma}_{g[u]}\right|_{g[u]}(x) & \text{ if } x\in\Sigma\,,\\
			0 & \text{ else}\,,
		\end{cases}
	\end{equation*}
	is a quasi-flatzoomer for every compact exhaustion of $M$.
\end{proposition}
\begin{proof}	
    Let $U$ be a $g$-tubular neighbourhood of $\Sigma$ of size $\delta:\Sigma\to(0,\infty)$, where we assume without loss of generality that $\delta$ is smooth. Let $(\nu_1,\ldots,\nu_\ell)$ be a global orthonormal frame of $N\Sigma$ and define, for $\alpha\in\R^\ell$ with $|\alpha|_2<1$, $\nu_\alpha\coloneqq\sum_i \alpha_i\nu_i$. Then $\Fcal = (\Sigma_\alpha)_{|\alpha|_2<1}$ with $\Sigma_\alpha\coloneqq \exp(\{(x,v)\in N\Sigma\,|\,v=\delta(x)\cdot \nu_\alpha(x) \})$ is a foliation of $U$ with $\Sigma_0 = \Sigma$. With
	\[
		\Phi^{(\Fcal,g)}_{\nabla^k{\secfund}}:\Csf^\infty(U,\R)\to \Csf^0(U,\R_{\geq 0})
	\]
	as in Example~\ref{ex:quasi-flatzoomers}(b), let $\tilde{b}\in \Csf^{0}(U,\R_{\geq 0})$, $\alpha > 0$ and $d\in\N$ such that
	\begin{equation*}
		\Phi^{(\Fcal,g)}_{\nabla^k{\secfund}}(\tilde{u})(x)\leq \e^{-\alpha\tilde{u}(x)}\tilde{b}(x)\left(1+\sum_{\ell=0}^{k}\left|\nabla^\ell_{g}\tilde{u}\right|_{g}(x)\right)^d 
	\end{equation*}
	for all $\tilde{u}\in\Csf^\infty(U,\R)$ and $x\in U$, where we note that the proof of \cite[Example~2.6]{nardmann} shows that this inequality is independent of any $u_0$ as in Definition~\ref{def:flatzoomer}. In particular, we have for any $u\in\Csf^\infty(M,\R)$ that
	\begin{equation}\label{eqn:submanifold-second-fund-form-quasi-flatzoomer}
		\Phi^{(\Sigma,g)}_{\nabla^k{\secfund}}(u)(x) = \Phi^{(\Fcal,g)}_{\nabla^k{\secfund}}(u_{|U})(x)\leq \e^{-\alpha u(x)}\tilde{b}(x)\left(1+\sum_{\ell=0}^{k}\left|\nabla^\ell_{g}u\right|_{g}(x)\right)^d	\quad \text{ for all } \quad x\in\Sigma\,.
	\end{equation}
	 By virtue of $\Sigma$ being properly embedded, we can extend $\tilde{b}|\Sigma$ to a function $b\in\Csf^0(M,\R_{\geq 0})$. Defining $P\in \Csf^0(M,\R\mathrm{Poly}^d_{k+1})$ by $P(x)(v_1,\ldots,v_{k+1})\coloneqq b(x)(1+v_1\ldots+v_{k+1})^d$, it follows from \eqref{eqn:submanifold-second-fund-form-quasi-flatzoomer} and the fact that $\Phi^{(\Sigma,g)}_{\nabla^k{\secfund}}(u)(x) = 0$ for $x\in M\setminus\Sigma$ that $\Phi^{(\Sigma,g)}_{\nabla^k{\secfund}}$ is a quasi-flatzoomer for every compact exhaustion of $M$.
\end{proof}

\begin{remark}\label{rem:almost-flatzoomer}
    Although $\Phi^{(\Sigma,g)}_{\nabla^k\secfund}$ is formally a quasi-flatzoomer since it takes values in \newline $\Fct(M,\R_{\geq 0})$, the proof above shows that it does satisfy the stronger pointwise estimate \eqref{eqn:flatzoomer-inequality}. We will call each quasi-flatzoomer with this property an \emph{almost-flatzoomer}.
\end{remark}



\begin{lemma}\label{lem:conf-trafo-scalsecfund}
    Let $\Gamma$ be a hypersurface in $(M,g)$ with unit normal $\nu$ and $u\in\Csf^\infty(M)$. Then for all $X$ tangent to $\Gamma$ we have
    \[
        \scalsecfund^{\Gamma,\nu}_g(X,X) = \e^{u}\scalsecfund^{\Gamma,\nu[u]}_{g[u]}(X[u],X[u])+\dd u(\nu)\cdot g(X,X)\,.
    \]
\end{lemma}
\begin{proof}
    Denoting $\nabla=\nabla_g$ and $\widetilde{\nabla} = \nabla_{g[u]}$ for the moment, we calculate straightforwardly,
    \begin{align*}
        \scalsecfund^{\Gamma,\nu[u]}_{g[u]}(X[u],X[u]) &= g[u]\left(\widetilde{\nabla}_{X[u]}X[u],\nu[u]\right) = g\left(\widetilde{\nabla}_{X}X[u],\nu\right)\\
        &=g\left(-\dd u(X)\e^{-u}X+\e^{-u}\widetilde{\nabla}_XX,\nu \right)\\
        &= \e^{-u}g\left(\nabla_XX + 2\dd u(X)X -g(X,X)\cdot\grad_g u,\nu \right)\\
        &= \e^{-u} \left(g(\nabla_XX,\nu) - \dd u(\nu)g(X,X)\right)\\
        &= \e^{-u} \left(\scalsecfund^{\Gamma,\nu}_g(X,X) - \dd u(\nu)g(X,X)\right)\,,
    \end{align*}
    where we have used in the second line the well-known formula for the Levi-Civita connection of a conformal metric, see \cite[Theorem~1.159 a)]{Be87}.
\end{proof}

\begin{lemma}\label{lem:actc-estimate}
    For $c,\lambda\geq 0$ one has
    \[
        \frac{1}{\compR^c_\lambda}\leq \sqrt{c}+\lambda\,.
    \]
\end{lemma}
\begin{proof}
    From the definition of $\compR^c_\lambda$, the inequality $1/\compR^0_\lambda\leq \lambda$ is immediate, which is why we assume $c>0$ in the sequel, i.e.\ we need to prove the inequality
    \[
        \frac{\sqrt{c}}{\arccot\left(\frac{\lambda}{\sqrt{c}}\right)}\leq \sqrt c+\lambda\,.
    \]
    After replacing $\lambda$ by $\sqrt{c}\lambda$ and dividing the inequality by $\sqrt c$, this is equivalent to $ 1/(1+\lambda)\leq \arccot(\lambda)$. Setting $\mu=1/(1+\lambda)$ and applying $\cot$ we arrive, after elementary equivalences, at
    \[
        \frac{\mu}{1-\mu} \geq \tan(\mu)\,.
    \]
    This is readily seen to be true by
    \[
        \frac{\mu}{1-\mu} = \int_0^\mu\!\frac{1}{(1-y)^2}\,\dd y  \geq \int_0^\mu\!\frac{1}{\cos^2(y)}\,\dd y = \tan(\mu)\,. 
    \]
\end{proof}

The following theorem shows that the inverse of the normal injectivity radius of a hypersurface is a quasi-flatzoomer. Of course, we use Theorem~\ref{thm:Klingenberg-hypersurface} to prove this. Looking at \eqref{eqn:pointwise-lower-estimate-ninj} and recalling that we already know $1/\conv$ to be a quasi-flatzoomer, we have to show that $1/\compR^c_\lambda$ and the inverse of the ``slice radius''  are quasi-flatzoomers. For $1/\compR^c_\lambda$ this is easy, since it is a mixture of the curvature of $(M,g)$ and second fundamental form of $\Sigma$ which we know to be quasi-flatzoomers. For the inverse slice radius, we proceed roughly as follows: since $\Sigma$ is assumed to be properly embedded, Remark~\ref{rem:continuity-seperation-property-and-slice-balls}(iii) asserts that we only have to show that for given $u\in\Csf^\infty(M,\R)$ and $x\in\Sigma$ there exists $s>0$ such that $\Bsf^{g[u]}_t(x)\cap\Sigma$ is connected for all $t\leq s$ and such that $1/s$ satisfies a quasi-flatzoomer estimate. We choose a $g$-tubular neighbourhood $U$ of $\Sigma$, a suitable real number $\lambda>0$ and $s>0$ such that $\overline\Bsf^{g[u]}_t(x)$ is a $g$-$\lambda$-convex domain that is contained in $U$ for all $t\leq s$. Then Corollary~\ref{cor:lambda-convex-body-in-Sigma-tubular-neighbourhood} tells us that $\Bsf^{g[u]}_s(x)$ is a regular $\Sigma$-slice ball.

\begin{theorem}\label{thm:tubular-neighbourhood-qfz}
	Let $(M,g)$ be a Riemannian manifold and $\Sigma\subseteq M$ a properly embedded two-sided hypersurface. Then, the map
	\begin{align*}
		\Phi^{\Sigma,g}_{\ninj}:\Csf^\infty(M,\R) &\to \Fct(M,\R_{\geq 0})\\
		\Phi^{\Sigma,g}_{\ninj}(u)(x) &:=	\begin{cases}
							1/\ninj^{\Sigma}_{g[u]}(x) & \text{ if } x\in \Sigma\,,\\
							0 & \text{ else}\,,
						\end{cases}
	\end{align*}
	is a quasi-flatzoomer for any compact exhaustion of $M$.
\end{theorem}
\begin{proof}
    Let $\Kcal = (\Kcal_i)_{i\in\N_0}$ be a compact exhaustion of $M$ and set $\Kcal_{-2}\coloneqq\Kcal_{-1}\coloneqq \emptyset$. Let $U$ be a $g$-tubular neighbourhood of $\Sigma$. For each $i\in \N_0$ define
    \begin{gather}\label{eqn:g-quantities}
    \begin{aligned}
        c_i^g &\coloneqq \inf\left\{ \seccurv_g(\sigma)\,|\, y\in \Kcal_{i+1}\setminus\Kcal_{i-2}, \sigma\in \operatorname{Gr}_2(T_yM)\right\}\,,\\
        r^g_i&\coloneqq \sup\left\{r > 0\,|\, \forall x\in\Sigma\cap(\Kcal_{i}\setminus\Kcal_{i-1}): r\leq \conv_g(x) \text{ and } \Bsf^g_r(x)\subseteq \Kcal_{i+1}\setminus\Kcal_{i-2}\right\}\,,\\
        \delta^g_i&\coloneqq \sup\left\{\delta > 0\,|\,\forall x\in\Sigma\cap(\Kcal_i\setminus\Kcal_{i-1}): \exp_g(U^g_{\delta}(x))\subseteq U\right\}\,,\\
        \lambda^g_i&\coloneqq \inf\big\{\lambda > 0\,|\,\compR^{c^g_i}_\lambda\leq r^g_i/4\big\}\,.
    \end{aligned}
    \end{gather}
    Note that these quantities are all taken with respect to the metric $g$, which we have emphasised in the notation.

    Since $(\mathring{\Kcal}_{i+1}\setminus\Kcal_{i-2})_{i\in\N_0}$ is a locally finite open covering of $M$, we can choose a function $\lambda\in\Csf^\infty(M,\R_{>0})$ such that $\lambda\geq \lambda_i$ on $\Kcal_i\setminus\Kcal_{i-1}$. In an analogous fashion, we choose a function $u_1\in\Csf^\infty(M,\R_{>0})$ such that for all $i\in\N_0$ and $x\in\Sigma\cap(\Kcal_i\setminus\Kcal_{i-1})$ one has
    \begin{equation}\label{eqn:compact-ball}
        \Bsf^{g[u_1]}_1(x)\subseteq \Bsf^g_{(1/2)\cdot\min\{r^g_i,\delta_i/2\}}(x)\,.
    \end{equation}
    Define flatzoomers $\Phi_{i}:\Csf^\infty(M,\R)\to \Csf^0(M,\R_{\geq 0})$, $i\in\{0,1,2\}$, by $\Phi_0(u)\coloneqq \e^{u_1-u}$, $\Phi_1(u)\coloneqq \Phi^{(M,g)}_{\mathrm{Riem}}(u)^{1/2}+\e^{-u}(\lambda+|\grad_gu|_g)$ and
    \[
        \Phi_2\coloneqq\max\{\Phi_0,2\Phi_1\}\,.
    \]
    Next, define an almost-flatzoomer, see Remark~\ref{rem:almost-flatzoomer}, by
    \[
        \Phi_3(u)\coloneqq2\left( \Phi^{(M,g)}_{\mathrm{Riem}}(u)^{1/2} + \Phi^{(\Sigma,g)}_{\secfund}(u)\right)\,.
    \]
    Finally, define a quasi-flatzoomer for $\Kcal$ by 
    \[
        \Phi:=6\max\left\{\Phi_2,\Phi^{(M,g)}_{\conv},\tfrac{5}{6}\Phi_3\right\}
    \]
    and let $\alpha > 0$, $d\in\N_0$, $u_0\in \Csf^0(M,\R)$ and $P\in \Csf^0(M,\R\mathrm{Poly}^d_{k+1})$ be the corresponding data so that we have
	\[
		\Phi(u)(x)\leq \sup \left\{\e^{-\alpha u(y)}P(y)\left(u(y), \left|\nabla^1_g u \right|_g(y) ,\left|\nabla^2_g u \right|_g(y) \right)\,\Big|\,y\in \mathcal K_{i+1}\setminus \mathcal K_{i-2}    \right\}
	\]
	for all $i\in\N_0$, $x\in \mathcal K_{i}\setminus\mathcal K_{i-1}$ and $u\in\Csf^\infty(M,\R)$ which satisfy $u>u_0$ on $\Kcal_{i+1}\setminus \Kcal_{i-2}$.
    
    We assume without loss of generality that $u_0$ is greater than or equal to $u_1$ and each of the functions ``$u_0$'' in the definition of $\Phi_2$, $\Phi_3$ and $\Phi^{(M,g)}_{\conv}$. Let $i\in\N_0$, $x\in\Kcal_i\setminus\Kcal_{i-1}$ and $u\in\Csf^\infty(M,\R)$ with $u>u_0$ on $\Kcal_{i+1}\setminus\Kcal_{i-2}$. Define
    \begin{align*}
        c&\coloneqq \sup\left\{ |\seccurv_{g[u]}(\sigma)|\,|\, y\in \Kcal_{i+1}\setminus\Kcal_{i-2}, \sigma\in \operatorname{Gr}_2(T_yM)\right\}\,,\\
        q&\coloneqq \inf_{\Kcal_{i+1}\setminus\Kcal_{i-2}}\e^{u-u_1}\,,\\
        s&\coloneqq \min\left\{\mathrm{act}_c\left(\sup_{\Kcal_{i+1}\setminus\Kcal_{i-2}}\e^{-u}(\lambda+|\grad_gu|_g)\right), q\right\}\,,\\
        C&\coloneqq \sup \left\{ |\seccurv_{g[u]}(\sigma)|\,|\, y\in \overline{\Bsf}^{g[u]}_s(x)\right\} \,,\\
        \Lambda & \coloneqq \sup \left\{|\secfund^\Sigma_{g[u]}|_{g[u]}(y)\,|\,y\in \overline{\Bsf}_s^{g[u]}(x)  \right\}\,,
    \end{align*}
    where $\mathrm{act}_c$ denotes the inverse of $\mathrm{ct}_c$, see \eqref{eqn:ctc}.
    
    First of all, by Lemma~\ref{lem:actc-estimate} we have
    \begin{multline}\label{eqn:inverse-slice-radius-qfz}
        \frac1s = \max\left\{\frac{1}{\mathrm{act}_c\left(\sup_{\Kcal_{i+1}\setminus\Kcal_{i-2}}\e^{-u}(\lambda+|\grad_gu|_g)\right)},\frac1q  \right\}\\ \leq \max\left\{\sqrt{c} + \sup_{\Kcal_{i+1}\setminus\Kcal_{i-2}}\e^{-u}(\lambda+|\grad_gu|_g), \sup_{\Kcal_{i+1}\setminus\Kcal_{i-2}}\e^{u_1-u}     \right\}\\ \leq \sup \left\{\Phi_2(u)(y)\,|\, y\in\Kcal_{i+1}\setminus\Kcal_{i-2}\right\}\,.
    \end{multline}
    Moreover,    
    \begin{align*}
        \Bsf_s^{g[u]}(x)\subseteq \Bsf_q^{g[u]}(x) = \Bsf_q^{\e^{2u-2u_1}g[u_1]}(x)\subseteq \Bsf_{q}^{q^2g[u_1]}(x)= \Bsf_{1}^{g[u_1]}(x)\,,
    \end{align*}
    and since the ball on the right-hand side of \eqref{eqn:compact-ball} is compactly contained in a $g$-strongly convex ball, $\Bsf_s^{g[u]}(x)$ is precompact.

    Next, note that by definition of $s$, each $t>0$ with $t\leq s$ satisfies
    \[
        t\leq \mathrm{act}_c\left(\sup_{\Kcal_{i+1}\setminus\Kcal_{i-2}}\e^{-u}(\lambda+|\grad_gu|_g)\right) \leq \frac{\pi}{2\sqrt{c}}\,.
    \]
    Hence, we can apply the principal curvature comparison theorem for level sets of distance functions, \cite[Corollary~11.8(a)]{lee}, to find
    \[
        \sup_{\Kcal_{i+1}\setminus\Kcal_{i-2}}\e^{-u}(\lambda+|\grad_gu|_g)\leq\mathrm{ct}_c(s)\leq \mathrm{ct}_c(t)\leq\scalsecfund^{\Gamma,\nu[u]}_{g[u]}\,,
    \]
    where $\Gamma = \Ssf_t^{g[u]}(x)$ and $\nu$ is an inward-pointing $g$-unit normal vector field to $\Gamma$. By Lemma~\ref{lem:conf-trafo-scalsecfund}, we now find that
    \[
        \lambda^g_i\leq\lambda\leq  \lambda + |\grad_gu|_g  +\dd u(\nu)\leq \scalsecfund^{\Gamma,\nu}_g \,,
    \]
    i.e.\ $\overline{\Bsf}_t^{g[u]}(x)$ is a $g$-$\lambda_i$-convex domain for all $t\leq s$. Definitions \eqref{eqn:g-quantities} and \eqref{eqn:compact-ball} ensure that the assumptions of Corollary~\ref{cor:lambda-convex-body-in-Sigma-tubular-neighbourhood} are satisfied, so that $\Bsf_t^{g[u]}(x)$ is a $\Sigma$-slice ball for all $t\leq s$, i.e.\ $\Bsf_s^{g[u]}(x)$ is a regular $\Sigma$-slice ball.

    Next, by Lemma~\ref{lem:actc-estimate} we have
    \begin{equation}\label{eqn:compR-qfz}
        \frac{1}{\compR^C_\Lambda} \leq \sqrt{C} + \Lambda\leq \sqrt{c}+\sup_{\Sigma\cap(\Kcal_{i+1}\setminus\Kcal_{i-2})}\left|\secfund^{\Sigma}_{g[u]}\right|_{g[u]} \leq \sup\{\Phi_3(u)(y)\,|\,y\in\Kcal_{i+1}\setminus\Kcal_{i-2}\}\,.
    \end{equation}

    With Theorem~\ref{thm:Klingenberg-hypersurface} we conclude from \eqref{eqn:inverse-slice-radius-qfz} and \eqref{eqn:compR-qfz}, taking into account that $\Phi_3$ is an almost-flatzoomer, that
    \begin{multline*}
        \Phi^{\Sigma,g}_{\ninj}(x) \leq 6 \max\left\{\frac{1}{s},\frac{1}{\conv^M_{g[u]}(x)},\frac 56 \cdot\frac 1{\compR^C_\Lambda}\right\}\\\leq 6 \max\left\{ \sup_{\Kcal_{i+1}\setminus\Kcal_{i-2}} \Phi_2(u),\frac{1}{\conv^M_{g[u]}(x)}, \frac56 \sup_{\Kcal_{i+1}\setminus\Kcal_{i-2}}\Phi_3(u) \right\} \\
        \leq \sup\left\{\e^{-\alpha u(y)}P(y)\left(u(y), \left|\nabla^1_g u \right|_g(y) ,\left|\nabla^2_g u \right|_g(y) \right)\,\Big|\,y\in \mathcal K_{i+1}\setminus \mathcal K_{i-2} \right\}\,.
    \end{multline*}
\end{proof}

\begin{remark}\label{rem:ninj-bounded-polynomial}
    Similar to Remark~\ref{rem:conv-bounded-polynomial}, the proof of Theorem~\ref{thm:tubular-neighbourhood-qfz} shows that if $g$ has positive convexity radius, bounded curvature, is such that $|\secfund^\Sigma_g|_g$ is bounded and $\Sigma$ possesses a uniform $g$-tubular neighbourhood, and the compact exhaustion $\Kcal$ satisfies \eqref{eqn:compact-exhaustion-boundary-spacing}, then $u_0$ and $P$ from Definition~\ref{def:quasi-flatzoomer} can be chosen bounded respectively with bounded coefficients.
\end{remark}

From the main flatzoomer Theorem~\ref{thm:main-flatzoomer-thm} we now obtain

\begin{theorem}\label{thm:bounded-submanifold-geometry}
	Let $M$ be a manifold, $\Sigma\subseteq M$ a properly embedded two-sided hypersurface and $g_0$ a Riemannian metric on $M$. Let $\Kcal = (\Kcal_i)_{i\in\N_0}$ be a smooth compact exhaustion of $M$, let $\iota, w\in\Csf^0(M,\R_{> 0})$ and $(\varepsilon_i)_{i\in\N_0}$ a sequence in $\Csf^0(M,\R_{>0})$. Then there exists a real-analytic $u:M\to\R$ with $u>w$ such that the metric $g\coloneqq g_0[u]=\e^{2u}g_0$ satisfies:
	\begin{enumerate}[(i)]
		\item $(M,g)$ is complete with $\mathrm{inj}^M_{g}\geq 2\cdot\mathrm{conv}^M_{g}> \iota$,
		\item $(\Sigma, g_{|\Sigma})$ is complete with $\mathrm{inj}^\Sigma_{g_{|\Sigma}}\geq 2\cdot\mathrm{conv}^\Sigma_{g_{|\Sigma}}> \iota_{|\Sigma}$,
		\item $\Sigma$ posses a $g$-uniform tubular neighbourhood in $M$ and $\ninj^{\Sigma}_{g}>2\cdot\iota_{|\Sigma}$,
		\item for every $i\in\N_0$ one has $\left| \nabla^i_g\mathrm{Riem}_g  \right|_g < \varepsilon_i$ on $M\setminus \Kcal_i$ and 
		\item for every $i\in\N_0$ one has $\left| \nabla^i\secfund^\Sigma_g\right|_g<\varepsilon_i$ on $\Sigma \setminus \Kcal_i$.
	\end{enumerate}
\end{theorem}
\begin{proof}
	We define a smooth compact exhaustion $\Kcal'=(\Kcal'_i)_{i\in\N_0}$ by $\Kcal'_0\coloneqq \emptyset$ and $\Kcal'_{i+1}\coloneqq \Kcal'_{i}$ for all $i\in\N_0$. Furthermore, we define $\varepsilon'_0\coloneqq \tfrac{1}{\iota+1}$ and $\varepsilon'_{i+1}\coloneqq \varepsilon_i$ for all $i\in\N_0$ and a sequence $(\Phi_i)_{i\in\N_0}$ of quasi-flatzoomers by
	\[
		\Phi_0\coloneqq \Phi^{(M,g)}_{\mathrm{conv}} + \Phi^{(\Sigma,g)}_{\mathrm{conv}} + \Phi^{\Sigma,g}_{\ninj}
	\]
	and
	\begin{equation*}
		\phantom{\quad \text{ for all }\quad i\in\N\,.}\Phi_i \coloneqq \Phi^{(M,g)}_{\nabla^{i-1}\mathrm{Riem}} + \Phi^{(\Sigma,g)}_{\nabla^{i-1}\secfund}\quad \text{ for all }\quad i\in\N\,.
	\end{equation*}
	Applying Theorem~\ref{thm:main-flatzoomer-thm} to $\Kcal'$, $(\Phi_i)_{i\in\N_0}$, $(\varepsilon'_i)_{i\in\N_0}$ and $w$ yields a real-analytic $u:M\to\R$ for which $g = g_0[u]$ satisfies (i)-(v) since $\conv^M_g>\iota +1>1$, $\conv^{\Sigma}_{g_{|\Sigma}}>\iota_{|\Sigma}+1>1$ and $\ninj^{\Sigma}_g>\iota_{|\Sigma}+1>1$, where we note that on a complete Riemannian manifold one always has $2\conv\leq \inj$.
\end{proof}

For the next corollary, we recall the definition of a Riemannian manifold with boundary of bounded geometry, see, e.g., \cite[Definition~2.2]{schick}.

\begin{definition}\label{boundedgeom} Let $(M,g)$ be a Riemannian manifold with boundary. Then $(M,g)$ is said to have \emph{bounded geometry} if the following holds:
\begin{itemize}
    \item With $\nu$ the unit inward-pointing normal to $\partial M$, there exists $r_0>0$ such that
    \[
        \partial M\times [0,r_0) \ni (x,r)\mapsto \exp_x(r\cdot\nu_x)\in M
    \]
    is a diffeomorphism onto its image $U_{\partial M}$, i.e.\ $U_{\partial M}$ is a uniform collar of $\partial M$ in $M$.
    \item The injectivity radius $\inj^{\partial M}$ of $\partial M$ is positive.
    \item There exists $r_1>0$ such that for all $x\in M\setminus U_{\partial M}$ one has $\inj^M(x)\geq r_1$.
    \item For every $i\in\N_0$ there exists $C_i>0$ such that
    \begin{equation*}
        \left|\nabla^i\mathrm{Riem}^M\right|\leq C_i\qquad \text{ and } \qquad\left|\nabla^i\secfund^{\partial M}\right|\leq C_i\,.
    \end{equation*}
\end{itemize}
\end{definition}

Note that no Riemannian manifold with boundary satisfies $\inj > 0$. Hence, the requirement that $\partial M$ possesses a uniform collar and $\inj$ is uniformly positive away from that collar.

\begin{corollary}\label{cor:existence-metric-bounded-geometry-boundary}
	Let $M$ be a manifold with boundary and $g_0$ a Riemannian metric on $M$. Then there exists $u\in\Csf^\infty(M)$ such that $(M,g)$ with $g\coloneqq g_0[u]$ is of bounded geometry.
\end{corollary}
\begin{proof}
	We realize $(M,g_0)$ as a domain in a Riemannian manifold $(N,h_0)$ without boundary, see \cite[Theorem~A]{stefano}. Applying Theorem~\ref{thm:bounded-submanifold-geometry} to $N$, $\Sigma=\partial M\subseteq N$, $h_0$, an arbitrary smooth compact exhaustion $\Kcal$ of $N$, $\iota = 1$, $w = 0$ and $(\varepsilon_i)_{i\in\N_0}$ with $\varepsilon_i\coloneqq C_i$ for arbitrary constants $C_i>0$ yields a metric $h\coloneqq h_0[u]$ on $N$ whose restriction $g\coloneqq h_{|M}$ to $M\hookrightarrow N$ has bounded geometry.
\end{proof}

\section{Convexification of the boundary}

The goal of this section is to strengthen Corollary~\ref{cor:existence-metric-bounded-geometry-boundary} such that one obtains a metric of bounded geometry and with convex boundary. To that end, we take a step back and let $\Sigma$ be a properly embedded hypersurface in a boundaryless Riemannian manifold $(N,h)$ which is two-sided and divides $N$ into two components. We will assume that $(N,h)$ has bounded geometry, $(\Sigma,g_{|\Sigma})$ has bounded geometry and possesses a uniform tubular neighbourhood $U$ of size $\omega\in\R_{>0}$.

We choose a unit normal vector field $\nu$ along $\Sigma$ and let $\signdist^\Sigma_h$ be the associated signed distance which takes positives values in the component into which $\nu$ points. Note that $\signdist^\Sigma_h$ is smooth on $U$.

Let $f:\R\to\R$ be the function
\begin{equation*}
    f(t) = \begin{cases}
            0& t\leq0\,,\\
            \e^{-1/t}& t>0
        \end{cases}
\end{equation*}
and define for $a,b > 0$
\begin{align*}
        \varphi_{a,b} : \R&\to\R\\
        t&\mapsto 2\cdot \frac{f(a\cdot( t+ b) )}{f(a\cdot( t+ b) )+f(a\cdot(b- t))}-1\,.
\end{align*}
Note that $\varphi_{a,b}$ is smooth, takes values in $[-1,1]$, satisfies $|\varphi_{a,b}(t)|\equiv 1$ for $|t|\geq b$ and $\varphi_{a,b}(0)=0$, $\varphi_{a,b}'(0) = \tfrac{1}{ab^2}$. In particular, $\varphi_{a,b}$ is $\Csf^\infty$-bounded.

With $\omega$ as above, we define for $\lambda > 0$ and $C\in\R$
\begin{gather}\label{eqn:convexification-conformal-factor}
    \begin{aligned}
        \upsilon = \upsilon_{\lambda,\omega,C} : M&\to \R\\
        x&\mapsto C - \varphi_{1/(\lambda\omega^2),\omega}(\signdist^\Sigma_h(x))\,.
    \end{aligned}
\end{gather}
From the properties of $\varphi_{a,b}$ we deduce that $\upsilon$ is smooth, $\Csf^\infty$-bounded and satisfies $\upsilon|\Sigma \equiv C$, $\dd\upsilon(\nu) = -\lambda$.
\begin{lemma}\label{lem:convexification}
If $-\lambda\leq \scalsecfund^{\Sigma,\nu}_h\leq \lambda$, then $0\leq \scalsecfund^{\Sigma,\nu[\upsilon]}_{h[\upsilon]} \leq 2\e^{-C}\lambda$.
\end{lemma}
\begin{proof}
    By Lemma~\ref{lem:conf-trafo-scalsecfund} we have
    \begin{align*}
        \scalsecfund^{\Sigma,\nu[\upsilon]}_{h[\upsilon]}(X[u],X[u]) &= \e^{-\upsilon}\cdot \left(\scalsecfund^{\Sigma,\nu}_h(X,X) - \dd\upsilon(\nu)\cdot h(X,X) \right) \\
        &= \e^{-C}\cdot \left(\scalsecfund^{\Sigma,\nu}_h(X,X) + \lambda\cdot h(X,X) \right)\,,
    \end{align*}
    from which the claim follows.
\end{proof}

\begin{theorem}\label{thm:bounded-geometry-convex-boundary}
    Let $M$ be a manifold with boundary and $g_0$ a Riemannian metric on $M$. Then there exists $u\in\Csf^\infty(M)$ such that $(M,g)$ with $g\coloneqq g_0[u]$ has bounded geometry with convex boundary, i.e.\ $\scalsecfund^{\partial M}_g\geq 0$ with respect to the inward pointing unit normal.
\end{theorem}
\begin{proof}
    As in the proof of Corollary~\ref{cor:existence-metric-bounded-geometry-boundary} we realise $(M,g_0)$ as a domain in a Riemannian manifold $(N,h_0)$ and apply Theorem~\ref{thm:bounded-submanifold-geometry} to obtain $\tilde{u}\in C^\infty(N)$ such that $(N,h)$ with $h\coloneqq h_0[\tilde{u}]$ has bounded geometry with $\conv^N_h>0$, $\Sigma\coloneqq\partial M\subseteq N$ has bounded geometry (with respect to $h$) with $\conv^\Sigma_{h|\Sigma}>0$ and such that $\partial M$ posseses a uniform $h$-tubular neighbourhood $U$ of size $\omega > 0$.

    With $\lambda\in\R_{>0}$ such that $|\scalsecfund^{\Sigma,\nu}_h|_h\leq \lambda$, where $\nu$ is an $h$-unit normal  pointing into $M$, and $C\in\R$ a constant to be determined later, we let $\upsilon=\upsilon_{\lambda,\omega,C}$ as in \eqref{eqn:convexification-conformal-factor}. By Lemma~\ref{lem:convexification}, we have $\scalsecfund^{\Sigma,\nu[\upsilon]}_{h[\upsilon]}\geq 0$. Moreover, since $\upsilon$ is $h$-$\Csf^\infty$-bounded, there exist constants $C_i>0$ for all $i\in\N_0$ such that    \begin{equation*}
            \left|\nabla^i_{h[\upsilon]}\mathrm{Riem}^N_{h[\upsilon]}\right|_{h[\upsilon]}\leq C_i \qquad \text{and} \qquad \left|\nabla^i_{h[\upsilon]}\secfund^\Sigma_{h[\upsilon]}\right|_{h[\upsilon]}\leq C_i 
    \end{equation*}
    for all $i\in\N_0$. If we can show that $\conv^M_{h[\upsilon]}>0$, $\conv^\Sigma_{h[\upsilon]|\Sigma}>0$ and that $\Sigma$ possesses a uniform $h[\upsilon]$-tubular neighbourhood, then $u\coloneqq (\tilde{u}+\upsilon){|M}$ is the desired conformal factor.

    To that end, fix $p\in N$ and let $\varepsilon > 0$. By \cite[Section~2]{GW78}, see also \cite[Theorem~1]{AFLR07}, we can choose an $\varepsilon$-regularisation $r:M\to\R$ of $\dist_p = (\dist^N_h)_p$, i.e.\ $r\in\Csf^\infty(M)$, for all $q\in M$ one has $|r(q)-\dist_p(q)|\leq \varepsilon$ and $\mathrm{Lip}_h(r)\leq \mathrm{Lip}_h(\dist_p)+\varepsilon=1+\varepsilon$. With $D>0$ an arbitrary constant, Sard's theorem garantuees the existence of a positive, strictly increasing sequence $(d_i)_{i\in\N_0}$ of real numbers with $d_{i+1}-d_i\geq D+3\varepsilon$ for all $i\in\N_0$ such that each $d_i$ is a regular value of $r$. Defining $\Kcal_i\coloneqq r^{-1}((-\infty,d_i])$, we thus obtain a smooth compact exhaustion $\Kcal=(\Kcal_i)_{i\in\N_0}$ of $N$.
    With $q\in\Kcal_i\setminus \Kcal_{i-1}$ and $o\in\partial\Kcal_{i+1}$ we have
    \[
        \dist(o,q)\geq \dist_p(o) - \dist_p(q) \geq r(0)-\varepsilon -r(q)-\varepsilon \geq d_{i+1}-d_i-2\varepsilon\geq D+\varepsilon
    \]
    and analogously $\dist(q,o)\geq D+\varepsilon$ for $q\in\Kcal_i\setminus \Kcal_{i-1}$ and $o\in\partial\Kcal_{i-2}$. Hence, $\Kcal$ satisfies \eqref{eqn:compact-exhaustion-boundary-spacing} with $(M,g)$ replaced by $(N,h)$.
    Since $\dist^\Sigma_{h|\Sigma}\geq {\dist^N_h}_{|\Sigma\times\Sigma}$, the compact exhaustion $\Kcal'=(\Kcal_i\cap\Sigma)_{i\in\N_0}$ of $\Sigma$ satisfies \eqref{eqn:compact-exhaustion-boundary-spacing} with $(M,g)$ replaced by $(\Sigma,h_{|\Sigma})$ and $\Kcal$ replaced by $\Kcal'$.

    Consider the quasi-flatzoomers $\Phi_0=\Phi^{(N,h)}_{\conv}$, $\Phi_1=\Phi^{(\Sigma,h)}_{\conv}$ and $\Phi_2=\Phi^{\Sigma,h}_{\ninj}$. By Remark~\ref{rem:conv-bounded-polynomial}, Remark~\ref{rem:submanifold-qfz-bounded-polynomial} and Remark~\ref{rem:ninj-bounded-polynomial}, we can choose $\alpha >0$, for each $j\in\{0,1,2\}$ 
 a bounded $u_j\in\Csf^0(M,\R)$ and $P_j\in\Csf^0(M,\R\mathrm{Poly}^d_{3})$ with bounded coefficients such that
	\begin{equation}\label{eqn:bounded-qfzs}
		\Phi_j(u)(x)\leq \sup \left\{\e^{-\alpha u(y)}P_j(y)\left(u(y), \left|\nabla^1_h u \right|_h(y) ,\left|\nabla^2_h u \right|_h(y) \right)\,\Big|\,y\in \mathcal K_{i+1}\setminus \mathcal K_{i-2}    \right\}
	\end{equation}
	for all $i\in\N_0$, $x\in \mathcal K_{i}\setminus\mathcal K_{i-1}$ and $u\in\Csf^\infty(M,\R)$ which satisfy $u>u_j$ on $\Kcal_{i+1}\setminus \Kcal_{i-2}$.

    Now choose the constant $C$ from above such that $\upsilon>\max\{u_0,u_1,u_2\}$. Since the $P_j$ have bounded coefficients and $\upsilon$ is $h$-$\Csf^\infty$-bounded, it follows from \eqref{eqn:bounded-qfzs} that each $\Phi_j(\upsilon)$ is bounded as desired.
\end{proof}

\section{Orientations of smooth metric measure spaces with boundary}

Let $M^n$ be a smooth connected manifold with (possibly empty) boundary $\partial M $ and set $\mathring{M}  = M\setminus \partial M$. Given a Riemannian metric $g$ on $M$ we denote with Riemannian volume measure with $\mu_g$, turning $(M,g)$ together with the geodesic distance $\dist_g$ into a metric measure space (mms). As above, $\nabla_g$ stands for the Levi-Civita connection and $\mathrm{Ric}_g$ the Ricci curvature. Likewise, $\mathrm{Hess}_g$ is the Hessian and $\Delta_g$ the \emph{nonnegative} Laplace-Beltrami operator.\vspace{1mm}

Let
$$
\mathsf{H}^{1,2}(\mathring{M},g|_{\mathring{M}}):= \{f\in \mathsf{L}^2(\mathring{M},g|_{\mathring{M}})\,|\, \dd f\in \Omega^1_{\mathsf{L}^2}(\mathring{M},g|_{\mathring{M}})\}
$$
denote the Sobolev space, defined using the distribution theory on the Riemannian manifold without boundary $\mathring{M}$. This becomes a Hilbert space with its natural norm
$$
\left\|f\right\|_{\mathsf{H}^{1,2}}:=\left\|f\right\|_{\mathsf{L}^{2}}+\left\|\dd f\right\|_{\mathsf{L}^2}.
$$

Note here that of course 
$$
\mathsf{L}^2(\mathring{M},g|_{\mathring{M}})=\mathsf{L}^2(M,g),\quad \Omega^1_{\mathsf{L}^2}(\mathring{M},g|_{\mathring{M}})=\Omega^1_{\mathsf{L}^2}(M,g),
$$
as the boundary of $M$ has measure zero.\vspace{1mm} 

Let $H_g$ denote the \emph{Neumann realization} of the Laplace-Beltrami operator in $\mathsf{L}^2(M,g)$, so that $\dom(H_g)$ is given by all $f\in \mathsf{H}^{1,2}(\mathring{M},g|_{\mathring{M}})$ such that $\Delta_g f\in L^2(\mathring{M},g|_{\mathring{M}})$. Concretely, $H_g$ is the uniquely determined nonnegative self-adjoint operator in $L^2(M,g)$ with $\dom(H_g)\subset \mathsf{H}^{1,2}(\mathring{M},g|_{\mathring{M}})$ and
\[
\left\langle H_gf_1,f_2\right\rangle_{\mathsf{L}^2} = \left\langle \dd f_1,\dd f_2 \right\rangle_{L^2} \quad\text{ for all $f_1\in \dom(H_g)$, $f_2\in H^{1,2}(\mathring{M},g|_{\mathring{M}})$.}
\]
In particular, for all $f\in \dom(H_g)$ one has $\Delta_g f\in L^2(\mathring{M},g|_{\mathring{M}})$ with $H_gf=-\Delta_g f$.\\
The space $\mathsf{H}^{1,2}(\mathring{M},g|_{\mathring{M}})$ is precisely the form domain $\dom(\sqrt{H_g})$ of $H_g$, which makes $H_g$ the nonnegative self-adjoint operator induced by the strongly local regular Dirichlet form in $\mathsf{L}^2(M,g)$ given by
$$
\mathsf{H}^{1,2}(\mathring{M},g|_{\mathring{M}})\times  \mathsf{H}^{1,2}(\mathring{M},g|_{\mathring{M}})\ni (f_1,f_2)\longmapsto \left\langle \dd f_1,\dd f_2 \right\rangle_{L^2} \in\R.
$$
It has been shown recently in \cite{bgs}, that if $(M,g)$ is metrically complete, then $\Delta_g$ is essentially self-adjoint on the domain given by all $\psi\in \mathsf{C}^{\infty}_c(M)$ with $N_g(\psi)=0$, where $N_g$ is the unit normal vector field. It follows that in this case $H_g$ is the unique self-adjoint extension of the latter symmetric operator. \\

Let $P^g_t:=\exp(-tH_g)$, $t>0$, denote the (Neumann) heat semigroup in $L^2(M,g)$. By the spectral theorem and elliptic regularity one has
$$
P^g_t: \mathsf{L}^2(M,g)\longrightarrow \bigcap_{r>0}\dom(H_g^r)\subset \mathsf{C}^{\infty}(M)\quad\text{for all $t>0$.}
$$ 
Moreover, $H_g P^g_t f=P^g_tH_g  f$ if $f\in\dom(H_g)$, and 
$$
P^{g}_t:\mathsf{L}^q(M,g)\longrightarrow \mathsf{L}^q(M,g)
$$
is a contraction for all $q\in [1,\infty]$, which follows from the Markovian property
\begin{align}\label{saas}
\int P^{g}_t(x,y) \dd\mu_g(y)\leq 1
\end{align}
of the Neumann heat kernel.\vspace{2mm}

We denote with $\mathsf{RCD}_K(M)$ the set of Riemannian metrics $g$ on $M$ such that $(M,g)$ is metrically complete with a convex boundary and $\mathrm{Ric}_g\geq K$. Then Theorem~\ref{thm:bounded-geometry-convex-boundary} and a simple scaling argument implies:

\begin{theorem} For all $K<0$ one has $\mathsf{RCD}_K(M)\neq \emptyset$.
\end{theorem}

It is stated in \cite{bxh} that for every Riemannian metric $g$ on $M$ and every number $K\in\R$ the following conditions are equivalent:
\begin{itemize}
    \item One has $g\in \mathsf{RCD}_K(M)$.
\item $(M,\dist_g, \mu_g)$ is an $\mathsf{RCD}(K,n)$ space,
\end{itemize}
justifying the notation $\mathsf{RCD}_K(M)$. We are, however, not going to use this result in the sequel. \vspace{1mm}

For later reference we also record the following result:

\begin{theorem}\label{main} Assume $g\in \mathsf{RCD}_K(M)$ for some $K\in\R$. \\
\emph{1)} One has the \emph{Bakry-Emery estimate} 
$$
|\dd P_t^g f|^2\leq \mathrm{e}^{-Kt} P_t^g(|\dd f|^2),\quad f\in  \mathsf{H}^{1,2}(\mathring{M},g|_{\mathring{M}})\cap \mathsf{L}^{\infty}(M,g).
$$
\emph{2)} For all $\epsilon>0$ one has the \emph{$\mathsf{L}^2$-Calderon-Zygmund inequality}, 
\begin{align*}
\left\|\mathrm{Hess}_g (f)\right\|_{\mathsf{L}^2}^2\leq \big(1+K^2/(2\epsilon^2)\big)\left\|\Delta_g f\right\|_{\mathsf{L}^2}+(K\epsilon^2)/2\left\|f\right\|_{\mathsf{L}^2}^2<\infty,\quad f\in\dom(H_g)\cap \mathsf{C}^\infty(M).
\end{align*}
\end{theorem}

The proof of part 1) of Theorem \ref{main} relies on the following simple auxiliary result, which holds without any assumptions on the geometry of $M$, and which allows to extend the validity of the asserted estimate for a sufficiently rich class of $f$'s:

\begin{lemma}\label{amb} For every Riemannian metric $g$ on $M$ and every $t> 0$, $f\in  \mathsf{H}^{1,2}(\mathring{M},g|_{\mathring{M}})$, $\psi\in  \mathsf{L}^{\infty}(M,g)$, the map
$$
[0,t)\longrightarrow \R, \quad s\longmapsto \int  |\dd P^g_{t-s} f|^2_g \; P^g_s \psi   \;\dd\mu_g 
$$
is continuous.
\end{lemma}

\begin{proof} Skipping $g$ in the notation, the map 
\begin{align}\label{aaa}
[0,\infty)\ni s\longmapsto P_s\in \ILL( \mathsf{L}^q(M))
\end{align}
is strongly continuous for all $q\in [1,\infty)$ and weak-*-continuous for $q=\infty$.
\vspace{1mm}
In addition, the map
\begin{align}\label{bbbb}
[0,\infty)\ni s\longmapsto \dd P_sf \in \Omega_{\mathsf{L}^2}(M)
\end{align}
is continuous (by Lemma 1.3.3 in \cite{fuku} and the semigroup property). It follows that 
$$
[0,t)\longrightarrow \mathsf{L}^1(M),\quad s\longmapsto |\dd P_{t-s} f|^2
$$
is continuous, and that 
$$
[0,t)\longrightarrow \mathsf{L}^{\infty}(M),\quad s\longmapsto P_s\psi
$$
is weakly continuous. Thus, if $s_n\to s$, then we can write
\begin{align*}
&\int  |\dd P_{t-s_n} f|^2 P_{s_n} \psi\;\dd\mu - \int  |\dd P_{t-s} f|^2  P_s\psi \;\dd\mu\\
&=\int  |\dd P_{t-s} f|^2 (P_{s_n} \psi-P_s\psi) \;\dd\mu- \int ( |\dd P_{t-s_n} f|^2-|\dd P_{t-s} f|^2)  P_s\psi \;\dd\mu,
\end{align*}
to conclude that this expression goes to $0$ as $n\to\infty$.
\end{proof}

\begin{proof}[Proof of Theorem \ref{main}] Again, we omit $g$ in the notation.\\
1)  In \cite{wang} (Corollary 3.2.6) this estimate has been shown for bounded $\mathsf{C}^1$-functions on $M$. It remains to extend the estimate to all $f\in  \mathsf{H}^{1,2}(\mathring{M},g|_{\mathring{M}})\cap \mathsf{L}^{\infty}(M,g)$: as for all $0<s<t$ we have $P_{s}f\in \mathsf{C}^1_b(M)$ by the smoothing property and the contraction property of $P_s$, given $0\leq \psi\in \mathsf{C}^{\infty}_c(\mathring{M})$ we get 
\begin{align*}
&\int |\dd P_t f|^2\psi \;\dd\mu= \int |\dd P_{t-s} P_{s} f|^2\psi \;\dd\mu\leq  \int \mathrm{e}^{-Kt} P_t(|\dd P_{t-s} f|^2)\psi\;\dd\mu\\
&= \mathrm{e}^{-Kt}\int  |\dd P_{t-s} f|^2 P_s \psi \;\dd\mu\to  \mathrm{e}^{-Kt}\int  |\dd P_{t} f|^2 \psi\;\dd\mu\quad\text{as $s\to 0+$},
\end{align*}
by Lemma \ref{amb}.\vspace{1mm}

2) We record the usual Bochner inequality 
$$
 |\mathrm{Hess}(\psi)|^2\leq (1/2)\Delta |\dd \psi|^2- (\dd \Delta \psi,\dd \psi)+K^2 |\dd \psi|^2 \quad\text{for all $\psi\in \mathsf{C}^{\infty}(M)$.}
$$
If in addition $N\psi=0$, then integrating by parts using Green's formula, we get
\begin{align*}
&\int(\dd \Delta \psi,\dd \psi)\dd\mu = \int  (\Delta \psi)^2 \dd\mu,\\
&\int|\dd \psi|^2\dd\mu=\int(\Delta \psi) \psi\dd\mu\leq 2\int(\Delta \psi)^2\dd\mu +2\int (\psi)^2\dd\mu,\\
&\int \Delta |\dd \psi|^2\dd\mu=0,
\end{align*}
which using $ab\leq \frac{a^2}{2\epsilon^2}+\frac{\epsilon^2b^2}{2}$ gives 
\begin{align}\label{aspsww}
\left\|\mathrm{Hess}(\psi)\right\|_{L^2}^2\leq (1+\frac{K^2}{2\epsilon^2})\left\|\Delta f\right\|_{L^2}+\frac{K\epsilon^2}{2}\left\|\psi\right\|_{L^22}^2
\end{align}
Since $\Delta$ with domain of definition $\{\psi\in \mathsf{C}^{\infty}_c(M):N\psi=0\}$ is essentially self-adjoint in $\mathsf{L}^2(M)$, given any smooth $f\in\dom(H)$ we can pick a sequence $f_j$ in $\{\psi\in \mathsf{C}^{\infty}_c(M):N\psi=0\}$ with 
$$
\left\|\Delta (f-f_j)\right\|_2^2 +\left\| f-f_j\right\|_2^2\to 0,
$$
showing that (\ref{aspsww}) holds with $\psi$ replaced by $f$. 
\end{proof}

Recall that $M$ is called \emph{orientable}, if there exists a nowhere vanishing smooth form $\widetilde{\omega}\in \Omega^n_{\mathsf{C}^{\infty}}  (M)$. In case this topological condition is satisfied, we shall call $M$ also \emph{smoothly orientable}, and refer to an $\widetilde{\omega}$ as above as a \emph{smooth orientation}. Any two smooth orientations are either equal or differ by a sign, and smooth orientations of course restrict to open subsets.

\begin{lemma}\label{asppqo} $M$ is smoothly orientable, if and only if $\mathring{M}$ is smoothly orientable, and an orientation on each of these uniquely determines an orientation on the other one.
\end{lemma}

\begin{proof}
		Let $M$ be smoothly orientable and let $\widetilde{\omega}\in \Omega^n_{\mathsf{C}^{\infty}} (M)$ be nowhere vanishing. Then $\iota^*\widetilde{\omega}$, where $\iota:\mathring{M}\to M$ is the canonical embedding, is a nowhere vanishing smooth $n$-form on $\mathring{M}$.\\
		Conversely, suppose $\mathring{M}$ is orientable and let $\widetilde{\omega}\in \Omega^n_{\mathsf{C}^{\infty}}(\mathring{M})$ be nowhere vanishing. After choosing an arbitrary Riemannian metric $g$ on $M$, we can assume that $|\widetilde{\omega}|_g=1$, i.e., if $(e_1,\ldots,e_n)$ is an $g$-orthonormal frame in an open neighbourhood $U$ of some boundary point of $M$, we have w.l.o.g.\ $\widetilde{\omega}(e_1,\ldots,e_n)(p)=1$ for all $p\in U\cap\mathring{M}$. Clearly, there exists a unique nowhere vanishing $\widehat{\omega}\in \Omega^m_{\mathsf{C}^{\infty}} (M)$ with $\iota^*\widehat{\omega} = \widetilde{\omega}$.
\end{proof}

Let
$$
\mathsf{Test}(M,g):= \{f\in \dom(H_g)\cap \mathsf{Lip}_b(M,g)\,|\,H_g f\in\dom(H_g) \}\subset \dom(H_g)\subset  \mathsf{H}^{1,2}(\mathring{M},g|_{\mathring{M}}).
$$
The following definition has been suggested in \cite{honda} in the context of Ricci limit spaces:

\begin{definition} One says that $(M,g)$ is \emph{mms-orientable}, if there exists $\omega\in \Omega^n_{\mathsf{L}^{\infty}}(M,g)$ with $|\omega|_g=1$ a.e. such that for all $f_1,\dots, f_n\in \mathsf{Test}(M,g)$ one has 
$$
 (\omega,\dd f_1\wedge\cdots\wedge  \dd f_n)\in \mathsf{H}^{1,2}(\mathring{M},g|_{\mathring{M}}),
$$
and we call such an $\omega$ an \emph{mms-orientation} on $(M,g)$.
\end{definition}

\begin{lemma}\label{sssaw} Let $g$ be an arbitrary Riemannian metric on $M$.\\
\textrm{1)} Let $U\subseteq M$ be open. Then every mms-orientation on $(M,g)$ restricts to an mms-orientation on $(U,g|_U)$.\\
\textrm{2)} Every mms-orientation on $(M,g)$ is an mms-orientation on $(\mathring{M},g|_{\mathring{M}})$, and vice versa. 
\end{lemma}

\begin{proof}1). This follows from the fact that 
$$
\dom(H_{g})|_U\subseteq \dom(H_{g|_U}),\quad \mathsf{Lip}_b(M,g)|_U\subseteq \mathsf{Lip}_b(U,g|_U)\quad\text{(because $(\dist_g)|_U\leq \dist_{g|_U})$},
$$
and so $\mathsf{Test}(M,g)|_{U}\subseteq\mathsf{Test}(U,g|_U)$.\\
2) The first statement follows from part 1). For the other direction, note that, since $\partial M$ is smooth and has measure zero, the underlying Sobolev spaces are equal, the underlying domains of the Neumann Laplacians are equal, and that $(\dist_g)|_{\mathring{M}}= \dist_{g|_{\mathring{M}}}$, so that each Lipschitz function on the interior of $M$ extends to a Lipschitz function on the whole of $M$. This shows that the test classes are equal, too.
\end{proof}

\begin{theorem}\label{main2} $M$ is smoothly orientable, if and only if there exists a Riemannian metric $g$ on $M$ such that $(M,g)$ is mms orientable. More precisely, if $\widetilde{\omega}$ is a smooth orientation on $M$, then for every $K\in\R$ and every metric $g\in \mathsf{RCD}_K(M)$ the form $\omega:=\widetilde{\omega}/|\widetilde{\omega}|_g$ is an mms-orientation on $(M,g)$. Conversely, given any metric $g$ on $M$ and any mms orientation $\omega$ of $(M,g)$, the form $\omega$ is automatically smooth. 
\end{theorem}

\begin{proof} $\Rightarrow$: Let $f_1,\dots, f_n\in \mathsf{Test}(M,g)$. Following the proof of Proposition 6.6 in \cite{honda}, we remark that since $H_g P^g_t f_i  = P^g_tH_g f_i$ converges in $\mathsf{L}^2(M,g)$ (to $H_g f_i$) and $P^g_t f_i$ converges in $ \mathsf{H}^{1,2}(\mathring{M},g|_{\mathring{M}})$ and thus in $\mathsf{L}^2(M,g)$ (to $f_i$) as $t\to 0+$, the $\mathsf{L}^2$-Calderon-Zygmund inequality implies 
\begin{align}\label{wsxcy}
\sup_{0<t<1} \left\|\mathrm{Hess}_g(P_t^gf_i)\right\|_{\mathsf{L}^2}<\infty.
\end{align}
Moreover, by the Bakry-Emery estimate and the Markovian property of $P^g_t$ we have 
\begin{align}\label{wsxcy2}
\sup_{0<t<1} \left\|\dd P^g_t f_i\right\|_{\mathsf{L}^\infty}<\infty.
\end{align}
Now we can conclude, using that $\nabla_g$ is a metric connection, that $\omega$ is $\nabla_g$-parallel, and using the derivation property of $\nabla_g$, the bound
\begin{align*}
&\left|\dd \big( \omega, \dd (P_t^g f_1)\wedge \cdots \wedge \dd (P^g_t f_n) \big)\right|\leq |\nabla_g \omega|\prod_i |\dd P_t^g f_i|+\sum_i |\mathrm{Hess}_g( P_t f_i)| \prod_{j\neq i} |\dd P^g_t f_j| \\
&= \sum_i |\mathrm{Hess}_g( P^g_t f_i)| \prod_{j\neq i} |\dd P_t f_j|, 
\end{align*}
and so by (\ref{wsxcy}) and (\ref{wsxcy2}),
$$
\sup_{0<t<1}\left\| \dd \big( \omega, \dd (P^g_t f_1)\wedge \cdots \wedge \dd (P_t^g f_n) \big)\right\|_{\mathsf{L}^2}<\infty.
$$
Moreover 
$$
\sup_{0<t<1}\left\|  \big( \omega, \dd(P_t^g f_1)\wedge \cdots \wedge \dd (P^g_t f_n) \big)\right\|_{\mathsf{L}^2}<\infty,
$$
by the contraction property of $P^g_t$, and so
$$
\sup_{0<t<1}\left\|  \big( \omega,\dd  (P^g_t f_1)\wedge \cdots \wedge \dd (P_t^g f_n) \big)\right\|_{\mathsf{H}^{1,2}}<\infty.
$$
Since 
$$
\big( \omega, \dd (P_t^g f_1)\wedge \cdots \wedge \dd (P_t^g f_n) \big)\to ( \omega, \dd  f_1\wedge \cdots \wedge \dd f_n )
$$
in $\mathsf{L}^2(M,g)$ as $t\to+$ (as $P^g_t f_i$ converges in $ \mathsf{H}^{1,2}(\mathring{M},g|_{\mathring{M}})$ to $f_i$), this shows that 
$$
( \omega, \dd  f_1\wedge \cdots \wedge  \dd f_n )\in  \mathsf{H}^{1,2}(\mathring{M},g|_{\mathring{M}}),
$$
completing the proof.\vspace{1mm}

$\Leftarrow$: In view of Theorem \ref{sssaw} and Lemma \ref{asppqo} we can assume $\partial M=\emptyset$. Pick an arbitrary connected relatively compact chart $((x^1,\dots,x^n),U)$ on $M$ such that $g|_U$ is quasi-isometric to the restriction $g^{\mathrm{Eucl}}_U$ of the Euclidean metric on $\IR^n$ to on $U$. Then we have  
$$
\omega|_U= \frac{\phi \dd x^1\wedge \cdots \wedge \dd x^n}{|\dd x^1\wedge \cdots \wedge \dd x^n|}
$$  
for some $\phi\in \mathsf{H}^{1,2}(U,g|_{U})=\mathsf{H}^{1,2}(U,g^{\mathrm{Eucl}}_U)$ with $|\phi|=1$ a.e., but the Euclidean Poincare inequality and Kato's inequality imply
$$
\int_{U} |\phi-\phi_{U}|^2 \; \dd x \lesssim  \int_{U} |\dd \phi|\dd x\leq   \int_{U} |\dd |\phi||\dd x=0,
$$
so that $\phi$ is constant a.e. on $U$ and $\omega$ can be chosen smooth on $U$.
\end{proof}

\begin{remark} If $\partial M\neq \emptyset$, then, even given the highly nontrivial fact $\mathsf{RCD}(M)\neq  \emptyset$, the '$\Rightarrow$' direction of Theorem \ref{main2} does not follow from Proposition 6.12 in \cite{honda}, as the Neumann-Sobolev capacity $\mathrm{Cap}_g(\partial M)$ of $\partial M$ is strictly positive for every metric Riemannian metric $g$ on $M$. While the fact $\mathrm{Cap}_g(\partial M)>0$ is clearly well-known to the experts, we have not been able to find a precise reference and so decided to include a simple prove for the sake of completeness: we can assume that $M$ is compact, and we recall from \cite{KM} that
\[
\mathrm{Cap}_g(\partial M) = \inf \left\{ \| \varphi \|_{\mathsf{H}^{1,2}}^{2} \,|\, \varphi \in \A(M,g) \right\},
\]
where
\[
\A(M,g) = \left\{ \varphi \in   \mathsf{H}^{1,2}(\mathring{M},g|_{\mathring{M}})\,|\, 0 \leq \varphi \leq 1, \, \varphi = 1 \text{ in a neighbourhood of } \partial M \right\}.
\]
Thus, with
\[
\A'(M,g): =\left\{ \varphi\in \mathsf{H}^{1,2}(\mathring{M},g|_{\mathring{M}})\,|\, 0 \leq \varphi \leq 1, \, \varphi = 1 \text{ on } \partial M \right\},
\]
where the boundary datum is understood in the sense of the trace map $ \mathsf{H}^{1,2}(\mathring{M},g|_{\mathring{M}})\to L^2(\partial M,g)$, we have
\[
\mathrm{Cap}_g(\partial M) \geq \inf \left\{  \| \varphi \|_{\mathsf{H}^{1,2}}\,|\, \varphi \in \A'(M,g) \right\}.
\]
Since $\A'(M,g)$ is a closed convex subset (hence weakly closed) of the reflexive Banach space $ \mathsf{H}^{1,2}(\mathring{M},g|_{\mathring{M}})$ and the Sobolev norm $\| \cdot \|_{ \mathsf{H}^{1,2}(\mathring{M},g|_{\mathring{M}})}$ is coercive and (weakly) lower semicontinuous in $\A'(M,g)$, we deduce that there exists $u \in \A'(M,g)$ such that
\[
\inf \{ \| \varphi \|_{\mathsf{H}^{1,2}} \,|\, \varphi \in \A'(M,g) \} = \| u \|_{ \mathsf{H}^{1,2}(\mathring{M},g|_{\mathring{M}})}.
\]
Since the trace map is bounded, we have $ \| u \|_{\mathsf{H}^{1,2}} >0$, so that
\[
\mathrm{Cap}_g(\partial M)  \geq   \| u \|_{\mathsf{H}^{1,2}}^{2} >0,
\]
as claimed.
\end{remark}

\end{document}